%% file: main.tex
\newif\ifMOR
\MORfalse      % arXiv/public preprint with authors
\ifMOR
  \documentclass[moor,dblanonrev]{informs4}
  \OneAndAHalfSpacedXI
\else
  \documentclass[11pt]{article}
\fi

\usepackage[T1]{fontenc}
\usepackage[utf8]{inputenc}
\usepackage{amsfonts,amsmath,amssymb,mathtools}
\ifMOR\else\usepackage{amsthm}\fi
\usepackage{algorithm,algpseudocode}
\usepackage{bm}
\usepackage{xcolor}
\usepackage{graphicx}
\usepackage{float}
\usepackage{nicematrix}
\usepackage{optidef}
\usepackage{mathabx}
\usepackage{pifont}
\usepackage{multirow}
\usepackage{listings}
\usepackage{natbib}

\ifMOR
  \NatBibNumeric
  \def\bibfont{\small}%
  \def\bibsep{\smallskipamount}%
  \bibpunct[, ]{[}{]}{,}{n}{}{,}%
\else
  \usepackage[margin=1in]{geometry}
  \usepackage{microtype}
\fi

\usepackage[colorlinks=true,breaklinks=true,bookmarks=true,urlcolor=blue,
  citecolor=blue,linkcolor=blue,bookmarksopen=false,draft=false]{hyperref}

\input{macros.tex}

\ifMOR
  \let\INFORMSproof\proof
  \let\INFORMSendproof\endproof
  \renewenvironment{proof}{\INFORMSproof{Proof.}}{\INFORMSendproof}

  \TheoremsNumberedThrough
  \EquationsNumberedThrough

  \makeatletter
  \@ifundefined{assumption}{\newtheorem{assumption}{Assumption}}{}
  \@ifundefined{definition}{\newtheorem{definition}{Definition}}{}
  \@ifundefined{remark}{\newtheorem{remark}{Remark}}{}
  \makeatother
\else
  \newtheorem{theorem}{Theorem}
  \newtheorem{lemma}[theorem]{Lemma}
  \newtheorem{proposition}[theorem]{Proposition}
  
  \newtheorem{assumption}[theorem]{Assumption}
  \theoremstyle{definition}
  \newtheorem{definition}[theorem]{Definition}
  \theoremstyle{remark}
  \newtheorem{remark}[theorem]{Remark}
  \numberwithin{equation}{section}
\fi

\newcommand{\compactdisplayspacing}{%
  \setlength{\abovedisplayskip}{0pt}%
  \setlength{\belowdisplayskip}{0pt}%
  \setlength{\abovedisplayshortskip}{0pt}%
  \setlength{\belowdisplayshortskip}{0pt}%
}

\begin{document}
\allowdisplaybreaks

\ifMOR
  % ---------------- MOOR anonymous metadata ----------------
  \RUNAUTHOR{Anonymous}
  \RUNTITLE{Variance-Reduced ConEx for Function Constraints}
  \TITLE{Variance-reduced First-order Methods for finite-sum function
  constrained convex optimization}

  \ARTICLEAUTHORS{%
    \AUTHOR{}
    \AFF{}
  }

  \ABSTRACT{%
    \input{abstract_moor.tex}
  }%

  \MSCCLASS{Primary: 90C25; secondary: 90C15, 90C30}
  \ORMSCLASS{Primary: Programming: nonlinear; secondary: Analysis of algorithms}
  %\HISTORY{}

  \maketitle
\else
% ---------------- arXiv/public metadata ----------------
\title{Variance-reduced First-order Methods for Finite-sum Function
Constrained Convex Optimization}

\author{%
  Yasin Gharehmohammadlou%
  \thanks{%
    \texttt{ygarehmohammadlou@smu.edu},
    Operations Research and Engineering Management,
    Southern Methodist University.
  }
  \qquad
  Digvijay Boob%
  \thanks{%
    \texttt{dboob@smu.edu},
    Operations Research and Engineering Management,
    Southern Methodist University.
  }
}

\date{}
\maketitle

  \begin{abstract}
    \input{abstract_moor.tex}
  \end{abstract}


\fi

\compactdisplayspacing
\input{sec-intro}
\input{sec-notation}
\input{sec-simple-algorithm}

\input{sec-full-algorithm}
\input{sec-conclude}

\ifMOR
  \bibliographystyle{informs2014}
\else
  \bibliographystyle{unsrtnat}
\fi
\bibliography{ref}

\end{document}

%% file: macros.tex
\global\long\def\vertiii#1{\left\vert \kern-0.25ex  \left\vert \kern-0.25ex  \left\vert #1\right\vert \kern-0.25ex  \right\vert \kern-0.25ex  \right\vert }%
\global\long\def\wt#1{\widetilde{#1}}%
\global\long\def\wh#1{\widehat{#1}}%
\global\long\def\wb#1{\widebar{#1}}%
\global\long\def\zero{\mathbf{0}}%

\global\long\def\norm#1{\lVert#1\rVert}%
\global\long\def\grad{\nabla}%
\global\long\def\argmin{\operatornamewithlimits{argmin}}%
\global\long\def\and{\text{and}}%
\global\long\def\Ebb{\mathbb{E}}%
\global\long\def\Rbb{\mathbb{R}}%
\DeclarePairedDelimiterX{\inprod}[2]{\langle}{\rangle}{#1, #2}
\DeclarePairedDelimiter\abs{\lvert}{\rvert}

\DeclarePairedDelimiter{\paran}{(}{)}

\DeclarePairedDelimiterX{\gnorm}[3]{\lVert}{\rVert_{#2}^{#3}}{#1}

\global\long\def\tsum{{\textstyle {\sum}}}
\newcommand{\vrconex}{\texttt{VaRAn-ConEx}}

\newcommand{\MeanFVR}
{\ensuremath{0.582591302172}}

\newcommand{\MeanGVR}
{\ensuremath{2.808908900035 \times 10^{-4}}}

\newcommand{\MeanViolVR}
{\ensuremath{2.808908900035 \times 10^{-4}}}

\newcommand{\MeanLambdaVR}
{\ensuremath{4.155604683294}}

\newcommand{\MeanFConEx}
{\ensuremath{0.578537872595}}

\newcommand{\MeanGConEx}
{\ensuremath{1.109530158969 \times 10^{-3}}}

\newcommand{\MeanViolConEx}
{\ensuremath{1.109530158969 \times 10^{-3}}}

\newcommand{\MeanYConEx}
{\ensuremath{4.110441332419}}

\newcommand{\MeanFYu}
{\ensuremath{0.583265437793}}

\newcommand{\MeanGYu}
{\ensuremath{2.323955180259 \times 10^{-3}}}

\newcommand{\MeanViolYu}
{\ensuremath{2.323955180259 \times 10^{-3}}}

\newcommand{\MeanQYu}
{\ensuremath{1700.883282583686}}

\newcommand{\MeanQOverVYu}
{\ensuremath{4.140238346771}}

%% file: abstract_moor.tex
This paper studies finite-sum convex optimization with convex function constraints, where both the objective $f(x)=\frac{1}{m}\sum_{i=1}^m f_i(x)$ and the constraint mapping $g(x)=\frac{1}{n}\sum_{j=1}^n g_j(x)$ may have a large number of components. We propose a variance-reduced anchored constraint extrapolation method (\vrconex), which combines stochastic component sampling with primal--dual extrapolation for nonlinear function constraints. Relative to stochastic approximation schemes such as ConEx \citep{boob2023stochastic}, \vrconex~method reduces gradient complexity by using epochwise full-gradient correction and variance-reduced estimators for both the objective and constraint components. In the semi-finite case where $n=1$, the method attains gradient complexities $O(\sqrt m/\epsilon)$ under constant steps and $O(m/\epsilon^{1/2}+1/\epsilon^{3/2})$ under nonconstant steps. In the fully finite-sum case, we establish boundedness of the Lagrange multiplier sequence and derive convergence and gradient-complexity bounds for constant and nonconstant step-size policies. Under these policies, we then prove similar gradient complexity with $m$ replaced by $\max\{m,n\}$. We also show acceleration when the finite-sum objective $f$ is strongly convex, irrespective of components $f_i$.

%% file: sec-intro.tex
\section{Introduction.}
In this paper, we consider the convex optimization problem with finite-sum structure as follows:
\begin{align}
    &\min_{x \in X} f(x) := \tfrac{1}{m}\tsum_{i=1}^{m}f_i(x) \nonumber\\
    &\text{s.t. } g(x) := \tfrac{1}{n}\tsum_{j=1}^{n}g_j(x) \le 0,\label{eq:c1}
\end{align}
where $f_i:\Rbb^d \to \Rbb$ and $g_j:\Rbb^d \to \Rbb^N$ are convex mappings with respect to $x$. The set $X \subset \Rbb^d$ is a compact convex set with a readily available projection oracle. Note that $g(x) \le 0$ models $N$ convex function constraints. For the sake of simplicity, we assume number of components in each function constraint is equal to $n$. Though our results easily extend to function constraints with differing number of components. %The algorithms in this paper have been designed in a way that the problem instances could have many finite-sum convex constraints, but for the sake of simplicity, we consider one constraints with finite-sum convex structure and one could easily realize that the algorithms could be used for many such constraints as well. 

Optimization problem with function constraints of type \eqref{eq:c1} arise in many modern applications including machine learning, signal processing, network optimization, resource allocation, online decision making, and statistical learning. First-order methods are widely used to solve them. These algorithms use the gradient of the objective function and constraints in each iteration. However, problem \eqref{eq:c1} has many function components $f_i(x)$ and $g_j(x)$. Thus, gradients of $f$ or $g$ require evaluating all components, which is computationally ineffective. To address this, variance-reduction methods have been explored in the literature. Instead of evaluating the full gradient, they take a random sample from the finite sum as a stochastic, unbiased estimator of the true gradient. By systematically updating the stochastic gradient, they show that the variance of the estimator reduces over the algorithmic run and hence achieve better computational complexity.  Our main goal here is to design such methods the for problem \eqref{eq:c1}.

The development of stochastic first-order methods for large-scale finite-sum optimization has attracted significant attention during the past decade. The vast majority of works address problem \eqref{eq:c1} where function constraint $g(x) \le 0$ is absent. One of the earliest breakthroughs for this problem was the stochastic variance reduced gradient (SVRG) method proposed by Johnson and Zhang \cite{johnson2013accelerating}. SVRG utilizes periodically updated reference gradients to construct low-variance stochastic estimators of $\grad f$ while maintaining inexpensive iterations. Shortly thereafter, Schmidt et al. \cite{schmidt2017minimizing} introduced the stochastic average gradient (SAG) method, while Defazio et al. \cite{defazio2014saga} proposed SAGA, which stores historical gradient information to achieve variance reduction without periodic full-gradient computations. Xiao et al. \cite{xiao2014proximal} further extended these ideas through the Prox-SVRG framework for composite optimization problems.

Subsequent developments focused on improving theoretical complexity guarantees and practical performance. Nguyen et al. \cite{nguyen2017sarah} introduced SARAH, which employs recursive gradient estimators and achieves improved theoretical properties compared to classical SVRG-type methods. Further advances include SPIDER \cite{fang2018spider} and SNVRG \cite{zhou2020stochastic}, which significantly improved gradient complexity for large-scale stochastic optimization, particularly in nonconvex settings. Allen-Zhu \cite{allen2017katyusha} proposed the Katyusha algorithm, which combines variance reduction and Nesterov-type acceleration to achieve optimal complexity guarantees for strongly convex finite-sum optimization. Along a related direction, Lan  et al. \cite{lan2019unified} developed the Variance-Reduced Accelerated Gradient (Varag) method and established unified accelerated complexity guarantees for smooth convex finite-sum optimization problems. Lan et al. \cite{lan2018optimal} further developed optimal randomized incremental gradient methods for structured finite-sum optimization and advanced the complexity theory of randomized first-order methods. Lan et al. \cite{lan2019accelerated} developed accelerated stochastic algorithms for nonconvex finite-sum and multi-block optimization, introducing the RapGrad and RapDual methods with improved complexity guarantees. Their work demonstrated that randomized acceleration can substantially reduce gradient and subproblem complexity for large-scale optimization problems. A comprehensive survey of variance reduction methods and their complexity guarantees is provided by Gower et al. \cite{gower2020variance}. These developments established variance reduction as one of the central tools in modern large-scale optimization.

While variance reduction methods were initially developed for finite-sum minimization problems, recent years have witnessed substantial efforts to extend these techniques to convex-concave saddle-point problems and monotone variational inequalities. One of the earliest works in this direction was the stochastic variance-reduced saddle-point framework of Palaniappan et al. \cite{palaniappan2016stochastic}, which extended SVRG- and SAGA-type estimators to convex-concave saddle-point optimization and established linear convergence guarantees under strong convexity-concavity assumptions. Yazdandoost Hamedani et al. \cite{hamedani2021stochastic} proposed a stochastic variance-reduced accelerated primal-dual method for convex-concave saddle-point optimization, while related doubly stochastic extensions were studied in \cite{jalilzadeh2021doubly}. More recently, Wang et al. \cite{wang2025near} combined stochastic extrapolation and variance reduction for finite-sum coupled compositional optimization. These developments further illustrate the growing role of variance reduction in saddle-point optimization, variational inequalities, and primal-dual algorithms.

The variance reduction technique is the subject of active research in the broader problem of monotone variational inequalities. Classical methods such as the extragradient method of Korpelevich \cite{korpelevich1976extragradient} and the Mirror-Prox method of Nemirovski \cite{nemirovski2004prox} provide powerful tools for solving monotone variational inequalities and saddle-point problems. Their stochastic counterparts were later developed by Juditsky et al. \cite{juditsky2011solving}. Building upon these foundations, Alacaoglu et al. \cite{alacaoglu2022variance} developed a general variance-reduction framework for finite-sum monotone variational inequalities and demonstrated that variance-reduced extragradient-type methods can substantially improve oracle complexity compared with classical stochastic VI methods.

%Variance-reduced primal-dual methods have also been developed for finite-sum saddle-point problems. 

Despite the success of variance reduction methods for finite-sum optimization and saddle-point problems, extending these techniques to function-constrained optimization remains significantly more challenging. Direct projections onto nonlinear function constraints are often computationally expensive or intractable, particularly when constraint functions possess finite-sum structure. Moreover, many variance-reduced saddle-point methods rely on assumptions such as bounded feasible sets or existence of uniform Lipschitz constant for coupling functions. These assumptions do not hold for the Lagrangian min-max formulation associated with function-constrained optimization \eqref{eq:c1}. The key difficulty is that the Lagrange reformulation $f(x) + \lambda^\top g(x)$ is not uniformly Lipschitz  smooth in $x$ for all $\lambda \ge 0$. In absence of global Lipschitz smoothness constant, the techniques used to show variance reduction for both optimization and saddle-point problems fail.

A closely related problem of \eqref{eq:c1} is the problem where both objective $f(x)$ and constraint functions $g(x)$ have an expectation form. Mahdavi et al. \cite{mahdavi2012trading} studied online convex optimization with long-term constraints and established regret-feasibility tradeoffs. Neely \cite{neely2010stochastic} proposed a drift-plus-penalty framework for this problem with applications in stochastic network optimization and queuing systems. Yu et al. \cite{yu2017online} further investigated online stochastic optimization with expectation constraints using primal-dual techniques. Lan et al. \cite{lan2020function} studied stochastic optimization with function and expectation constraints and developed first-order algorithms with provable complexity guarantees. These works established important foundations for stochastic optimization under nonlinear functional constraints.

A major advance in stochastic function-constrained optimization was the Constraint Extrapolation (ConEx) framework proposed by Boob et al. \cite{boob2023stochastic}. ConEx combines primal-dual updates with extrapolation of constraint information and establishes convergence guarantees for stochastic convex optimization with function constraints without requiring projections onto the feasible set or bounded dual domains. The framework provides a unified treatment of deterministic, semi-stochastic (assumed deterministic $g$ function), and fully stochastic (where both $f$ and $g$ are stochastic) constrained optimization problems. %Subsequent developments include level-constrained first-order methods \cite{boob2024level}, 
These techniques are further extended for stochastic variational inequalities with stochastic function constraints \cite{boob2026first}. The augmented constraint extrapolation method by \cite{boob2024optimal} showed optimal last-iterate guarantees for semi-stochastic problem. \cite{boob2024level} considered nonconvex version and showed convergence rates for semi-stochastic problems. Related developments have also appeared in constrained stochastic compositional optimization. Yang et al. \cite{yang2025stochastic} and Thomdapu et al. \cite{thomdapu2023stochastic} studied optimization problems with nonlinear compositional constraints. However, their methods do not exploit finite-sum objective and constraint structures through variance-reduced primal-dual extrapolation. Among existing methods, ConEx is particularly relevant to the present work since it applies to fully stochastic problems. It applies stochastic primal-dual updates with extrapolation mechanisms for nonlinear functional constraints to handle the aforementioned non-existence of uniform Lipschitz constant due to unbounded dual set $\{\lambda \ge 0\}$.

More recently, Lu et al. \cite{lu2025first} developed accelerated and variance-reduced first-order methods for stochastic and finite-sum convex optimization with deterministic constraints. In essence, their problem is equivalent to \eqref{eq:c1} when $n = 1$. Their framework is based on penalty reformulations and provides complexity guarantees for constrained optimization problems where they can show that the stochastic solution obtained is deterministically feasible up to some error tolerance.

Despite these advances, relatively few methods simultaneously exploit finite-sum structure in both the objective and nonlinear constraint functions while combining variance reduction, primal-dual extrapolation, and explicit gradient complexity analysis. As discussed, most existing variance-reduced methods focus either on unconstrained optimization, saddle-point problems, or linearly constrained formulations. On the other hand, many stochastic function-constrained optimization methods do not exploit finite-sum structure and therefore cannot fully benefit from modern variance reduction techniques.

The present work aims to bridge this gap by introducing a variance-reduced anchored constraint extrapolation framework for convex optimization with finite-sum nonlinear function constraints. The proposed method combines variance-reduced estimators with primal-dual extrapolation and proximal regularization, explicitly exploiting finite-sum structure in both the objective and constraint functions. We establish improved gradient complexity guarantees for both semi-finite and fully finite settings and provide a unified analysis of variance-reduced constrained extrapolation methods. %In contrast to existing variance-reduced saddle-point and stochastic constrained optimization methods, the proposed framework explicitly exploits finite-sum structures in both the objective and constraint functions while establishing improved gradient complexity guarantees for semi-finite and fully finite stochastic constrained optimization settings.
%In this work, we address the fully finite structure for both the objective function and function constraints. 
We studied both convex and strongly-convex cases separately. In the both cases mentioned, we considered two general cases as semi-finite-sum case (finite-sum structure for the objective function only) and fully-finite-sum case (finite-sum structure for both the objective function and the function constraints). In case of semi-finite-sum for convex case, we achieved $O(\tfrac{m^{1/2}}{\epsilon})$ and $O(\tfrac{m}{\epsilon^{1/2}}+\tfrac{1}{\epsilon^{1.5}})$ gradient computations complexity for constant and non-constant step-size policy separately. For the non-constant policy for latter case, we needed to prove an upper-bound for the expected Lagrangian terms. For fully-finite-sum case, bounding the expected Lagrangian terms was the challenging part. We achieved acceleration after we proved a constant upper-bound for the expected Lagrangian terms with different parameters setups for each of them. After bounding the expected Lagrangian terms, the gradient computations complexity for constant and non-constant step-size policies will be $O(\tfrac{\max\{m,n\}^{1/2}}{\epsilon})$ and $O(\tfrac{\max\{m,n\}}{\epsilon^{1/2}}+\tfrac{1}{\epsilon^{1.5}})$. For strongly-convex case, both semi-finite-sum and fully-finite-sum cases have non-constant step-size policies. For semi-finite-sum case the gradient computations complexity is $O(\tfrac{m}{\epsilon^{1/3}} + \tfrac{1}{\epsilon^{2/3}})$. For Fully-finite-sum case in strongly-convex case, similar to convex case, we achieve the acceleration after we prove a bound on the expected Lagrangian terms and the step-size policies for bounding the expected Lagrangian terms and achieving acceleration are different. The gradient computations complexity is $O(\tfrac{\max\{m,n\}}{\epsilon^{1/3}} + \tfrac{1}{\epsilon^{2/3}})$. 
Table~\ref{tab:complexity} summarizes the complexity guarantees established in this paper and compares the semi-finite and fully finite settings under convex and strongly convex assumptions.

\begin{table}[htbp]
\centering
\caption{Gradient complexity comparison}
\label{tab:complexity}
\begin{tabular}{|c|c|cc|}
\hline
&
\multicolumn{1}{c|}{Strongly convex}
&
\multicolumn{2}{c|}{Convex}
\\
\cline{2-4}
&
Non-constant step-size
&
Constant step-size
&
Non-constant step-size
\\
\hline
Semi-finite
&
$O(\tfrac{m}{\epsilon^{1/3}} + \tfrac{1}{\epsilon^{2/3}})$
&
$O(\tfrac{m^{1/2}}{\epsilon})$
&
$O(\tfrac{m}{\epsilon^{1/2}}+\tfrac{1}{\epsilon^{1.5}})$
\\
\hline
Fully-finite
&
$O(\tfrac{\max\{m,n\}}{\epsilon^{1/3}} + \tfrac{1}{\epsilon^{2/3}})$
&
$O(\tfrac{\max\{m,n\}^{1/2}}{\epsilon})$
&
$O(\tfrac{\max\{m,n\}}{\epsilon^{1/2}}+\tfrac{1}{\epsilon^{1.5}})$
\\
\hline
\end{tabular}
\end{table}

\paragraph{Outline.} This paper is organized as follows. Section \ref{sec-notation}  provides the notation and preliminaries. Section \ref{sec-method} describes the algorithm and assumptions necessary for its convergence. Section \ref{sec-method2.2} is about the analysis of the method to construct the gap function for the convergence analysis. Sections~\ref{sec-method2.2.1} and \ref{sec-method2.2.2} elaborate the gap convergence of the algorithm for the semi-finite and fully-finite cases, respectively. In both cases, we present constant and nonconstant step-size policies. Section \ref{sec-strong} discusses the strongly convex case of Problem \eqref{eq:c1} as well as the algorithm to tackle it. Sections \ref{sec-semi} and \ref{sec-full} study the semi-finite and fully-finite form of Problem \eqref{eq:c1} with a strongly-convex objective function with non-constant step-size policies. In Section \ref{sec-con}, we perform numerical tests to verify our convergence results. %for method of discussion and its comparison to the two other method, and 
We conclude in Section \ref{sec-con}.% gives an overall conclusion for this study.    

%% file: sec-notation.tex
\subsection{Notation and preliminaries.}\label{sec-notation}
%Throughout the paper, we use the following notations. 
We denote Euclidean norm as $\norm{\cdot}$. For any positive integer $l$, we denote the set $\{1,\dots, l\}$ as $[l]$. %All expectation signs in this work are with respect to random variables $x_t^k$ and $\lambda_t^k$ for iteration $t$ which are associated with sampling random variables $i_t^k$, $j_t^k$ and $p_t^k$. The sample space for random variables $j_t^k$ and $p_t^k$ is $[n]$ and for $i_t^k$ is $[m]$. 
%For the sake of simplicity, we consider only one constraint, but the reader could notice that, this study could be applied for multiple function constraints.
Recall that $X \subset \Rbb^d$ is a compact convex set. We use $\Rbb^k_+$ to denote $k$-dimensional nonnegative vectors. We denote the Lagrange multipliers associated with function constraint $g(x) \le 0$ by variable $\lambda \in \Rbb^N_+$. 

We assume all component functions $f_i, i \in [m]$ and $g_j, j \in [n]$ are differentiable. We denote gradient of $f_i$ as $\grad f_i$ and Jacobian of $g_j$ as $\grad g_j$. We further assume that all component functions of both objective and function constraint mappings are Lipchitz smooth as well as Lipschitz continuous, i.e., for $f_i(x) , i \in [m]$ and $g_j(x) , j \in [n]$ we have 
\begin{align}
    \abs{f_i(x_1) - f_i(x_2)} &\leq M_f \norm{x_1 - x_2} &
    \norm{g_j(x_1) - g_j(x_2)} &\leq M_g \norm{x_1 - x_2}. \label{eq:p1}\\
%\end{align}
%We, also consider a $f(x)$ and $g(x)$ to be Lipchitz continuous
%\begin{align}
    \norm{\grad f_i(x_1) - \grad f_i(x_2)} &\leq L_f \norm{x_1 - x_2} & \norm{\grad g_j(x_1) - \grad g_j(x_2)} &\leq L_g \norm{x_1 - x_2}. \label{eq:p2}
\end{align}
The gradient norms of these functions are bounded, meaning
\begin{align}
    \norm{\grad f_i(\cdot)} \leq M_f && \norm{\grad g_j(\cdot)} \leq M_g. \label{eq:p3}
\end{align} 

%% file: sec-simple-algorithm.tex
\section{Variance-reduced anchored constraint extrapolation method.}
It is well known that for a compact set $X$, we get existence of an optimal solution $x^*$ for \eqref{eq:c1} immediately. We consider the following equivalent reformulation of \eqref{eq:c1}
\begin{align}
    \min_{x \in X}\max_{\lambda\geq 0} \Gamma(x,\lambda) := f(x) + \lambda^\top g(x), \label{eq:p14}
\end{align}
where, $x$ is the primal variable, and $\lambda$ is the Lagrangian dual variable. Throughout this paper, we assume existence of the optimal dual solution. We state this formally below.
\begin{assumption}
    We assume existence of Lagrange multiplier $\lambda^* \ge 0$ such that $(x^*,\lambda^*)$ is saddle point solution of \eqref{eq:p14} satisfying 
    \begin{align*}
    \Gamma(x^*,\lambda) \leq \Gamma(x^*,\lambda^*) \leq \Gamma(x,\lambda^*). 
\end{align*}
\end{assumption}
%According to KKT conditions, $(x^*,\lambda^*)$ is the saddle point solution of \eqref{eq:p14}, if
This motivates us to consider the following gap function as a convergence criterion.
\begin{definition}
    For any given $\epsilon > 0$, we say a solution pair $(\wb{x},\wb{\lambda})$ is an $\epsilon$-gap solution if
\begin{align*}
    \Gamma(\wb{x},\lambda^*) - \Gamma(x^*,\wb{\lambda}) \leq \epsilon.
\end{align*}
\end{definition}
 We use $\epsilon$-gap criterion to measure the progress towards convergence. Clearly, $\epsilon = 0$ iff we are at a saddle point of \eqref{eq:p14}. We are now ready to present the \vrconex~method.%In this study, we assume that Problem \eqref{eq:c1} has a primal-dual solution $(x^*,\lambda^*)$ and we can use $\epsilon$-gap criterion for the method of discussion .   
\subsection{Variance-reduced Anchored ConEx method.} \label{sec-method}
We present Algorithm \ref{alg:alg2} for the general convex case. As the name suggests, it includes a variance reduction part, which is the part for the sampling from the finite-sum objective and function constraint, and convex extrapolation, which uses a linear approximation of the constraint function and uses it in the proximal operator for the Lagrangian update. The linear approximation for $g_j(x_t)$ with respect to $x_{t-1}$ is
\begingroup
\setlength{\abovedisplayskip}{4pt}
\setlength{\abovedisplayshortskip}{4pt}
\setlength{\belowdisplayskip}{4pt}
\setlength{\belowdisplayshortskip}{4pt}
\begin{equation}
    \ell_{g_j}^{t-1}(x_t)
    =
    g_j(x_{t-1})
    +
    \inprod{\grad g_j(x_{t-1})}{x_t-x_{t-1}}
    \label{eq:p5}
\end{equation}
\endgroup
We use the linear approximation operator \eqref{eq:p5} as an extrapolation.

We can now describe what the algorithm does. In Algorithm \ref{alg:alg2}, we start with the initialization in line 1. Lines 4 and 5 calculate the constraint value and the gradient of the Lagrangian function with respect to $x$, which will be used for the variance reduction steps. Step 7 gives a sample from the constraint. Step 8 gives a convex combination of the previous Lagrangian iterate and the average of the previous epoch, which will be used for extrapolation. It acts as an anchor point in line 11. Step 9 is the variance reduction for sampling from the constraint. Step 10 gives a linear form of the sampled constraint for iteration $t$ that will be used for extrapolation. Step 11 is the prox update for $\lambda$, using extrapolation and the anchor point. Step 12 gives samples from both the objective and the constraint for variance reduction. Step 13, similar to step 8, gives a convex combination of the previous iteration and the average of the previous epoch, which will act as a prox center. Step 14 is the variance reduction for the gradient of the primal-dual function with respect to $x$. Step 15 is the prox operation over the feasible region $X$. At the end of each epoch, in step 17, we average over all iterations of the epoch and we will use it for the next epoch. Step 18 is the initialization for the next epoch for the first iteration.
\begin{algorithm}[H]
	\begin{algorithmic}[1]
		\State {\bf Input: }$\wt{x}^0 \in X, \wt{\lambda}^0 = \zero$. 
		\State Set $(x_0^1, \lambda_0^1) =(\wt{x}^0,\wt{\lambda}^0)$
		\For{$k=1,2,\dots,K$}
        \State $G_{\lambda} \gets g(\tilde{x}^{k-1})$
		\State $G_x \gets \grad f(\tilde{x}^{k-1}) + \grad g(\tilde{x}^{k-1}) \wt{\lambda}^{k-1}$
        %\State In epoch $k$
        \For{ $t = 0,1, \dots, T^k-1$}
        \State {Pick a random subset $j_t^k \in [n]$ for function constraint
        }
        \State $\wh{\lambda}^k_t \gets (1-\gamma^k_\lambda)\lambda^k_t + \gamma^k_\lambda \wt{\lambda}^{k-1}$
		\State $\xi_t^k \gets g_{j_t^k}(x^k_t) - g_{j_t^k}(\wt{x}^{k-1}) + G_\lambda$
		\State $q^k_t \gets \ell{g_{j_t^k}}(x^k_t) - g_{j_t^k}(x^k_{t-1})$
        \State $\lambda_{t+1}^k \gets \argmin_{\lambda \ge 0} -\inprod{\xi_t^k + q^k_t}{\lambda} + \frac{1}{2\sigma^k}\norm{\lambda - \wh{\lambda}^k_t}^2$
        \State {Pick random subsets $i_t^k \in [m]$ and $p_t^k \in [n]$ for objective function and function constraint}
        %\State{function and function constraint}
        \State $\wh{x}_t^k \gets (1-\gamma^k_x)x^k_t + \gamma^k_x \wt{x}^{k-1}$
        \State $\beta_t^k \gets \grad f_{i_t^k}(x^k_t)+ \lambda_{t+1}^k \grad g_{p_t^k}(x^k_t)  - [\grad f_{i_t^k}(\wt{x}^{k-1}) + \grad g_{p_t^k}(\wt{x}^{k-1})\wt{\lambda}^{k-1}] + G_x$
        \State $x_{t+1}^k \gets \argmin_{x \in X} \inprod{\beta_t^k}{x} + \frac{1}{2\tau^k} \norm{x-\wh{x}^k_t}^2$
        \EndFor
        \State {\bf Output: }$(\wt{x}^k, \wt{\lambda}^k) := \tfrac{1}{T^k}\tsum_{t =1}^{T^k}(x^k_t, \lambda^k_t)$
        \State $(x^{k+1}_0, \lambda^{k+1}_0, x^{k+1}_{-1}, \lambda^{k+1}_{-1}) = (x^k_{T^k},\lambda^k_{T^k}, x^k_{T^k-1}, \lambda^k_{T^k-1})$
		\EndFor
	\end{algorithmic}
	\caption{Variance-reduced anchored constraint extrapolation (\vrconex) method} \label{alg:alg2}
\end{algorithm}
First, we state useful propositions and inequalities that will be used frequently.
\begin{proposition}
    [Mirrored Descent Theorem]\label{theorem:1}
    Given $\wb{x} \in X \cap dom f$ and $s \in \mathcal{X}$ and $t \in \mathrm{R}^+$ and let $x^+ = \argmin_{x \in X} f(x) - \inprod{s}{x} + t\norm{x - \wb{x}}^2$. Then for all $x \in X$, the following inequality holds:
    \begin{align}
        f(x) \geq f(x^+) + \inprod{s}{x - x^+} + 2t\inprod{\wb{x} - x^+}{x - x^+} \label{eq:a3}
    \end{align}
\end{proposition}
    \begin{proof}
        Define $\varphi(x) := f(x) - \inprod{s}{x} + t\norm{x - \wb{x}}^2$. Then, $x^+ = \argmin_{x \in X} \varphi(x)$. Using  %and $x^+$ is the minimizer of this function then by first order optimality, 
        optimality of $x^+$ for convex $\varphi(x)$, we have $\inprod{\grad\varphi(x^+)}{x - x^+} \geq 0$. 
        Using $\grad\varphi(x^+) = \grad f(x^+) -s + 2t(x^+ - \wb{x})$ and rearranging terms, we have %Also $\grad\varphi(x^+) = \grad f(x^+) -s + 2t(x^+ - \wb{x})$ and the optimality condition will be 
        %$\inprod{\grad f(x^+) -s + 2t(x^+ - \wb{x})}{x - x^+} \geq 0$. This will give 
        $\inprod{\grad f(x^+)}{x - x^+} \geq \inprod{s + 2t(x^+ - \wb{x})}{x - x^+}$. Using convexity of $f$, we get  $f(x) \geq f(x^+) + \inprod{\grad f(x^+)}{x - x^+}$. Combining these relations, we obtain 
        %will be $f(x) \geq f(x^+) + \inprod{s + 2t(x^+ - \wb{x})}{x - x^+} \geq f(x^+) + \inprod{s}{x - x^+} + \inprod{2t(x^+ - \wb{x})}{x - x^+}$ which is the final 
        the concluding inequality \eqref{eq:a3}.
    \end{proof}

\begin{proposition}
    %[Euclidean distance inequality bounding]
    \label{theorem:2}
    Let $x^+$ be as defined in Proposition \ref{theorem:1}. Then, we have
    \begin{align*}
        f(x^+) -f(x) + \inprod{s}{x - \wb{x}} \leq t(\norm{x - \wb{x}}^2 - \norm{x - x^+}^2) + \tfrac{1}{2t}\norm{s}^2. 
    \end{align*} 
    \begin{proof}
        Using Arithmetic-Geometric-mean inequality, we have $\inprod{s}{x^+ - \wb{x}} \leq \frac{1}{2t}\norm{s}^2 + 2t\norm{x^+ - \wb{x}}^2$. We also have the fact $2t\inprod{\wb{x} - x^+}{x - x^+} = -t(\norm{\wb{x} - x^+}^2 + \norm{x - x^+}^2 - \norm{x - \wb{x}}^2)$. Summing these two relations and \eqref{eq:a3} from Proposition \ref{theorem:1}, and rearranging terms, we obtain the concluding relation.
        %From Proposition \ref{theorem:1}, adding $\inprod{s}{x - \wb{x}}$ to both sides of \eqref{eq:a3}, which gives $f(x^+) - f(x) + \inprod{s}{x - \wb{x}}\leq \inprod{s}{x^+ - \wb{x}} - 2t\inprod{\wb{x} - x^+}{x - x^+}$, and  inner product $\inprod{s}{x^+ - \wb{x}}$ could be upper bounded as $\inprod{s}{x^+ - \wb{x}} \leq \frac{1}{2t}\norm{s}^2 + 2t\norm{x^+ - \wb{x}}^2$ via Arithmetic-Geometric mean inequality. Inner product $2t\inprod{\wb{x} - x^+}{x - x^+}$ could be upper bounded as $2t\inprod{\wb{x} - x^+}{x - x^+} \leq -t(\norm{\wb{x} - x^+}^2 + \norm{x - x^+}^2 - \norm{x - \wb{x}}^2$ which comes from Pythagorean theorem for vectors. Using upper bounds in the original inequality, the result is $f(x^+) - f(x) + \inprod{s}{x - \wb{x}}\leq \frac{1}{2t}\norm{s}^2 + 2t\norm{x^+ - \wb{x}}^2 + -t\norm{\wb{x} - x^+}^2 -t\norm{x - x^+}^2 + t\norm{x - \wb{x}}^2$. Since $\norm{x^+ - \wb{x}}^2$ has a negative coefficient, this term could be canceled from the right hand side of the inequality and the result will be: $f(x^+) -f(x) + \inprod{s}{x - \wb{x}} \leq t(\norm{x - \wb{x}}^2 - \norm{x - x^+}^2) + \frac{1}{2t}\norm{s}^2$.
    \end{proof}
\end{proposition}
We define gap function 
\begin{align*}
    Q(z_{t+1}^k,z) = f(x_{t+1}^k) + \lambda g(x_{t+1}^k) - \big(f(x) + {\lambda_{t+1}^k}^\top g(x)\big).
\end{align*}
Also
\begin{align}
    D_x^k(x) &= \tfrac{1}{T^{k-1}}\tsum_{l=1}^{T^{k-1}}\norm{x - x_l^{k-1}}^2 & D_{\lambda}^k(\lambda) &= \tfrac{1}{T^{k-1}}\tsum_{l=1}^{T^{k-1}}\norm{\lambda - \lambda_l^{k-1}}^2.  \label{eq:p6}
\end{align}
And
\begin{align}
    \norm{x - \wt{x}^{k-1}}^2 &\leq D_x^k(x) & \norm{\lambda - \wt{\lambda}^{k-1}}^2 &\leq D_{\lambda}^k(\lambda), \label{eq:p7}
\end{align}
where, \eqref{eq:p7} is the direct result of applying Jensen's inequality and triangular inequality for \eqref{eq:p6}. We also use Jensen's inequality for the gap function with its expectation form as
\begin{align}
    E[Q(\wt{z}^k,z)] &\leq \tfrac{1}{T^k}\tsum_{t=1}^{T^k}E[Q(z_t^k,z)] & E[Q(\wt{Z}^K,z)] &\leq \tfrac{1}{S^K}\tsum_{k=1}^{K}E[Q(\wt{z}^k,z)], \label{eq:p8}
\end{align}
for a fixed $z$, where, $\wt{z}^k = \tfrac{1}{T^k}\tsum_{t=1}^{T^k}z_t^k$, $S^K =\tsum_{k=1}^{K}T^k$, and $\wt{Z}^K = \tfrac{1}{S^K}\tsum_{k=1}^{K}T^k\wt{z}^k$. We also use Cauchy-Schwarz inequality frequently in this work. 
\begin{align}
    \inprod{a}{b} \leq \norm{a}\norm{b} \label{eq:p9}
\end{align}
and Arithmetic-Geometric mean inequalities
\begin{align}
    \norm{a}\norm{b} \leq \tfrac{\alpha}{2}\norm{a}^2 + \tfrac{1}{2\alpha}\norm{b}^2 \label{eq:p10}
\end{align}
frequently for any vectors $a$ and $b$, in this work.
Young's inequality for scalar variable $x$ where $\beta > 0$ 
\begin{align}
    \gamma x - \beta x^2 \leq \tfrac{\gamma^2}{4\beta} \label{eq:p11}
\end{align}

%% file: sec-full-algorithm.tex
\subsection{Convergence analysis for the convex objective.} 
\label{sec-method2.2}
Algorithm \ref{alg:alg2} is a primal-dual method. We construct a gap function and derive an upper bound on it in order to establish convergence rates and gradient-complexity bounds.
\begin{lemma}
    Let $\{x_t^k,\lambda_t^k\}$ be the sequence generated by Algorithm \ref{alg:alg2}. The gap inequality is
\begin{align}
    & Q(z_{t+1}^k, z) 
    %\nonumber\\&
    \leq \inprod{\delta_t^{\lambda^k}}{\lambda_{t+1}^k - \lambda} 
    - \inprod{\bar{q}_{t+1}^k}{\lambda_{t+1}^k - \lambda} 
    + \inprod{\bar{q}_t^k}{\lambda_t^k - \lambda} 
    %\nonumber \\ & 
    + \inprod{q_t^k - \bar{q}_t^k}{\lambda_t^k - \lambda} 
    \nonumber \\
    & + \tfrac{1 - \gamma_\lambda^k}{2\sigma^k} \paran[\big]{\norm{\lambda - \lambda_t^k}^2 
    - \norm{\lambda - \lambda_{t+1}^k}^2 
    - \norm{\lambda_t^k - \lambda_{t+1}^k}^2} 
    %\nonumber \\& 
    + \tfrac{\gamma_\lambda^k}{2\sigma^k} \paran[\big]{ D_\lambda^k(\lambda) 
    - \norm{\lambda - \lambda_{t+1}^k}^2 
    - D_\lambda^k(\lambda_{t+1}^k) } \nonumber \\
    & + \tfrac{L_g\norm{\lambda}}{2} \norm{x_t^k - x_{t+1}^k}^2  + \inprod{\delta_t^{x^k}}{x - x_{t+1}^k} 
    + \tfrac{L_f}{2} \gnorm{x^k_{t+1} - x^k_t}{}{2} + M_g \norm{x_t^k - x_{t-1}^k} \norm{\lambda_{t+1}^k - \lambda_t^k} \nonumber \\
    & + \tfrac{1 - \gamma_x^k}{2\tau^k} \big( \norm{x - x_t^k}^2 
    - \norm{x - x_{t+1}^k}^2 
    - \norm{x_{t+1}^k - x_t^k}^2 \big) 
    %\nonumber \\& 
    + \tfrac{\gamma_x^k}{2\tau^k} \big( D_x^k(x) 
    - \norm{x - x_{t+1}^k}^2 
    - D_x^k(x_{t+1}^k) \big), \label{eq:1}
\end{align}
where $\delta_t^{x^k} := \beta_t^k -  \grad \varphi(x^k_t,  \lambda^k_{t+1})$, $\delta_t^{\lambda^k} := \xi_t^k - g(x^k_t)$ and $\wb{q}_t^k := \ell_g(x^k_t) - g(x^k_{t-1})$.
\end{lemma}
\begin{proof}
    Using the optimality condition for line 12 of Algorithm \ref{alg:alg2}:
\begin{align}
    &\inprod{\xi_t^k + q_t^k}{\lambda - \lambda_{t+1}^k} \leq \tfrac{1}{\sigma^k} \inprod{\lambda_{t+1}^k - \wh{\lambda}_t^k}{\lambda - \lambda^k_{t+1}} = \tfrac{1}{\sigma^k} \inprod{\lambda_{t+1}^k - (1-\gamma_\lambda^k) \lambda_t^k - \gamma_\lambda^k \wt{\lambda}^{k-1}}{\lambda - \lambda_{t+1}^k} \nonumber\\
	& = \tfrac{1 - \gamma_\lambda^k}{\sigma^k} \inprod{\lambda_{t+1}^k - \lambda_t^k}{\lambda - \lambda_{t+1}^k} + \tfrac{\gamma_\lambda^k}{\sigma^k T^{k-1}} \tsum_{l=1}^{T^{k -1}} \inprod{\lambda_{t+1}^k - \lambda_l^{k-1}}{\lambda - \lambda_{t+1}^k} \nonumber\\
	& \leq \tfrac{1 - \gamma_\lambda^k}{2\sigma^k}\big( \norm{\lambda - \lambda_t^k}^2 - \norm{\lambda - \lambda_{t+1}^k}^2 - \norm{\lambda_t^k - \lambda_{t+1}^k}^2 \big) \nonumber\\
    & + \tfrac{\gamma_\lambda^k}{2\sigma^k T^{k-1}} \tsum_{l=1}^{T^{k -1}} \big(\norm{\lambda - \lambda_l^{k-1}}^2 - \norm{\lambda - \lambda_{t+1}^k}^2 - \norm{\lambda_l^{k-1} - \lambda_{t+1}^k}^2 \big), \label{eq:a8}
\end{align}
where the first inequality in \eqref{eq:a8} comes from Proposition \ref{theorem:1}, the first equality comes from the definition of $\wh{\lambda}_t^k$, the second equality comes from the definition of $\wt{\lambda}^{k-1}$, and the second inequality comes from the Pythagorean relation for inner products.
\begin{align}
	&\inprod{\xi_t^k + q_t^k}{\lambda_{t+1}^k - \lambda} = \inprod{\delta_t^{\lambda^k}}{\lambda_{t+1}^k - \lambda} + \inprod{g(x_t^k)}{\lambda_{t+1}^k - \lambda} + \inprod{q_t^k}{\lambda_t^k - \lambda} + \inprod{q_t^k}{\lambda_{t+1}^k - \lambda_t^k} \nonumber \\
	\Rightarrow & \inprod{\ell_g(x_{t+1}^k)}{\lambda - \lambda_{t+1}^k} \leq \inprod{\delta_t^{\lambda^k}}{\lambda_{t+1}^k - \lambda} - \inprod{\bar{q}_{t+1}^k}{\lambda_{t+1}^k - \lambda} + \inprod{q_t^k - \bar{q}_{t}^k}{\lambda_t^k - \lambda} \nonumber\\
    &+ \inprod{\bar{q}_t^k}{\lambda_t^k - \lambda} + \inprod{q_t^k}{\lambda_{t+1}^k - \lambda_t^k} + \tfrac{1 - \gamma_\lambda^k}{2\sigma^k} \left( \norm{\lambda - \lambda_t^k}^2 - \norm{\lambda - \lambda_{t+1}^k}^2 - \norm{\lambda_t^k - \lambda_{t+1}^k}^2 \right) \nonumber \\
	& + \tfrac{\gamma^k_\lambda}{2\sigma^k} \left( D_\lambda^k(\lambda) - \norm{\lambda - \lambda_{t+1}^k}^2 - D_\lambda^k(\lambda_{t+1}^k) \right), \label{eq:a9}
\end{align}  
The first inequality in \eqref{eq:a9} comes from $\delta_t^{\lambda^k}$ and $\bar{q}^k_t$, and applying \eqref{eq:a8}. Recall \eqref{eq:p6}. Also, with \eqref{eq:p5}
\begin{align}
	\norm{q_t^k} &= \norm{\ell_{g_{j_t^k}}(x_t^k) - g_{j_t^k}(x_{t-1}^k)}
    %\nonumber \\&
    = \norm{g_{j_t^k}(x_{t-1}^k) + \inprod{\grad g_{j_t^k}(x_{t-1}^k)}{x_t^k - x_{t-1}^k} - g_{j_t^k}(x_{t-1}^k)} \leq M_g \norm{x_t^k - x_{t-1}^k}. \label{eq:a10} 
\end{align}  
Using \eqref{eq:a10} in \eqref{eq:a9} and \eqref{eq:p9} for $\inprod{q_t^k}{\lambda_{t+1}^k - \lambda_t^k}$
\begin{align}
    &\inprod{\ell_g(x_{t+1}^k)}{\lambda - \lambda_{t+1}^k} \leq \inprod{\delta_t^{\lambda^k}}{\lambda_{t+1}^k - \lambda} - \inprod{\bar{q}_{t+1}^k}{\lambda_{t+1}^k - \lambda} + \inprod{\bar{q}_t^k}{\lambda_t^k - \lambda} + \inprod{q_t^k - \bar{q}_{t}^k}{\lambda_t^k - \lambda} \nonumber\\
    & + M_g \norm{x_t^k - x_{t-1}^k}\norm{\lambda_{t+1}^k - \lambda_t^k} + \tfrac{1 - \gamma_\lambda^k}{2\sigma^k} \big( \norm{\lambda - \lambda_t^k}^2 - \norm{\lambda - \lambda_{t+1}^k}^2 - \norm{\lambda_t^k - \lambda_{t+1}^k}^2 \big) \nonumber\\
    & + \tfrac{\gamma^k_\lambda}{2\sigma^k} \big( D_\lambda^k(\lambda) - \norm{\lambda - \lambda_{t+1}^k}^2 - D_\lambda^k(\lambda_{t+1}^k) \big). \label{eq:a11}
\end{align}
Inner product $\inprod{\ell_g(x_{t+1}^k)}{\lambda - \lambda_{t+1}^k}$ is $\inprod{\ell_g(x_{t+1}^k)}{\lambda} + \inprod{\ell_g(x_{t+1}^k)}{-\lambda_{t+1}^k}$ and
\begin{align}
    \inprod{g(x_{t+1}^k) - \ell_g(x_{t+1}^k)}{\lambda} &\leq \tfrac{\norm{\lambda}L_g}{2}\norm{x_t^k - x_{t+1}^k}^2,\label{eq:a12}
\end{align}
where in \eqref{eq:a12}, we use \eqref{eq:p9} and \eqref{eq:p2}. Taking $\inprod{\ell_g(x_{t+1}^k)}{- \lambda_{t+1}^k}$ to the right hand side of \eqref{eq:a11} and using \eqref{eq:a12}
\begin{align}
    & \inprod{g(x_{t+1}^k)}{\lambda} \leq \inprod{\delta_t^{\lambda^k}}{\lambda_{t+1}^k - \lambda} 
    - \inprod{\bar{q}_{t+1}^k}{\lambda_{t+1}^k - \lambda} 
    + \inprod{\bar{q}_t^k}{\lambda_t^k - \lambda} + \inprod{q_t^k - \bar{q}_t^k}{\lambda_t^k - \lambda} \nonumber\\ 
    &+ M_g \norm{x_t^k - x_{t-1}^k} \norm{\lambda_{t+1}^k - \lambda_t^k}  + \tfrac{1 - \gamma_\lambda^k}{2\sigma^k} ( \norm{\lambda - \lambda_t^k}^2 
    - \norm{\lambda - \lambda_{t+1}^k}^2 
    - \norm{\lambda_t^k - \lambda_{t+1}^k}^2 ) \nonumber \\
    & + \tfrac{\gamma_\lambda^k}{2\sigma^k} \big(D_\lambda^k(\lambda) 
    - \norm{\lambda - \lambda_{t+1}^k}^2 
    - D_\lambda^k(\lambda_{t+1}^k)\big) + \inprod{\ell_g(x_{t+1}^k)}{\lambda_{t+1}^k} 
    + \tfrac{\norm{\lambda}L_g}{2} \norm{x_t^k - x_{t+1}^k}^2. \label{eq:a13}
\end{align}
Using the optimality of $x_{t+1}^k$ from line 16 of Algorithm \ref{alg:alg2}
\begin{align}
	&\inprod{\beta_t^k}{x_{t+1}^k - x} \leq \tfrac{1}{\tau^k} \inprod{x_{t+1}^k - \wh{x}_t^k}{x - x_{t+1}^k}  = \tfrac{1 - \gamma_x^k}{\tau^k} \inprod{x_{t+1}^k - x_t^k}{x - x_{t+1}^k} 
	 + \tfrac{\gamma_x^k}{\tau^k} \inprod{x_{t+1}^k - \wt{x}^{k-1}}{x - x_{t+1}^k} \nonumber \\
	& \leq \tfrac{1 - \gamma_x^k}{2\tau^k} \big( \norm{x - x_t^k}^2 - \norm{x - x_{t+1}^k}^2 - \norm{x_t^k - x_{t+1}^k}^2 \big) \nonumber\\
    & + \tfrac{\gamma_x^k}{2\tau^k T^{k-1}} \tsum_{l=1}^{T^{k-1}} \big( \norm{x - x_l^{k-1}}^2 - \norm{x - x_{t+1}^k}^2 - \norm{x_l^{k-1} - x_{t+1}^k}^2 \big), \label{eq:a14}
\end{align}
where  the first equality comes from the definition of $\wh{x}_t^k$, the first inequality in \eqref{eq:a14} comes from  Proposition \ref{theorem:1}, and the second inequality comes from Proposition \ref{theorem:2}. Using $ \beta_t^k = \delta_t^{x^k} + \nabla_x\varphi(x_t^k,\lambda_{t+1}^k)$
\begin{align}
	& \inprod{\delta_t^{x^k}}{x_{t+1}^k - x} \leq \inprod{\nabla_x \varphi(x_t^k,\lambda_{t+1}^k)}{x - x_{t+1}^k} + \tfrac{1-\gamma_x^k}{2\tau^k} (\norm{x - x_t^k}^2 - \norm{x - x_{t+1}^k}^2 - \norm{x_{t+1}^k - x_t^k}^2) \nonumber\\
    & + \tfrac{\gamma_x^k}{2\tau^k} \big(D_x^k(x) - \norm{x - x_{t+1}^k}^2 - D_x^k(x_{t+1}^k)\big), \label{eq:a15}
\end{align}
where \eqref{eq:a15} follows from \eqref{eq:a14},
the identity $\beta_t^k=\delta_t^{x^k}+\nabla_x\varphi(x_t^k,\lambda_{t+1}^k)$, and \eqref{eq:p6}. 
\begin{align}
    &\inprod{\grad_x \varphi(x_t^k, \lambda_{t+1}^k)}{x - x_{t+1}^k} = \inprod{\grad_x \varphi(x_t^k, \lambda_{t+1}^k)}{x - x_t^k} 
    + \inprod{\grad_x \varphi(x_t^k, \lambda_{t+1}^k)}{x_t^k - x_{t+1}^k} \nonumber \\
    & \leq \varphi(x, \lambda_{t+1}^k) - \varphi(x_t^k, \lambda_{t+1}^k)  + \inprod{\grad_x \varphi(x_t^k, \lambda_{t+1}^k)}{x_t^k - x_{t+1}^k}, \label{eq:a16}
\end{align}
where the inequality in \eqref{eq:a16} follows from the convexity of $\varphi(x,\lambda_{t+1}^k)$ with respect to $x$. Using \eqref{eq:a16} in the gap function and rearranging terms
\begin{align}
	& \inprod{\delta_t^{x^k}}{x_{t+1}^k - x} - \varphi(x, \lambda_{t+1}^k) + \varphi(x_t^k, \lambda_{t+1}^k) \leq \inprod{\nabla_x \varphi(x_t^k, \lambda_{t+1}^k)}{x_t^k - x_{t+1}^k} \nonumber\\
    & + \tfrac{1-\gamma_x^k}{2\tau^k} (\norm{x - x_t^k}^2 - \norm{x - x_{t+1}^k}^2 - \norm{x_{t+1}^k - x_t^k}^2) + \tfrac{\gamma_x^k}{2\tau^k} \big(D_x^k(x) - \norm{x - x_{t+1}^k}^2 - D_x^k(x_{t+1}^k)\big). \label{eq:a18}  
\end{align}
Upper bound of $\inprod{\nabla_x \varphi(x_t^k, \lambda_{t+1}^k)}{x_t^k - x_{t+1}^k}$
\begin{align}
    & \inprod{\nabla_x \varphi(x_t^k, \lambda_{t+1}^k)}{x_t^k - x_{t+1}^k} = \inprod{\nabla f(x_t^k)}{x_t^k - x_{t+1}^k} 
    + \inprod{\lambda_{t+1}^k{^\top} \nabla g(x_t^k)}{x_t^k - x_{t+1}^k} \nonumber \\
    &= \inprod{\nabla f(x_t^k)}{x_t^k - x_{t+1}^k} 
    + \inprod{\lambda_{t+1}^k }{\nabla g(x_t^k)^\top(x_t^k - x_{t+1}^k)} =  \inprod{\nabla f(x_t^k)}{x_t^k - x_{t+1}^k} 
    + \inprod{\lambda_{t+1}^k }{g(x^k_t)-\ell_g(x_{t+1}^k)} \nonumber \\
    &\le f(x^k_t) - f(x^k_{t+1}) + \tfrac{L_f}{2} \gnorm{x^k_{t+1}- x^k_t}{}{2} + \inprod{\lambda_{t+1}^k }{g(x^k_t)-\ell_g(x_{t+1}^k)} \nonumber\\
    & = \varphi(x^k_t, \lambda^k_{t+1}) - f(x^k_{t+1}) + \tfrac{L_f}{2} \gnorm{x^k_{t+1}- x^k_t}{}{2} - \inprod{\lambda_{t+1}^k }{\ell_g(x_{t+1}^k)}, \label{eq:a300}
\end{align}
where in \eqref{eq:a300}, second-last equality, we use \eqref{eq:p5}, and the inequality comes from \eqref{eq:p2} and $f(x^k_t) + \lambda_{t+1}^k{^\top} g(x_t^k) = \varphi(x^k_t, \lambda^k_{t+1})$. Finally, for \eqref{eq:a18}
\begin{align}
    &\inprod{\delta_t^{x^k}}{x_{t+1}^k - x} - \varphi(x, \lambda_{t+1}^k) + \varphi(x_t^k, \lambda_{t+1}^k) \leq \varphi(x^k_t, \lambda^k_{t+1}) - f(x^k_{t+1}) + \tfrac{L_f}{2} \gnorm{x^k_{t+1}- x^k_t}{}{2} \nonumber\\
    & - \inprod{\lambda_{t+1}^k }{\ell_g(x_{t+1}^k)} + \tfrac{1-\gamma_x^k}{2\tau^k} (\norm{x - x_t^k}^2 - \norm{x - x_{t+1}^k}^2 - \norm{x_{t+1}^k - x_t^k}^2) \nonumber\\
    & + \tfrac{\gamma_x^k}{2\tau^k} \big(D_x^k(x) - \norm{x - x_{t+1}^k}^2 - D_x^k(x_{t+1}^k)\big). \label{eq:a19}
\end{align}
Adding \eqref{eq:a13} and \eqref{eq:a19} and rearranging the resulting terms, yields \eqref{eq:1}.
\end{proof}
For further analysis, we define new auxiliary sequences for the inner products that are related to $\delta_t^{x^k}$,$\delta_t^{\lambda^k}$ and $q_t^k - \wb{q}_t^k$. 
\begin{definition} \label{eq:a20}
    Define the auxiliary sequences $u_t^k$, $v_t^k$, and $w_t^k$ as follows:
    \begin{align*}
        u_{t+1}^k &= \arg\min_{x \in X} {-\inprod{\delta_t^{x^k}}{x} + \tfrac{1}{\eta^k} \|x - u_t^k\|^2 }, \nonumber \\
        v_{t+1}^k &= \arg\min_{\lambda \geq0}{\inprod{\delta_t^{\lambda^k}}{\lambda} + \tfrac{1}{\eta^k} \|\lambda - v_t^k\|^2 }, \nonumber \\
        w_{t+1}^k &= \arg\min_{\lambda \geq0}{\inprod{q_t^k - \bar{q}_t^k}{\lambda} + \tfrac{1}{\eta^k} \|\lambda - w_t^k\|^2 }.
    \end{align*}
\end{definition}
\begin{lemma} \label{theorem:3}
    The following inequalities hold based on the definition of the auxiliary variables:
    \begin{align}
        \inprod{\delta_t^{x^k}}{x - x_{t+1}^k} &\leq \tfrac{1}{\eta^k} (\norm{x - u_t^k}^2 - \norm{x - u_{t+1}^k}^2) +\eta^k \norm{\delta_t^{x^k}}^2 + \inprod{\delta_t^{x^k}}{u_t^k - x_t^k} + \tfrac{1}{\eta^k} \norm{x_t^k - x_{t+1}^k}^2, \label{eq:a21} 
    \end{align}
    \begin{align}
        \inprod{\delta_t^{\lambda^k}}{\lambda_{t+1}^k - \lambda} &\leq \tfrac{1}{\eta^k} (\norm{\lambda - v_t^k}^2 - \norm{\lambda - v_{t+1}^k}^2) +\eta^k \norm{\delta_t^{\lambda^k}}^2 + \inprod{\delta_t^{\lambda^k}}{\lambda_t^k - v_t^k} + \tfrac{1}{\eta^k} \norm{\lambda_t^k - \lambda_{t+1}^k}^2, 
    \end{align}
    \begin{align}
        \inprod{q_t^k - \bar{q}_t^k}{\lambda_t^k - \lambda} &\leq \tfrac{1}{\eta^k} (\norm{\lambda - w_t^k}^2 - \norm{\lambda - w_{t+1}^k}^2) + \tfrac{\eta^k}{2} \norm{q_t^k - \bar{q}_t^k}^2 + \inprod{q_t^k - \bar{q}_t^k}{\lambda_t^k - w_t^k}.  
    \end{align}
\end{lemma}
\begin{proof} We have 
    $\inprod{\delta_t^{x^k}}{x - x_{t+1}^k}
    = 
    %\inprod{\delta_t^{x^k}}{x - x_t^k} + 
    \inprod{\delta_t^{x^k}}{x_t^k - x_{t+1}^k}
    % $. Now for expression $\inprod{\delta_t^{x^k}}{x - x_t^k}$, $\inprod{\delta_t^{x^k}}{x - x_t^k} 
    %= 
    +
    \inprod{\delta_t^{x^k}}{x - u_t^k} + \inprod{\delta_t^{x^k}}{u_t^k - x_t^k}$ and
    \begin{align*}
        \inprod{\delta_t^{x^k}}{x - u_t^k} + \inprod{\delta_t^{x^k}}{u_t^k - x_t^k} \leq \tfrac{1}{\eta^k}(\norm{x - u_t^k}^2 - \norm{x - u_{t+1}^k}^2) + \tfrac{\eta^k}{2}\norm{\delta_t^{x^k}}^2 + \inprod{\delta_t^{x^k}}{u_t^k - x_t^k},
    \end{align*}
    where %in the inequality above, 
    we use Lemma \ref{theorem:2} for $\inprod{\delta_t^{x^k}}{x - u_t^k}$. For $\inprod{\delta_t^{x^k}}{x_t^k - x_{t+1}^k}$ the upper bound is $\frac{\eta^k}{2}\norm{\delta_t^{x^k}}^2 + \frac{1}{\eta^k}\norm{x_t^k - x_{t+1}^k}^2$, which we use \eqref{eq:p9}, \eqref{eq:p10} and strong convexity of Euclidean norm function. Combining the above bounds yields the desired estimate for
$\inprod{\delta_t^{x^k}}{x - x_{t+1}^k}$. The proofs of the other two inequalities are analogous.
\end{proof}
\begin{lemma}
    The expected gap generated by Algorithm \ref{alg:alg2} satisfies, 
\begin{align}
    & E[Q(\wt{Z}^K,z^*)] \leq \tfrac{1}{S^K}\big(\tfrac{1}{\eta^1}E[\norm{\lambda^* - v_0^1}^2] + \tsum_{k=1}^{K-1}(\tfrac{1}{\eta^{k+1}} - \tfrac{1}{\eta^k})E[\norm{\lambda^* - v_{T^k}^k}^2] - \tfrac{1}{\eta^K}E[\norm{\lambda^* - v_{T^K}^K}^2]  \nonumber\\
    & + \tsum_{k=1}^{K}\eta^k\tsum_{t=0}^{T^k-1}E[\norm{\delta_t^{\lambda^k}}^2] + \tsum_{k=1}^{K}\tsum_{t=0}^{T^k-1}E[\inprod{\delta_t^{\lambda^k}}{\lambda_t^k - v_t^k}] + E[\inprod{\bar{q}_0^1}{\lambda_0^1 - \lambda^*}] - E[\inprod{\bar{q}_{T^K}^K}{\lambda_{T^K}^K - \lambda^*}] \nonumber\\
    & + \tfrac{1}{\eta^1}E[\norm{\lambda^* - w_0^1}^2] + \tsum_{k=1}^{K-1}(\tfrac{1}{\eta^{k+1}} - \tfrac{1}{\eta^k})E[\norm{\lambda^* - w_{T^k}^k}^2] - \tfrac{1}{\eta^K}E[\norm{\lambda^* - w_{T^K}^K}^2] \nonumber\\
    & + \tsum_{k=1}^{K}\tfrac{\eta^k}{2}\tsum_{t=0}^{T^k-1}E[\norm{q_t^k - \bar{q}_t^k}^2] +  \tsum_{k=1}^{K}\tsum_{t=0}^{T^k-1}E[\inprod{q_t^k - \bar{q}_t^k}{\lambda_t^k - w_t^k}] \nonumber\\
    & + \tsum_{k=1}^{K}(\tfrac{1}{\eta^k} + \tfrac{M_g}{2} - \tfrac{1 - \gamma_\lambda^k}{2\sigma^k})\tsum_{t=0}^{T^k-1}E[\norm{\lambda_t^k - \lambda_{t+1}^k}^2] -  \tsum_{k=1}^{K}\tfrac{\gamma_\lambda^k}{2\sigma^k}\tsum_{t=0}^{T^k-1}E[D_\lambda^k(\lambda_{t+1}^k)]\nonumber\\
    & + \tfrac{1 - \gamma_\lambda^1}{2\sigma^1}E[\norm{\lambda^* - \lambda_0^1}^2] + \tsum_{k=1}^{K-1}(\tfrac{1 - \gamma_\lambda^{k+1}}{2\sigma^{k+1}} - \tfrac{1 - \gamma_\lambda^k}{2\sigma^k})E[\norm{\lambda^* - \lambda_{T^k}^k}^2] - \tfrac{1 - \gamma_\lambda^K}{2\sigma^K}E[\norm{\lambda^* - \lambda_{T^K}^K}^2] \nonumber \\
    & + \tfrac{\gamma_\lambda^1T^1}{2\sigma^1}E[D_\lambda^1(\lambda^*)] + \tsum_{k=1}^{K-1}(\tfrac{\gamma_\lambda^{k+1}T^{k+1}}{2\sigma^{k+1}} - \tfrac{\gamma_\lambda^kT^k}{2\sigma^k})E[D_\lambda^{k+1}(\lambda^*)] - \tfrac{\gamma_\lambda^KT^K}{2\sigma^K}E[D_\lambda^{K+1}(\lambda^*)] \nonumber\\
    & + \tsum_{k=1}^{K}(\tfrac{\norm{\lambda^*}L_g}{2} + \tfrac{M_g}{2} + \tfrac{1}{\eta^k} + \tfrac{L_f}{2} - \tfrac{1 - \gamma_x^k}{2\tau^k})\tsum_{t=0}^{T^k-1}E[\norm{x_t^k - x_{t+1}^k}^2]  + \tsum_{k=1}^{K}\eta^k\tsum_{t=0}^{T^k-1}E[\norm{\delta_t^{x^k}}^2] \nonumber\\
    & + \tfrac{1}{\eta^1}E[\norm{x^* - u_0^1}^2] + \tsum_{k=1}^{K-1}(\tfrac{1}{\eta^{k+1}} - \tfrac{1}{\eta^k})E[\norm{x^* - u_{T^k}^k}^2] - \tfrac{1}{\eta^K}E[\norm{x^* - u_{T^K}^K}^2] + \tsum_{k=1}^{K}\tsum_{t=0}^{T^k-1}E[\inprod{\delta_t^{x^k}}{u_t^k - x_t^k}]\nonumber \\
    & + \tfrac{1 - \gamma_x^1}{2\tau^1}E[\norm{x^* - x_0^1}^2] 
    + \tsum_{k=1}^{K-1}(\tfrac{1 - \gamma_x^{k+1}}{2\tau^{k+1}} - \tfrac{1 - \gamma_x^k}{2\tau^k})E[\norm{x^* - x_{T^k}^k}^2] - \tfrac{1 - \gamma_x^K}{2\tau^K}E[\norm{x^* - x_{T^K}^K}^2] \nonumber \\
    & + \tfrac{\gamma_x^1T^1}{2\tau^1} E[D_x^1(x^*)] + \tsum_{k=1}^{K-1}(\tfrac{\gamma_x^{k+1}T^{k+1}}{2\tau^{k+1}} - \tfrac{\gamma_x^kT^k}{2\tau^k})E[D_x^{k+1}(x^*)] - \tfrac{\gamma_x^KT^K}{2\tau^K}E[D_x^{K+1}(x^*)]  \nonumber\\
    & - \tsum_{k=1}^{K}\tfrac{\gamma_x^k}{2\tau^k}\tsum_{t=0}^{T^k-1}E[D_x^k(x_{t+1}^k)]\big). \label{eq:a26}  
\end{align}
\end{lemma}
\begin{proof}
Using inequalities from Lemma \ref{theorem:3} in  \eqref{eq:1}, and setting $(x,\lambda,z) = (x^*,\lambda^*,z^*)$, and taking expectation on both sides of the inequality  
\begin{align}
    & E[Q(z_{t+1}^k, z^*)] \leq \tfrac{1}{\eta^k} (E[\norm{\lambda^* - v_t^k}^2] - E[\norm{\lambda^* - v_{t+1}^k}^2]) +\eta^k E[\norm{\delta_t^{\lambda^k}}^2] + E[\inprod{\delta_t^{\lambda^k}}{\lambda_t^k - v_t^k}] \nonumber\\
    & + \tfrac{1}{\eta^k} E[\norm{\lambda_t^k - \lambda_{t+1}^k}^2] - E[\inprod{\bar{q}_{t+1}^k}{\lambda_{t+1}^k - \lambda^*}] 
    + E[\inprod{\bar{q}_t^k}{\lambda_t^k - \lambda^*}] \nonumber \\
    & + \tfrac{1}{\eta^k} (E[\norm{\lambda^* - w_t^k}^2] - E[\norm{\lambda^* - w_{t+1}^k}^2]) + \tfrac{\eta^k}{2} E[\norm{q_t^k - \bar{q}_t^k}^2] \nonumber\\
    & + E[\inprod{q_t^k - \bar{q}_t^k}{\lambda_t^k - w_t^k}] 
    + M_g E[\norm{x_t^k - x_{t-1}^k}\norm{\lambda_{t+1}^k - \lambda_t^k}] \nonumber \\
    &+ \tfrac{1 - \gamma_\lambda^k}{2\sigma^k} \big( E[\norm{\lambda^* - \lambda_t^k}^2] 
    - E[\norm{\lambda^* - \lambda_{t+1}^k}^2] 
    - E[\norm{\lambda_t^k - \lambda_{t+1}^k}^2] \big) \nonumber \\
    & + \tfrac{\gamma_\lambda^k}{2\sigma^k} \big( E[D_\lambda^k(\lambda^*)] 
    - E[\norm{\lambda^* - \lambda_{t+1}^k}^2] 
    - E[D_\lambda^k(\lambda_{t+1}^k)] \big) \nonumber \\
    & + \tfrac{\norm{\lambda^*}L_g}{2} E[\norm{x_t^k - x_{t+1}^k}^2]  + \tfrac{1}{\eta^k} (E[\norm{x^* - u_t^k}^2] - E[\norm{x^* - u_{t+1}^k}^2]) +\eta^k E[\norm{\delta_t^{x^k}}^2] \nonumber\\
    & + E[\inprod{\delta_t^{x^k}}{u_t^k - x_t^k}] + \tfrac{1}{\eta^k} E[\norm{x_t^k - x_{t+1}^k}^2] + \tfrac{L_f}{2}E[\norm{x^k_{t+1} - x^k_t}^2] \nonumber \\
    & + \tfrac{1 - \gamma_x^k}{2\tau^k} \big( E[\norm{x^* - x_t^k}^2] 
    - E[\norm{x^* - x_{t+1}^k}^2] 
    - E[\norm{x_{t+1}^k - x_t^k}^2] \big) \nonumber \\
    & + \tfrac{\gamma_x^k}{2\tau^k} \big( E[D_x^k(x^*)] 
    - E[\norm{x^* - x_{t+1}^k}^2] 
    - E[D_x^k(x_{t+1}^k)] \big). \label{eq:a23}
\end{align}
Now, summing both sides of \eqref{eq:a23} from $t =0$ to $T^k-1$ for epoch $k$ and using \eqref{eq:p8}, we have
\begin{align}
    & T^kE[Q(\wt{z}^k,z^*)] \leq \tfrac{1}{\eta^k} (E[\norm{\lambda^* - v_0^k}^2] - E[\norm{\lambda^* - v_{T^k}^k}^2]) +\eta^k \tsum_{t=0}^{T^k-1}E[\norm{\delta_t^{\lambda^k}}^2]  \nonumber\\
    & + \tsum_{t=0}^{T^k-1}E[\inprod{\delta_t^{\lambda^k}}{\lambda_t^k - v_t^k}] + \tfrac{1}{\eta^k} \tsum_{t=0}^{T^k-1}E[\norm{\lambda_t^k - \lambda_{t+1}^k}^2] + E[\inprod{\bar{q}_0^k}{\lambda_0^k - \lambda^*}] - E[\inprod{\bar{q}_{T^k}^k}{\lambda_{T^k}^k - \lambda^*}] \nonumber \\
    & + \tfrac{1}{\eta^k}(E[\norm{\lambda^* - w_0^k}^2] - E[\norm{\lambda^* - w_{T^k}^k}^2]) + \tfrac{\eta^k}{2} \tsum_{t=0}^{T^k-1}E[\norm{q_t^k - \bar{q}_t^k}^2] \nonumber\\
    & + \tsum_{t=0}^{T^k-1}E[\inprod{q_t^k - \bar{q}_t^k}{\lambda_t^k - w_t^k}] 
    + \tsum_{t=0}^{T^k-1}M_gE[\norm{x_t^k - x_{t-1}^k} \norm{\lambda_{t+1}^k - \lambda_t^k}] \nonumber \\
    & + \tfrac{1 - \gamma_\lambda^k}{2\sigma^k} \bigg(E[\norm{\lambda^* - \lambda_0^k}^2] - E[\norm{\lambda^* - \lambda_{T^k}^k}^2] 
    - \tsum_{t=0}^{T^k-1}E[\norm{\lambda_t^k - \lambda_{t+1}^k}^2] \bigg) \nonumber \\
    & + \tfrac{\gamma_\lambda^k}{2\sigma^k} \bigg(T^kE[D_\lambda^k(\lambda^*)] 
    - \tsum_{t=0}^{T^k-1}E[\norm{\lambda^* - \lambda_{t+1}^k}^2] 
    - \tsum_{t=0}^{T^k-1}E[D_\lambda^k(\lambda_{t+1}^k)] \bigg) \nonumber \\
    & + \tsum_{t=0}^{T^k-1}(\tfrac{\norm{\lambda^*}L_g}{2}E[\norm{x_t^k - x_{t+1}^k}^2] ) + \tfrac{1}{\eta^k}(E[\norm{x^* - u_0^k}^2] - E[\norm{x^* - u_{T^k}^k}^2]) +\eta^k \tsum_{t=0}^{T^k-1}E[\norm{\delta_t^{x^k}}^2] \nonumber\\
    & + \tsum_{t=0}^{T^k-1}E[\inprod{\delta_t^{x^k}}{u_t^k - x_t^k}] + \frac{1}{\eta^k} \tsum_{t=0}^{T^k-1}E[\norm{x_t^k - x_{t+1}^k}^2] + \tfrac{L_f}{2}\tsum_{t=0}^{T^k-1}E[\norm{x^k_{t+1} - x^k_t}^2] \nonumber \\
    & + \tfrac{1 - \gamma_x^k}{2\tau^k}\bigg(E[\norm{x^* - x_0^k}^2] 
    - E[\norm{x^* - x_{T^k}^k}^2] 
    - \tsum_{t=0}^{T^k-1}E[\norm{x_{t+1}^k - x_t^k}^2] \bigg) \nonumber \\
    & + \tfrac{\gamma_x^k}{2\tau^k}\bigg(\tsum_{t=0}^{T^k-1}E[D_x^k(x^*)] 
    - \tsum_{t=0}^{T^k-1}E[\norm{x^* - x_{t+1}^k}^2] 
    - \tsum_{t=0}^{T^k-1}E[D_x^k(x_{t+1}^k)] \bigg). \label{eq:a25}
\end{align}
Summing \eqref{eq:a25} over k, dividing both sides by $S^K$, using \eqref{eq:p8}, the epoch initialization, and applying \eqref{eq:p10}, we obtain \eqref{eq:a26}.
\end{proof}
Henceforth, we study two cases separately: (1) semi-finite-sum, and (2) fully-finite-sum. In the next section, we consider the semi-finite-sum case.%We will see that for non constant step size policies, we need to find a constant upper bound for some terms associated with Lagrangian.      
\subsubsection{Semi-finite-sum case.} \label{sec-method2.2.1}
In this case, we assume $n=1$, meaning that the constraint mapping $g$ consists of a single component. We can only project onto the set $X$ and treat $\{g(x) \le 0\}$ as a function constraint.  %l has function constraints of the form $g(x)$, which does not have a finite-sum structure, so the gradient calculation is not computationally costly yet, it is still hard to do the projection operations on the function constraint. 
The objective $f$ still has a finite-sum structure. First, we derive upper bounds on the expected values of the terms $\norm{\delta_t^{x^k}}^2$, $\norm{\delta_t^{\lambda^k}}^2$, and $\norm{q_t^k - \wb{q}_t^k}^2$. %Notice that, all expectations are with respect to the random variables $i_t^k$, $j_t^k$, and $p_t^k$ from Algorithm \ref{alg:alg2}. % which come from the selection of the function components from the objective function or the function constraints.   
\begin{lemma} \label{eq:C1}
    The following relations hold in expectation with respect to the randomness in Algorithm \ref{alg:alg2}:
    \begin{align*}
        E[\norm{\delta_t^{x^k}}^2] & \leq 2L_f^2(E[D_x^k(x_{t+1}^k)] + E[\norm{x^k_t -x^k_{t+1}}^2]) & E[\norm{\delta_t^{\lambda^k}}^2] & = 0 & E[\norm{q_t^k - \wb{q}_t^k}^2] = 0,
    \end{align*}
    \begin{align}
        E[- \inprod{\bar{q}_{T^K}^K}{\lambda_{T^K}^K - \lambda^*} -\tfrac{1 - \gamma_{\lambda}^K}{4\sigma^K}\norm{\lambda^* - \lambda_{T^K}^K}^2] \leq \tfrac{M_g^2\sigma^K}{1-\gamma_{\lambda}^K}E[\norm{x_{T^K-1}^K - x_{T^K}^K}^2]. \label{eq:p15}
    \end{align}
\end{lemma}
\begin{proof}
    Recall that $\delta_t^{x^k} = \beta_t^k - \grad \varphi(x_t^k,  \lambda_{t+1}^k)$. Using $\beta_t^k$ from Algorithm \ref{alg:alg2}, noting that 
    $G_x = \grad f(\tilde{x}^{k-1}) + \grad g(\tilde{x}^{k-1}) \wt{\lambda}^{k-1}$ and $g_{p_t^k}(x) = g_{j_t^k}(x) = g(x)$ for the semi-finite-sum case, we have
     \begin{align*}
         \delta_t^{x^k} = \grad f_{i_t^k}(x^k_t) - \grad f_{i_t^k}(\wt{x}^{k-1}) - [\grad f(x_t^k) - \grad f(\wt{x}^{k-1})].
     \end{align*}
     By taking expectation of $\norm{\delta_t^{x^k}}^2$, we obtain 
    \begin{align}
        &E[\norm{\delta_t^{x^k}}^2] = E\big[\norm{\grad f_{i_t^k}(x^k_t) - \grad f_{i_t^k}(\wt{x}^{k-1}) - [\grad f(x_t^k) - \grad f(\wt{x}^{k-1})]}^2\big] \leq E[\norm{\grad f_{i_t^k}(x^k_t) - \grad f_{i_t^k}(\wt{x}^{k-1})}^2] \nonumber\\
        &  \leq L_f^2E[\norm{x^k_t - x^k_{t+1} +x^k_{t+1} -\wt{x}^{k-1}}^2] \leq 2L_f^2(E[D_x^k(x_{t+1}^k)] +E[\norm{x^k_t -x^k_{t+1}}^2]), 
    \end{align}
    where the first inequality is due to $Var[X] \leq E[X^2]$ for any random variable $X$, the second inequality comes from \eqref{eq:p2}, and the third inequality comes from the triangular inequality, and \eqref{eq:p7}. For $\delta_t^{\lambda^k}$, 
    \begin{align*}
        \delta_t^{\lambda^k} = \xi_t^k - g(x^k_t) = g(x^k_t) - g(\wt{x}^{k-1}) + G_{\lambda} - g(x^k_t) = 0,
    \end{align*}
    so $E[\norm{\delta_t^{\lambda^k}}^2] = 0$ and
    \begin{align*}
        q_t^k - \wb{q}_t^k = \ell_g(x^k_t) - g(x^k_{t-1}) - [\ell_g(x^k_t) - g(x^k_{t-1})] = 0,
    \end{align*}
    so $E[\norm{q_t^k - \wb{q}_t^k}^2] = 0$. Also,
    \begin{align*}
        & E[- \inprod{\bar{q}_{T^K}^K}{\lambda_{T^K}^K - \lambda^*} -\tfrac{(1 - \gamma_{\lambda}^K)}{4\sigma^K}\norm{\lambda^* - \lambda_{T^K}^K}^2] \leq E[ \norm{\bar{q}_{T^K}^K}\norm{\lambda_{T^K}^K - \lambda^*} -\tfrac{(1 - \gamma_{\lambda}^K)}{4\sigma^K}\norm{\lambda^* - \lambda_{T^K}^K}^2]\nonumber\\
        & = E[M_g\norm{x_{T^K-1}^K - x_{T^K}^K}\norm{\lambda_{T^K}^K - \lambda^*} -\tfrac{1 - \gamma_{\lambda}^K}{4\sigma^K}\norm{\lambda^* - \lambda_{T^K}^K}^2] \leq \tfrac{M_g^2\sigma^K}{1-\gamma_{\lambda}^K}E[\norm{x_{T^K-1}^K - x_{T^K}^K}^2],
    \end{align*}
    where the first inequality comes from \eqref{eq:p9}, the second inequality comes from \eqref{eq:p1}, and the last inequality comes from \eqref{eq:p11}.
\end{proof}
Using the upper bounds in Lemma \ref{eq:C1} and accumulating similar terms in \eqref{eq:a26}
\begin{align}
    &E[Q(\wt{Z}^K,z^*)] \leq  \tfrac{1}{S^K}\bigg(\tfrac{1 - \gamma_\lambda^1}{2\sigma^1}E[\norm{\lambda^* - \lambda_0^1}^2] + \tsum_{k=1}^{K-1}(\tfrac{1 - \gamma_\lambda^{k+1}}{2\sigma^{k+1}} - \tfrac{1 - \gamma_\lambda^k}{2\sigma^k})E[\norm{\lambda^* - \lambda_{T^k}^k}^2] - \tfrac{1 - \gamma_\lambda^K}{4\sigma^K}E[\norm{\lambda^* - \lambda_{T^K}^K}^2]\nonumber \\
    & + \tsum_{k=1}^{K}(\tfrac{1}{\eta^k} + \tfrac{M_g}{2} - \tfrac{1 - \gamma_\lambda^k}{2\sigma^k})\tsum_{t=0}^{T^k - 1}E[\norm{\lambda_{t+1}^k - \lambda_t^k}^2] -  \tsum_{k=1}^{K}\tfrac{\gamma_\lambda^k}{2\sigma^k}\tsum_{t=0}^{T^k-1}E[D_\lambda^k(\lambda_{t+1}^k)] \nonumber\\
    & + \tfrac{\gamma_\lambda^1T^1}{2\sigma^1}E[D_\lambda^1(\lambda^*)] + \tsum_{k=1}^{K-1}\big(\tfrac{\gamma_\lambda^{k+1}T^{k+1}}{2\sigma^{k+1}} - \tfrac{\gamma_\lambda^kT^k}{2\sigma^k}\big)E[D_\lambda^{k+1}(\lambda^*)] - \tfrac{\gamma_\lambda^KT^K}{2\sigma^K}E[D_\lambda^{K+1}(\lambda^*)]  \nonumber \\
    &  + \tsum_{k=1}^{K}\tsum_{t=0}^{T^k - 1}\big((\tfrac{\norm{\lambda^*}L_g + M_g + L_f}{2} + \tfrac{1}{\eta^k} + 2\eta^k L_f^2 - \tfrac{1 - \gamma_x^k}{2\tau^k})E[\norm{x_t^k - x_{t+1}^k}^2] + (2\eta^k L_f^2 - \tfrac{\gamma_x^k}{2\tau^k})E[D_x^k(x_{t+1}^k)] \big) \nonumber \\
    &  + \tfrac{1}{\eta^1}E[\norm{x^* - u_0^1}^2] + \tsum_{k=1}^{K-1}(\tfrac{1}{\eta^{k+1}} - \tfrac{1}{\eta^k})E[\norm{x^* - u_{T^k}^k}^2] - \tfrac{1}{\eta^K}E[\norm{x^* - u_{T^K}^K}^2] + \tfrac{M_g^2\sigma^K}{1-\gamma_{\lambda}^K}E[\norm{x_{T^K-1}^K - x_{T^K}^K}^2] \nonumber\\
    & + \tfrac{1 - \gamma_x^1}{2\tau^1}E[\norm{x^* - x_0^1}^2] + \tsum_{k=1}^{K-1}(\tfrac{1 - \gamma_x^{k+1}}{2\tau^{k+1}} - \tfrac{1 - \gamma_x^k}{2\tau^k})E[\norm{x^* - x_{T^k}^k}^2] - \tfrac{1 - \gamma_x^K}{2\tau^K}E[\norm{x^* - x_{T^K}^K}^2] \nonumber \\
    & + \tfrac{\gamma_x^1T^1}{2\tau^1}E[D_x^1(x^*)] + \tsum_{k=1}^{K-1}\big(\tfrac{\gamma_x^{k+1}T^{k+1}}{2\tau^{k+1}} - \tfrac{\gamma_x^kT^k}{2\tau^k}\big)E[D_x^{k+1}(x^*)] - \tfrac{\gamma_x^KT^K}{2\tau^K}E[D_x^{K+1}(x^*)] \bigg). \label{eq:a27}
\end{align}
\begin{remark} \label{eq:th3}
The following step-size policy holds for Algorithm \ref{alg:alg2} for the semi-finite case to achieve acceleration via the gap inequality \eqref{eq:a27}:
\begin{align}
    \tfrac{1}{\eta^k} + \tfrac{M_g}{2} - \tfrac{1 - \gamma_\lambda^k}{2\sigma^k} &\leq 0 
    &\tfrac{\norm{\lambda^*}L_g + M_g + L_f}{2} + \tfrac{1}{\eta^k} + 2\eta^k L_f^2 - \tfrac{1 - \gamma_x^k}{2\tau^k} \leq 0, \label{eq:a28}
\end{align}
\begin{align}
    2\eta^k L_f^2 - \tfrac{\gamma_x^k}{2\tau^k} &\leq 0 
    &\tfrac{\sigma^KM_g^2}{1-\gamma_{\lambda}^K} - \tfrac{1 - \gamma_x^K}{2\tau^K} & \leq 0, \label{eq:a32}
\end{align}
\end{remark}   
\subsubsection*{Constant step size policy.}
Under the constant step-size policy, $\{\eta^k,\sigma^k,\tau^k,\gamma_{\lambda}^k,\gamma_x^k,T^k\} = \{\eta,\sigma,\tau,\gamma_{\lambda},\gamma_x,T\}$. 
Applying conditions in Remark \ref{eq:th3} on \eqref{eq:a27} 
\begin{align}
    &E[Q(\wt{Z}^K,z^*)]  \leq  \tfrac{1}{S}\big(\tfrac{1 - \gamma_\lambda}{2\sigma}E[\norm{\lambda^* - \lambda_0^1}^2]  + \tfrac{\gamma_\lambda}{2\sigma}TE[D_\lambda^1(\lambda^*)] \nonumber\\
    & + \tfrac{1}{\eta}E[\norm{x^* - u_0^1}^2] + \tfrac{1 - \gamma_x}{2\tau}E[\norm{x^* - x_0^1}^2]  + \tfrac{\gamma_x}{2\tau}TE[D_x^1(x^*)] \big). \label{eq:d1}
\end{align}
\begin{theorem}\label{eq:th7}
Assume that the parameters are chosen as follows:
\begin{align}
    \eta &= \tfrac{1}{m^{1/2}} &\tau&= \tfrac{1}{4m^{1/2} + 2c} &\sigma &=\tfrac{1}{4(m^{1/2} + M_g)}
   & \gamma_x &= \tfrac{1}{m} &\gamma_\lambda &= \tfrac{1}{m} &T=m, \label{eq:a306}
\end{align}
where  $c = \tfrac{2\norm{\lambda^*}L_g + 2M_g + 2L_f + 4L_f^2 + 2M_g^2 + 2}{2}$. Under the parameter choice \eqref{eq:a306}, the expected gap inequality \eqref{eq:a27} will be:
\begin{align}
    & E[Q(\wt{Z}^K,z^*)] \leq  \tfrac{1}{mK}\bigg(\tfrac{2(M_g + m^{1/2})(m-1)}{m}E[\norm{\lambda^* - \lambda_0^1}^2] + 2(m^{1/2} + M_g)E[D_\lambda^1(\lambda^*)] + m^{1/2}E[\norm{x^* - u_0^1}^2]  \nonumber\\
    & + \tfrac{(2m^{1/2} + c)(m-1)}{m}E[\norm{x^* - x_0^1}^2] + (2m^{1/2} + c)E[D_x^1(x^*)] \bigg) = \epsilon. \label{eq:a305}
\end{align}
The number of gradient computations required to achieve an $\epsilon$-gap is $O(\frac{m^{1/2}}{\epsilon})$. 
\end{theorem}
\begin{proof}
    It can be verified directly that the parameters in \eqref{eq:a306} satisfy the conditions of Remark \ref{eq:th3}, and therefore inequality \eqref{eq:a305} holds. Also \eqref{eq:p15} will cancel using \eqref{eq:a28}. From the inequality \eqref{eq:a305}, $K = O(\tfrac{1}{\epsilon m^{1/2}})$. The number of gradient computations is:
    %\begin{align*}
     $   mK+\tsum_{k=1}^{K}T^k=2Km=O(\tfrac{m^{1/2}}{\epsilon})
    $.%\end{align*}
\end{proof}
\begin{remark}\label{rem:manage_lambda_in_step_size}
    In Theorem \ref{eq:th7}, for the parameter $\tau$, the expression $\norm{\lambda^*}$ appears which is unknown and finite. For a variance-reduced method to be effective, we must have $m$ and/or $n$ be large. If these are large, e.g., $\sqrt{m} \gg L_g\norm{\lambda^*}$ then the driving factor deciding the order of $\tau$ is the $\sqrt{m}$ term in the denominator. Hence, a large $m$ (or $n$ as we  will see in fully-finite-sum case later) would not only allows us to see the impact of variance reduction but also assist us to ignore the $L_g\norm{\lambda^*}$ term in $\tau$. It is also not hard to see that $m$ and $n$ do not affect $\lambda^*$ in the limiting sense. In particular, one can think of $f$ as a sample average approximation of an unknown stochastic function $\wt{f}(x) = \Ebb_\xi[F(x, \xi)]$ and each $f_i(x) = F(x, \xi_i)$ where $\{\xi_i\}_{i=1}^m$ are i.i.d. samples of $\xi$. As $m$ increases, we see that the $f(x) \approx \wt{f}(x)$ which does not cause increase of $\lambda^*$ as a function of $m$. Hence, a large value of $m$ (or $n$ in fully-finite-sum case later) helps in mitigating the impact of $\lambda^*$ on the step-size. On the other hand, if both $m$ and $n$ are small, we may as well treat the problem as deterministic. Henceforth, throughout the paper we treat $m$ or $n$ large enough to balance $L_g\norm{\lambda^*}$. Finally, we set $\tau = O(\tfrac{1}{\sqrt{m}})$ in our numerical experiments, and do not observe any impact of $\lambda^*$ on the convergence behavior.%    This will be true for both constant and non-constant step-size policies and in the latter case, $k$ increases with the epochs changing.  
\end{remark}

\subsubsection*{Non-constant step-size policy.}
The following approach provides a constant upper bound for certain Lagrangian-related terms.
\begin{lemma} (Discrete Grönwall's inequality) \label{eq:a35}
    Assume that $\{a_i\}_{i \in \mathrm{N}}$ is a non-negative sequence, and $m_1,m_2 \geq 0$ are constants such that $a_0 \leq m_1$ and $a_t \leq m_1 + m_2\sum_{i=0}^{t-1}a_i$, then $a_t \leq m_1(1+m_2)^t$.
\end{lemma}
\begin{proof}
    With induction on $t$, the statement is true for $t=0$ since $a_0 \leq m_1$. Assume that the statement is also true for an arbitrary $t \geq 1$. Then we need to prove that it is true for $t+1$:
    \begin{align*}
        a_{t+1} &\leq m_1 + m_2\tsum_{i=0}^{t}a_i \leq m_1\big[1 + m_2\tsum_{i=0}^{t}(1+m_2)^i\big]  \leq m_1\big[1 + m_2\tfrac{(1+m_2)^{t+1} - 1}{m_2}\big] = m_1(1+m_2)^{t+1}.
    \end{align*}
\end{proof} 
\begin{lemma} \label{eq:l6}
    Assume that the parameters for Algorithm \ref{alg:alg2} are as follows:
    \begin{align*}
        \eta^k &= \tfrac{1}{k+1} &\tau^k &= \tfrac{1}{4k + 1 + 2c} & \sigma^k &= \tfrac{1}{2(k+1)(M_g + 1)} &
        \gamma_x^k &= \tfrac{1}{(k+1)^2} & \gamma_\lambda^k &= \tfrac{1}{(k+1)^2} &T^k = (k+1)^2.
    \end{align*}
    Lagrangian terms $E[\norm{\lambda - \lambda_{T^K}^K}^2] + E[D_\lambda^{K+1}(\lambda)]$ will have the following constant upper bound:
    \begin{align}
    E[\norm{\lambda - \lambda_{T^K}^K}^2] + E[D_\lambda^{K+1}(\lambda)] \leq \sup_{K \geq 1} \tfrac{ (M_g+1)( E[\norm{\lambda^*-\lambda_0^1}^2] + E[D_\lambda^1(\lambda^*)]) + 5KD_x^2 } {\tfrac{(K+1)(M_g+1)}{2}( 1-\tfrac{1}{(K+1)^2}) } e^4 = D_{\lambda}^2.\label{eq:d5}
\end{align}
\end{lemma}
\begin{proof}
Applying conditions in Remark \ref{eq:th3}, multiplying both sides of \eqref{eq:a27} by $S^K$ and since $E[Q(\wt{Z}^K,z)] \geq 0$, transferring $ - \tfrac{1 - \gamma_\lambda^K}{4\sigma^K}E[\norm{\lambda - \lambda_{T^K}^K}^2]$ and $ - \tfrac{\gamma_\lambda^{K+1}T^{K+1}}{2\sigma^{K+1}}E[D_\lambda^{K+1}(\lambda)]$ to the left-hand side of \eqref{eq:a27} and $M_1 = R - \tsum_{k=1}^{K-1}(\tfrac{\gamma_\lambda^{k+1}T^{k+1}}{2\sigma^{k+1}} - \tfrac{\gamma_\lambda^kT^k}{2\sigma^k})E[D_\lambda^{k+1}(\lambda^*)] - \tsum_{k=1}^{K-1}(\tfrac{1 - \gamma_\lambda^{k+1}}{2\sigma^{k+1}} - \tfrac{1 - \gamma_\lambda^k}{2\sigma^k})E[\norm{\lambda - \lambda_{T^k}^k}^2]$ and 
\begin{align*}
    M_2 &=\min\{\tfrac{1 - \gamma_\lambda^K}{4\sigma^K}, \tfrac{\gamma_\lambda^{K+1}T^{K+1}}{2\sigma^{K+1}}\}& , & & M_3^k = \max \{\tfrac{1 - \gamma_\lambda^{k+1}}{2\sigma^{k+1}} - \tfrac{1 - \gamma_\lambda^k}{2\sigma^k},\tfrac{\gamma_\lambda^{k+1}T^{k+1}}{2\sigma^{k+1}} - \tfrac{\gamma_\lambda^kT^k}{2\sigma^k}\},
\end{align*}
and $M_3=\max_{k \in \{1,...,K-1\}}\{M_3^k\}$, and $R$ is the right-hand side of \eqref{eq:a27}.
\begin{align}
    & M_2\big(E[\norm{\lambda - \lambda_{T^K}^K}^2] + E[D_\lambda^{K+1}(\lambda)]\big) \leq M_1 + \tsum_{k=1}^{K-1}M_3^k\big(E[D_\lambda^{k+1}(\lambda)] + E[\norm{\lambda - \lambda_{T^k}^k}^2]\big) \nonumber\\
    & \leq M_1 + M_3\tsum_{k=1}^{K-1}\big(E[D_\lambda^{k+1}(\lambda)] + E[\norm{\lambda - \lambda_{T^k}^k}^2]\big). \label{eq:a36} 
\end{align}
Using Lemma \ref{eq:a35} for \eqref{eq:a36}, $E[\norm{\lambda - \lambda_{T^K}^K}^2] + E[D_\lambda^{K+1}(\lambda)] \leq \tfrac{M_1}{M_2}(1+\tfrac{M_3}{M_2})^K$, and with the parameter setup in the lemma we obtain \eqref{eq:d5}. 
\end{proof}
\begin{theorem}
    The expected gap inequality for semi-finite case, with the parameter setup in Lemma \ref{eq:l6} will be: 
    \begin{align}
        E[Q(\wt{Z}^K,z^*)] \leq \tfrac{ 6\Big( C_0 + 5KD_x^2 + 2KD_\lambda^2 \Big) } {K(K+1)(2K+1)} = \epsilon, \label{eq:a407}
    \end{align}
    where $C_0 = (M_g+1)(\tfrac{3}{2}E[\norm{\lambda^*-\lambda_0^1}^2] + 2E[D_\lambda^1(\lambda^*)]) + cD_x^2$. The number of gradient computations of the objective function will be $O(\tfrac{m}{\sqrt{\epsilon}}+\tfrac{1}{\epsilon\sqrt{\epsilon}})$ to reach $\epsilon$-gap.
\end{theorem}
\begin{proof}
It can be verified directly that the gap inequality for non-constant step size with the parameters in Lemma \ref{eq:l6} will be \eqref{eq:a407} and $K = O(\tfrac{1}{\epsilon^{1/2}})$ in $\epsilon$-gap. 
The gradient computations to reach $\epsilon$-gap is $mK+\tsum_{k=1}^{K}T^k=mK+S^K=O(\tfrac{m}{\sqrt{\epsilon}}+\tfrac{1}{\epsilon\sqrt{\epsilon}})$ for the objective function.
 \end{proof}
 Similar to Remark \ref{rem:manage_lambda_in_step_size}, we also see that dominating term of $\tau$ becomes $k$ whenever $k \ge L_g\norm{\lambda^*}$. This is a small number compared to the large majority of $k$-values that get as large as $O(\tfrac{1}{\sqrt{\epsilon}})$.
\subsubsection{Fully-finite-sum case.} \label{sec-method2.2.2}
\begin{lemma} \label{eq:a38} The following inequalities hold:
    \begin{align}
        &E[\norm{\delta_t^{\lambda^k}}^2] \leq 2M_g^2(E[D_x^k(x_{t+1}^k)]  + E[\norm{x_{t+1}^k - x_{t}^k}^2]), \label{eq:a39}\\
        &E[\norm{q_t^k - \wb{q}_t^k}^2] \leq M_g^2E[\norm{x_t^k - x_{t-1}^k}^2]. \label{eq:a41}
    \end{align}
    \begin{proof}
        For $E[\norm{\delta_t^{\lambda^k}}^2]$:
        \begin{align*}
            &E[\norm{\delta_t^{\lambda^k}}^2] = E[\norm{g_{j_t^k}(x^k_t) - g_{j_t^k}(\wt{x}^{k-1}) - [g(x_t^k) - g(\wt{x}^{k-1})]}^2] \leq 2M_g^2E[D_x^k(x_{t+1}^k)]  + 2M_g^2E[\norm{x_{t+1}^k - x_{t}^k}^2],
        \end{align*}
        where the  inequality comes from $Var(X) \leq E[X^2]$, triangle inequality, \eqref{eq:p1}, and \eqref{eq:p7}. For $E[\norm{q_t^k - \wb{q}_t^k}^2]$:
        \begin{align*}
           & E[\norm{q_t^k - \bar{q}_t^k}^2] =  E[\norm{\ell_{g_{j_t^k}}(x_t^k) - g_{j_t^k}(x_{t-1}^k) - [\ell_g(x^k_t) - g(x^k_{t-1})]}^2] \leq M_g^2E[\norm{x_t^k - x_{t-1}^k}^2],
        \end{align*}
        where we use $Var(X) \leq E[X^2]$ and \eqref{eq:p1}.
    \end{proof}
\end{lemma}
\begin{remark}
    For bounding the Lagrangian terms the following conditions hold for the gap inequality \eqref{eq:b27}:
    \begin{align*}
        \tfrac{M_g\eta^k}{2} - \tfrac{1 - \gamma_\lambda^k}{2\sigma^k} & \leq 0 & M_g\eta^k - \tfrac{\gamma_\lambda^k}{2\sigma^k} & \leq 0 & \tfrac{M_g^2\sigma^K}{1-\gamma_{\lambda}^K} - \tfrac{\eta^KM_g^2}{2} - \tfrac{M_g}{2\eta^K} & \leq 0
    \end{align*}
    \begin{align*}
        \tfrac{\norm{\lambda^*}L_g + 5L_f}{2} + \tfrac{M_g + L_gD_x + 2L_f}{\eta^k} + 2M_g^2\eta^k + \tfrac{\eta^kM_g^2}{2} + \tfrac{M_g}{2\eta^k} - \tfrac{1 - \gamma_x^k}{2\tau^k} & \leq 0 & (2L_f + 2M_g^2)\eta^k - \tfrac{\gamma_x^k}{2\tau^k} & \leq 0
    \end{align*}
\end{remark}
\begin{lemma}
    Lagrangian term $E[\norm{\lambda^* - v_{T^k}^k}^2] + E[\norm{\lambda^* - w_{T^k}^k}^2] + E[\norm{\lambda^* - \lambda_{T^k}^k}^2] + E[D_\lambda^{k+1}(\lambda^*)]$ has the constant upper bound
    \begin{align*}
        D_{\lambda}^2 := \sup_{K \geq 1} \tfrac{R^K}{U(K)}e^{5},
    \end{align*}
    where $R^K = 2c^{\prime}E[\norm{\lambda^* - v_0^1}^2] + 2c^{\prime}E[\norm{\lambda^* - w_0^1}^2] + (2c^{\prime} + 1)E[D_\lambda^1(\lambda^*)] + (2c' + c_x + \tfrac{M_gc'}{2} + \tfrac{6c'+c_x}{6})KD_x^2 + K\norm{\lambda^*}^2$, $U(K) = \min\{(K+1)c',\tfrac{c^{\prime}K}{2},c^{\prime}(K+1)\}$, $c_x = c'(\tfrac{2M_g + 2L_gD_x + 4L_f +M_g}{2}) + \norm{\lambda^*}L_g + 5L_f + 5M_g^2$ and $c^{\prime} =2 + M_g + 2L_gD_x$. 
\end{lemma} 
\begin{proof}
For $\inprod{\delta_t^{x^k}}{x_t^k - x_{t+1}^k}$:
\begin{align*}
    & \inprod{\delta_t^{x^k}}{x_t^k - x_{t+1}^k} = \inprod{\grad f_{i_t^k}(x^k_t) - \grad f_{i_t^k}(\wt{x}^{k-1}) - [\grad f(x^k_t)  -\grad f(\wt{x}^{k-1})]}{x_t^k - x_{t+1}^k} \nonumber\\
    & + \inprod{\lambda_{t+1}^k\grad g_{p_t^k}(x^k_t) - \grad g_{p_t^k}(\wt{x}^{k-1})\wt{\lambda}^{k-1} - [\lambda_{t+1}^k\grad g(x^k_t)  - \grad g(\wt{x}^{k-1})\wt{\lambda}^{k-1}]}{x_t^k - x_{t+1}^k}.
\end{align*}
Moreover, we have
\begin{align*}
    &\inprod{\lambda_{t+1}^k\grad g_{p_t^k}(x^k_t) - \grad g_{p_t^k}(\wt{x}^{k-1})\wt{\lambda}^{k-1} - [\lambda_{t+1}^k\grad g(x^k_t)  - \grad g(\wt{x}^{k-1})\wt{\lambda}^{k-1}]}{x_t^k - x_{t+1}^k} \nonumber\\
    & \leq 2M_g\norm{\lambda_{t+1}^k - \wt{\lambda}^{k-1}}\norm{x_t^k - x_{t+1}^k} + 2L_g\norm{\wt{\lambda}^{k-1}}\norm{x_t^k - \wt{x}^{k-1}}\norm{x_t^k - x_{t+1}^k} \nonumber\\
    & \leq 2M_g\norm{\lambda_{t+1}^k - \wt{\lambda}^{k-1}}\norm{x_t^k - x_{t+1}^k} + 4L_g(\norm{\wt{\lambda}^{k-1} - \lambda^*} + \norm{\lambda^*})\norm{x_t^k - \wt{x}^{k-1}}\norm{x_t^k - x_{t+1}^k} \nonumber\\
    & \leq 2M_g\norm{\lambda_{t+1}^k - \wt{\lambda}^{k-1}}\norm{x_t^k - x_{t+1}^k} + 4L_gD_x\norm{\wt{\lambda}^{k-1} - \lambda^*}\norm{x_t^k - x_{t+1}^k} + 4L_g\norm{\lambda^*}\norm{x_t^k - \wt{x}^{k-1}}\norm{x_t^k - x_{t+1}^k} \nonumber\\
    & \leq M_g\eta^kD_{\lambda}^k(\lambda_{t+1}^k) + 2L_gD_x\eta^kD_{\lambda}^k(\lambda^*) + 2L_gD_x\eta^k\norm{\lambda^*}^2 + \tfrac{M_g + L_gD_x }{\eta^k}\norm{x_t^k - x_{t+1}^k}^2 ,
\end{align*}
where the first inequality comes from  \eqref{eq:p9}, \eqref{eq:p2}, \eqref{eq:p3} and triangle inequality, and the second inequality comes from \eqref{eq:p10}, and \eqref{eq:p7}. Also
\begin{align*}
    &\inprod{\grad f_{i_t^k}(x^k_t) - \grad f_{i_t^k}(\wt{x}^{k-1}) - [\grad f(x^k_t)  -\grad f(\wt{x}^{k-1})]}{x_t^k - x_{t+1}^k}  \leq \big(2L_f + \tfrac{L_f}{\eta^k}\big)\norm{x_t^k - x_{t+1}^k}^2 + L_f\eta^kD_x^k(x_{t+1}^k),
\end{align*}
where the inequality comes from \eqref{eq:p9}, \eqref{eq:p10}, \eqref{eq:p2}, and \eqref{eq:p7}. The upper bound for $\inprod{\delta_t^{x^k}}{x_t^k - x_{t+1}^k}$
\begin{align}
     &\inprod{\delta_t^{x^k}}{x_t^k - x_{t+1}^k} \leq M_g\eta^kD_{\lambda}^k(\lambda_{t+1}^k) + 2L_gD_x\eta^kD_{\lambda}^k(\lambda^*) + 2L_gD_x\eta^k\norm{\lambda^*}^2\nonumber\\ 
     &  + (\tfrac{M_g + L_gD_x + L_f}{\eta^k} + 2L_f)\norm{x_t^k - x_{t+1}^k}^2 + L_f\eta^kD_x^k(x_{t+1}^k). \label{eq:b24}
\end{align}
We apply\eqref{eq:b24} to \eqref{eq:1}, take expectation on both sides of the inequality and rearrange the terms
\begin{align*}
    & E[Q(z_{t+1}^k, z^*)] \leq \tfrac{1}{\eta^k} (E[\norm{\lambda^* - v_t^k}^2] - E[\norm{\lambda^* - v_{t+1}^k}^2]) + E[\inprod{\delta_t^{\lambda^k}}{\lambda_t^k - v_t^k}] - E[\inprod{\bar{q}_{t+1}^k}{\lambda_{t+1}^k - \lambda^*}] \nonumber\\
    & + E[\inprod{\bar{q}_t^k}{\lambda_t^k - \lambda^*}]  + \tfrac{1}{\eta^k} (E[\norm{\lambda^* - w_t^k}^2] - E[\norm{\lambda^* - w_{t+1}^k}^2]) + (\tfrac{\eta^kM_g^2}{2} + \tfrac{M_g}{2\eta^k})E[\norm{x_t^k - x_{t-1}^k}^2] \nonumber\\
    & + E[\inprod{q_t^k - \bar{q}_t^k}{\lambda_t^k - w_t^k}] + (\tfrac{M_g\eta^k}{2} + \tfrac{1}{\eta^k} - \tfrac{1 - \gamma_\lambda^k}{2\sigma^k})E[\norm{\lambda_t^k - \lambda_{t+1}^k}^2] + \tfrac{1 - \gamma_\lambda^k}{2\sigma^k} \big(E[\norm{\lambda^* - \lambda_t^k}^2] - E[\norm{\lambda^* - \lambda_{t+1}^k}^2] \big) \nonumber \\
    & + \tfrac{\gamma_\lambda^k}{2\sigma^k} \big(E[D_\lambda^k(\lambda^*)] 
    - E[\norm{\lambda^* - \lambda_{t+1}^k}^2]\big) + (M_g\eta^k - \tfrac{\gamma_\lambda^k}{2\sigma^k})E[D_\lambda^k(\lambda_{t+1}^k)] + 2L_gD_x\eta^kE[D_{\lambda}^k(\lambda^*)] + 2L_gD_x\eta^k\norm{\lambda^*}^2\nonumber\\ 
    & + \tfrac{1 - \gamma_x^k}{2\tau^k} \big( E[\norm{x^* - x_t^k}^2] 
    - E[\norm{x^* - x_{t+1}^k}^2] \big) + \tfrac{\gamma_x^k}{2\tau^k} \big(E[D_x^k(x^*)]- E[\norm{x^* - x_{t+1}^k}^2]\big) \nonumber\\
    & + (\tfrac{\norm{\lambda^*}L_g + 5L_f}{2} + \tfrac{M_g + L_gD_x + L_f}{\eta^k} + 2M_g^2\eta^k - \tfrac{1 - \gamma_x^k}{2\tau^k})E[\norm{x_{t+1}^k - x_t^k}^2] + (L_f\eta^k + 2M_g^2\eta^k - \tfrac{\gamma_x^k}{2\tau^k})E[D_x^k(x_{t+1}^k)]. 
\end{align*}
Summing up the last inequality with respect to index $t$ and \eqref{eq:p8}
\begin{align}
    & T^k E[Q(\wt{z}^k,z)] \leq \tfrac{1}{\eta^k} (E[\norm{\lambda^* - v_0^k}^2] - E[\norm{\lambda^* - v_{T^k}^k}^2]) + \tsum_{t=0}^{T^k - 1}E[\inprod{\delta_t^{\lambda^k}}{\lambda_t^k - v_t^k}] \nonumber\\
    & + E[\inprod{\bar{q}_0^k}{\lambda_0^k - \lambda^*}] - E[\inprod{\bar{q}_{T^k}^k}{\lambda_{T^k}^k - \lambda^*}] + \tfrac{1}{\eta^k} (E[\norm{\lambda^* - w_0^k}^2] - E[\norm{\lambda^* - w_{T^k}^k}^2]) \nonumber\\
    & + \tsum_{t=0}^{T^k - 1}(\tfrac{M_g\eta^k}{2} - \tfrac{1 - \gamma_\lambda^k}{2\sigma^k})E[\norm{\lambda_t^k - \lambda_{t+1}^k}^2] + \tfrac{1 - \gamma_\lambda^k}{2\sigma^k} \big(E[\norm{\lambda^* - \lambda_0^k}^2] - E[\norm{\lambda^* - \lambda_{T^k}^k}^2] \big) \nonumber \\
    & + \tfrac{\gamma_\lambda^k T^k}{2\sigma^k}E[D_\lambda^k(\lambda^*)] 
    - \tfrac{\gamma_\lambda^k T^k}{2\sigma^k}E[D_\lambda^{k+1}(\lambda^*)] + \tsum_{t=0}^{T^k - 1}(M_g\eta^k - \tfrac{\gamma_\lambda^k}{2\sigma^k})E[D_\lambda^k(\lambda_{t+1}^k)]\nonumber \\
    & +  2L_gD_x\eta^kT^kE[D_{\lambda}^k(\lambda^*)] + 2L_gD_x\eta^kT^k\norm{\lambda^*}^2 + \tfrac{1 - \gamma_x^k}{2\tau^k} \big( E[\norm{x^* - x_0^k}^2] - E[\norm{x^* - x_{T^k}^k}^2] \big) \nonumber\\
    & + \tsum_{t=0}^{T^k - 1}(\tfrac{\norm{\lambda^*}L_g + 5L_f}{2} + \tfrac{M_g + L_gD_x + L_f}{\eta^k} + 2M_g^2\eta^k + \tfrac{\eta^kM_g^2}{2} + \tfrac{M_g}{2\eta^k} - \tfrac{1 - \gamma_x^k}{2\tau^k})E[\norm{x_{t+1}^k - x_t^k}^2] \nonumber \\
    & + (\tfrac{\eta^kM_g^2}{2} + \tfrac{M_g}{2\eta^k})E[\norm{x_0^k - x_{-1}^k}^2] - (\tfrac{\eta^kM_g^2}{2} + \tfrac{M_g}{2\eta^k})E[\norm{x_{T^k-1}^k - x_{T^k}^k}^2] \nonumber\\
    & + \tfrac{\gamma_x^k T^k}{2\tau^k}E[D_x^k(x^*)]- \tfrac{\gamma_x^k T^k}{2\tau^k}E[D_x^{k+1}(x^*)] + \tsum_{t=0}^{T^k - 1}(L_f\eta^k + 2M_g^2\eta^k - \tfrac{\gamma_x^k}{2\tau^k})E[D_x^k(x_{t+1}^k)].\label{eq:b26}
\end{align}
Summing up inequality \eqref{eq:b26} with respect to index $k$, and \eqref{eq:p8}, and since $(x_0^{K+1},\lambda_0^{k+1})= (x_{T^k}^k,\lambda_{T^k}^k)$:
\begin{align}
    & E[Q(\wt{Z}^K,z^*)] \leq \tfrac{1}{S^K}\bigg(\tfrac{1}{\eta^1}E[\norm{\lambda^* - v_0^1}^2] + \tsum_{k=1}^{K-1}(\tfrac{1}{\eta^{k+1}} - \tfrac{1}{\eta^k})E[\norm{\lambda^* - v_{T^k}^k}^2] - \tfrac{1}{\eta^K}E[\norm{\lambda^* - v_{T^K}^K}^2]\nonumber\\
    &  + \tsum_{k=1}^{K-1}\tsum_{t=0}^{T^k - 1}(\tfrac{M_g\eta^k}{2} - \tfrac{1 - \gamma_\lambda^k}{2\sigma^k})E[\norm{\lambda_t^k - \lambda_{t+1}^k}^2] \nonumber\\
    & + \tfrac{1}{\eta^1}E[\norm{\lambda^* - w_0^1}^2] + \tsum_{k=1}^{K-1}(\tfrac{1}{\eta^{k+1}} - \tfrac{1}{\eta^k})E[\norm{\lambda^* - w_{T^k}^k}^2] - \tfrac{1}{\eta^K}E[\norm{\lambda^* - w_{T^K}^K}^2] \nonumber\\
    & + \tfrac{1 - \gamma_\lambda^1}{2\sigma^1}E[\norm{\lambda^* - \lambda_0^1}^2] + \tsum_{k=1}^{K-1}(\tfrac{1 - \gamma_\lambda^{k+1}}{2\sigma^{k+1}} - \tfrac{1 - \gamma_\lambda^k}{2\sigma^k})E[\norm{\lambda^* - \lambda_{T^k}^k}^2] - \tfrac{1 - \gamma_\lambda^K}{4\sigma^K}E[\norm{\lambda^* - \lambda_{T^K}^K}^2] \nonumber \\
    & + (\tfrac{\gamma_\lambda^1 T^1}{2\sigma^1} + 2L_gD_x\eta^1T^1)E[D_\lambda^1(\lambda^*)] + \tsum_{k=1}^{K-1}(\tfrac{\gamma_\lambda^{k+1} T^{k+1}}{2\sigma^{k+1}} - \tfrac{\gamma_\lambda^k T^k}{2\sigma^k} + 2L_gD_x\eta^kT^k)E[D_\lambda^{k+1}(\lambda^*)] \nonumber\\
    & - \tfrac{\gamma_\lambda^K T^K}{2\sigma^K}E[D_\lambda^{K+1}(\lambda^*)] + \tsum_{k=1}^{K}\tsum_{t=0}^{T^k - 1}(M_g\eta^k - \tfrac{\gamma_\lambda^k}{2\sigma^k})E[D_\lambda^k(\lambda_{t+1}^k)] + \tsum_{k=1}^{K}2L_gD_x\eta^kT^k\norm{\lambda^*}^2\nonumber \\
    & + \tfrac{1 - \gamma_x^1}{2\tau^1}E[\norm{x^* - x_0^1}^2] + \tsum_{k=1}^{K-1}(\tfrac{1 - \gamma_x^{k+1}}{2\tau^{k+1}} - \tfrac{1 - \gamma_x^k}{2\tau^k})E[\norm{x^* - x_{T^k}^k}^2] - \tfrac{1 - \gamma_x^K}{2\tau^K}E[\norm{x^* - x_{T^K}^K}^2] \nonumber\\
    & + \tsum_{k=1}^{K}\tsum_{t=0}^{T^k - 1}(\tfrac{\norm{\lambda^*}L_g + 5L_f}{2} + \tfrac{M_g + L_gD_x + 2L_f}{\eta^k} + 2M_g^2\eta^k + \tfrac{\eta^kM_g^2}{2} + \tfrac{M_g}{2\eta^k} - \tfrac{1 - \gamma_x^k}{2\tau^k})E[\norm{x_{t+1}^k - x_t^k}^2] \nonumber \\
    & + (\tfrac{\eta^1M_g^2}{2} + \tfrac{M_g}{2\eta^1})E[\norm{x_0^1 - x_{-1}^1}^2] + \tsum_{k=1}^{K-1}\big((\tfrac{\eta^{k+1}M_g^2}{2} + \tfrac{M_g}{2\eta^{k+1}}) - (\tfrac{\eta^kM_g^2}{2} + \tfrac{M_g}{2\eta^k})\big)E[\norm{x_{T^k-1}^k - x_{T^k}^k}^2]  \nonumber\\
    & + (\tfrac{M_g^2\sigma^K}{1-\gamma_{\lambda}^K} - \tfrac{\eta^KM_g^2}{2} - \tfrac{M_g}{2\eta^K})E[\norm{x_{T^K-1}^K - x_{T^K}^K}^2] + \tsum_{k=1}^{K}\tsum_{t=0}^{T^k - 1}((2L_f + 2M_g^2)\eta^k - \tfrac{\gamma_x^k}{2\tau^k})E[D_x^k(x_{t+1}^k)] \nonumber\\
    & + \tfrac{\gamma_x^1 T^1}{2\tau^1}E[D_x^1(x^*)] + \tsum_{k=1}^{K-1} (\tfrac{\gamma_x^{k+1} T^{k+1}}{2\tau^{k+1}}- \tfrac{\gamma_x^k T^k}{2\tau^k})E[D_x^{k+1}(x^*)] - \tfrac{\gamma_x^K T^K}{2\tau^K}E[D_x^{K+1}(x^*)] \bigg).\label{eq:b27}
\end{align}
The following parameter setup for \eqref{eq:b27}.
\begin{align}
    \eta^k &= \tfrac{1}{c'(k+1)} & \gamma_{\lambda}^k &= \tfrac{1}{k+1} & \sigma^k &= \tfrac{1}{2c^{\prime}(k+1)} &\gamma_x^k &= \tfrac{1}{k+1} &\tau^k &= \tfrac{1}{2c^{\prime}(k+1) + 2c_x} & T^k &= k+1 \label{eq:b28}
\end{align}
Notice that Remark \ref{rem:manage_lambda_in_step_size} holds true for this parameters setup. The new expected gap inequality after applying \eqref{eq:b28}, multiplying both sides by $S^K$, and $\tfrac{2L_gD_x}{c'} < 1$ 
\begin{align}
    & S^K E[Q(\wt{Z}^K,z^*)] \leq 2c'E[\norm{\lambda^* - v_0^1}^2] + \tsum_{k=1}^{K-1}c'E[\norm{\lambda^* - v_{T^k}^k}^2] - (K+1)c'E[\norm{\lambda^* - v_{T^K}^K}^2]\nonumber\\
    & + 2c'E[\norm{\lambda^* - w_0^1}^2] + \tsum_{k=1}^{K-1}c'E[\norm{\lambda^* - w_{T^k}^k}^2] - (K+1)c'E[\norm{\lambda^* - w_{T^K}^K}^2] \nonumber\\
    & + c'E[\norm{\lambda^* - \lambda_0^1}^2]+ \tsum_{k=1}^{K-1}c^{\prime}E[\norm{\lambda^* - \lambda_{T^k}^k}^2] - \tfrac{c^{\prime}K}{2}E[\norm{\lambda^* - \lambda_{T^K}^K}^2] \nonumber \\
    & + (2c^{\prime} + 1)E[D_\lambda^1(\lambda^*)] + \tsum_{k=1}^{K-1}(2c^{\prime} + 1)E[D_\lambda^{k+1}(\lambda^*)] - c^{\prime}(K+1)E[D_\lambda^{K+1}(\lambda^*)]\nonumber\\
    & + K\norm{\lambda^*}^2 + \tfrac{2c'+c_x}{2}E[\norm{x^*-x_0^1}^2]  + \tsum_{k=1}^{K-1}\tfrac{6c'+c_x}{6}E[\norm{x^* - x_{T^k}^k}^2]\nonumber \\
    & + (\tfrac{M_g^2}{4c'} + M_gc')E[\norm{x_0^1 - x_{-1}^1}^2] + \tsum_{k=1}^{K-1}\tfrac{M_gc'}{2}E[\norm{x_{T^k-1}^k - x_{T^k}^k}^2]\nonumber\\
    & + (2c'+c_x)E[D_x^1(x^*)] + \tsum_{k=1}^{K-1} (2c'+c_x)E[D_x^{k+1}(x^*)].\label{eq:b29}
\end{align}
If $U(K) = \min\{(K+1)c',\tfrac{c^{\prime}K}{2},c^{\prime}(K+1)\}$ then
\begin{align}
    & E[\norm{\lambda^* - v_{T^K}^K}^2] + E[\norm{\lambda^* - w_{T^K}^K}^2] + E[\norm{\lambda^* - \lambda_{T^K}^K}^2] + E[D_\lambda^{K+1}(\lambda^*)] \leq \tfrac{1}{U(K)}\bigg( 2c^{\prime}E[\norm{\lambda^* - v_0^1}^2] \nonumber\\
    &+ 2c^{\prime}E[\norm{\lambda^* - w_0^1}^2] + (2M_g + 1)E[D_\lambda^1(\lambda^*)] + \tsum_{k=1}^{K-1}(2c^{\prime} + 1)\big(E[\norm{\lambda^* - v_{T^k}^k}^2]  \nonumber\\
    & + E[\norm{\lambda^* - w_{T^k}^k}^2] + E[\norm{\lambda^* - \lambda_{T^k}^k}^2] + E[D_\lambda^{k+1}(\lambda^*)]\big) + (2c' + c_x + \tfrac{M_gc'}{2} + \tfrac{6c'+c_x}{6})KD_x^2 + K\norm{\lambda^*}^2 \bigg) .\label{eq:b30}
\end{align}
Naming the right-hand-side of \eqref{eq:b30} $R^K$, applying Lemma \ref{eq:a35} to \eqref{eq:b30}%. The upper bound for the expression
\begin{align}
    & E[\norm{\lambda^* - v_{T^K}^K}^2] + E[\norm{\lambda^* - w_{T^K}^K}^2] + E[\norm{\lambda^* - \lambda_{T^K}^K}^2] + E[D_\lambda^{K+1}(\lambda^*)] \nonumber\\
    &\leq \tfrac{R^K - \tsum_{k=1}^{K-1}(2c^{\prime} + 1 )\big(E[\norm{\lambda^* - v_{T^k}^k}^2] + E[\norm{\lambda^* - w_{T^k}^k}^2] + E[\norm{\lambda^* - \lambda_{T^k}^k}^2] + E[D_\lambda^{k+1}(\lambda^*)]\big)}{U(K)}\bigg(1 + \tfrac{2c^{\prime} + 1}{U(K)}\bigg)^K. \label{eq:b31} %what is U(K)
\end{align}
\begin{align*}
    \lim_{K \rightarrow \infty}\bigg(1 + \tfrac{2c^{\prime} + 1}{U(K)}\bigg)^K = \lim_{K \rightarrow \infty} e^{\tfrac{K(2c^{\prime} + 1)}{U(K)}} \leq  e^{5},
\end{align*}
\begin{align*}
    & E[\norm{\lambda^* - v_{T^K}^K}^2] + E[\norm{\lambda^* - w_{T^K}^K}^2] + E[\norm{\lambda^* - \lambda_{T^K}^K}^2] + E[D_\lambda^{K+1}(\lambda^*)]\nonumber\\
    & \leq \lim_{K \rightarrow \infty} \tfrac{R^K - \tsum_{k=1}^{K-1}(2c^{\prime} + 1)\big(E[\norm{\lambda^* - v_{T^k}^k}^2] + E[\norm{\lambda^* - w_{T^k}^k}^2] + E[\norm{\lambda^* - \lambda_{T^k}^k}^2] + E[D_\lambda^{k+1}(\lambda^*)]\big)}{U(K)}e^{5} = D_{\lambda}^2.
\end{align*}
Hence, we conclude the proof. 
\end{proof}
 %\label{rem :new_assumption_for_bounding_the_error}
    Note that when we bound the expected dual iterate by $D_\lambda$, we only obtain a constant variance on the random variable $\delta_t^{x^k}$. This is because variance due to stochastic estimator of $\grad g$ is scaled by factor $\wt{\lambda}^{k}$ which remains bounded in expectation. This is similar to the stochastic approximation algorithm. However, we cannot get accelerated complexity with a constant variance. To achieve variance reduction, we need to bound the expectation of product term $\norm{\wt{\lambda}^{k-1}} \norm{x^k_t-x^k_{t+1}}$. Unfortunately, these two terms have bias that does not seems to have an apparent upper bound. %To achieve accelerated complexity, we must show the variance reduction. 
    We make the following assumption on root mean squared (RMS) quantity $(E[\norm{x^k_t-x^k_{t+1}}^2])^{1/2}$ to get variance reduction for $E[\norm{\delta_t^{x^k}}]$.
    \begin{assumption}
        For the fully-finite-sum case, we assume
    \begin{equation}
        E[\norm{x_t^k - x_{t+1}^k}^2]^{1/2} \leq M E[\norm{x_t^k - x_{t+1}^k}^2] \label{eq:RMS}\tag{RMS-bound}
    \end{equation}
    \end{assumption}
    The \eqref{eq:RMS} assumption is not restrictive. It is well known that for any nonnegative random variable $X$, if we know $E[X] > 0$ then there exists an $M$ that satisfies $E[X^2]^{1/2} \le M E[X^2]$, Since $\norm{x^k_t-x^k_{t+1}}$ is positive a.s., we get existence of $M$. % tells us that $E[\norm{x_t^k - x_{t+1}^k}^2]^{1/2}$ is bounded below by some constant number which is not a strong assumption for our study since the method of study is an inexact method and we do not have a guarantee for $\lim_{k \Rightarrow \infty} E[\norm{x_t^k - x_{t+1}^k}^2] = 0$.

\begin{lemma}
    Suppose the \ref{eq:RMS} holds. Then, we get the following
    \begin{align}
        & E[\inprod{\delta_t^{x^k}}{x - x_{t+1}^k}] \leq M_g \eta^k E[D_{\lambda}^k(\lambda_{t+1}^k)] + \eta^k L_f^2E[D_x^k(x_{t+1}^k)] \nonumber\\
        & + \big(\tfrac{1}{2\eta^k} + \eta^k L_f^2 + \tfrac{M_g}{\eta^k} + 4L_gD_xD_{\lambda}M\big)E[\norm{x_t^k - x_{t+1}^k}^2]. \label{eq:b35}
    \end{align}
\end{lemma}
\begin{proof}
    For $E[\inprod{\delta_t^{x^k}}{x - x_{t+1}^k}]$
    \begin{align*}
        E[\inprod{\delta_t^{x^k}}{x - x_{t+1}^k}] = E[\inprod{\delta_t^{x^k}}{x_t^k - x_{t+1}^k}],
    \end{align*}
    because $E[\inprod{\delta_t^{x^k}}{x - x_t^k}] = 0$ for any fixed $x$. Recall that 
    \begin{align*}
        & \delta_t^{x^k} = \grad f_{i_t^k}(x_t^k) - \grad f_{i_t^k}(\wt{x}^{k-1}) - [\grad f(x_t^k) - \grad f(\wt{x}^{k-1})]  + \lambda_{t+1}^k\grad g_{p_t^k}(x_t^k) - \wt{\lambda}^{k-1}\grad g_{p_t^k}(\wt{x}^{k-1}) \nonumber\\
        & - \lambda_{t+1}^k\grad g(x_t^k) + \wt{\lambda}^{k-1}\grad g(\wt{x}^{k-1}).
    \end{align*}
    For $E[\inprod{\grad f_{i_t^k}(x_t^k) - \grad f_{i_t^k}(\wt{x}^{k-1}) - [\grad f(x_t^k) - \grad f(\wt{x}^{k-1})]}{x_t^k - x_{t+1}^k}]$:
    \begin{align*}
        & E[\inprod{\grad f_{i_t^k}(x_t^k) - \grad f_{i_t^k}(\wt{x}^{k-1}) - [\grad f(x_t^k) - \grad f(\wt{x}^{k-1})]}{x_t^k - x_{t+1}^k}] \nonumber\\
        & \leq \tfrac{\eta^k}{2}E[\norm{\grad f_{i_t^k}(x_t^k) - \grad f_{i_t^k}(\wt{x}^{k-1}) - [\grad f(x_t^k) - \grad f(\wt{x}^{k-1})]}^2] + \tfrac{1}{2\eta^k}E[\norm{x_t^k - x_{t+1}^k}^2] \nonumber\\
        & \leq \eta^k L_f^2E[D_x^k(x_{t+1}^k)] + (\tfrac{1}{2\eta^k} + \eta^k L_f^2)E[\norm{x_t^k - x_{t+1}^k}^2], 
    \end{align*}
    where the first inequality comes from \eqref{eq:p9} and  \eqref{eq:p10}. The last inequality comes from $Var(X) \leq E[X^2]$, triangle inequality and \eqref{eq:p7}. Also
    \begin{align}
        & E[\inprod{\lambda_{t+1}^k\grad g_{p_t^k}(x_t^k) - \wt{\lambda}^{k-1}\grad g_{p_t^k}(\wt{x}^{k-1})}{x_t^k - x_{t+1}^k}] \leq E[\norm{\lambda_{t+1}^k\grad g_{p_t^k}(x_t^k) - \wt{\lambda}^{k-1} \grad g_{p_t^k} (x_t^k)}\norm{x_t^k - x_{t+1}^k}] \nonumber\\
        & + E[\norm{\wt{\lambda}^{k-1} \grad g_{p_t^k} (x_t^k) - \wt{\lambda}^{k-1}\grad g_{p_t^k}(\wt{x}^{k-1})}\norm{x_t^k - x_{t+1}^k}] \nonumber\\
        & \leq M_g E[\norm{\lambda_{t+1}^k - \wt{\lambda}^{k-1}}\norm{x_t^k - x_{t+1}^k}] + E[\norm{\wt{\lambda}^{k-1}} \norm{\grad g_{p_t^k} (x_t^k) - \grad g_{p_t^k}(\wt{x}^{k-1})}\norm{x_t^k - x_{t+1}^k}] \nonumber\\
        & \leq \tfrac{M_g \eta^k}{2} E[\norm{\lambda_{t+1}^k - \wt{\lambda}^{k-1}}^2] + \tfrac{M_g}{2\eta^k}E[\norm{x_t^k - x_{t+1}^k}^2] + L_gE[\norm{\wt{\lambda}^{k-1}} \norm{x_t^k - \wt{x}^{k-1}}\norm{x_t^k - x_{t+1}^k}] \nonumber\\
        & \leq \tfrac{M_g \eta^k}{2} E[\norm{\lambda_{t+1}^k - \wt{\lambda}^{k-1}}^2] + \tfrac{M_g}{2\eta^k}E[\norm{x_t^k - x_{t+1}^k}^2] + L_gD_x(E[\norm{\wt{\lambda}^{k-1}}]E[\norm{x_t^k - x_{t+1}^k}] \nonumber\\
        & + \sqrt{E[\norm{\wt{\lambda}^{k-1}}^2]E[\norm{x_t^k - x_{t+1}^k}^2]}) \leq \tfrac{M_g \eta^k}{2} E[D_{\lambda}^k(\lambda_{t+1}^k)] + \big(\tfrac{M_g}{2\eta^k} + 2L_gD_xD_{\lambda} M\big)E[\norm{x_t^k - x_{t+1}^k}^2], \label{eq:b33}
    \end{align}
    where the first inequality in \eqref{eq:b33} comes from \eqref{eq:p9} and \eqref{eq:p3}, the second inequality comes from \eqref{eq:p1}, \eqref{eq:p2}, and \eqref{eq:p10}, the third inequality comes from upper bounding $\norm{x_t^k - \wt{x}^{k-1}}$ with $D_x$ and also $E[XY] \leq E[X]E[Y] + (E[Y^2]E[X^2])^{1/2}$, and the last inequality comes from upper bounding $E[\norm{\wt{\lambda}^{k-1}}]$ and $\sqrt{E[\norm{\wt{\lambda}^{k-1}}^2]}$ with $D_{\lambda}$, \ref{eq:RMS}, and \eqref{eq:p7}. $E[\inprod{\lambda_{t+1}^k\grad g_{p_t^k}(x_t^k) - \wt{\lambda}^{k-1}\grad g_{p_t^k}(\wt{x}^{k-1})}{x_t^k - x_{t+1}^k}]$ and $E[\inprod{\lambda_{t+1}^k\grad g(x_t^k) - \wt{\lambda}^{k-1}\grad g(\wt{x}^{k-1})}{x_t^k - x_{t+1}^k}]$
    have the same upper bound so 
    \begin{align*}
        & E[\inprod{\delta_t^{x^k}}{x - x_{t+1}^k}] \leq M_g \eta^k E[D_{\lambda}^k(\lambda_{t+1}^k)] + \eta^k L_f^2E[D_x^k(x_{t+1}^k)] \nonumber\\
        & + \big(\tfrac{1}{2\eta^k} + \eta^k L_f^2 + \tfrac{M_g}{\eta^k} + 4L_gD_xD_{\lambda}M\big)E[\norm{x_t^k - x_{t+1}^k}^2].
    \end{align*}
\end{proof}
We use the upper bounds in Lemma \ref{eq:a38} and \eqref{eq:b35} for \eqref{eq:a26} and rearrange the terms
\begin{align}
    & E[Q(\wt{Z}^K,z^*)] \leq \tfrac{1}{S^K}\bigg(\tfrac{1}{\eta^1}E[\norm{\lambda^* - v_0^1}^2] + \tsum_{k=1}^{K-1}(\tfrac{1}{\eta^{k+1}} - \tfrac{1}{\eta^k})E[\norm{\lambda^* - v_{T^k}^k}^2] - \tfrac{1}{\eta^K}E[\norm{\lambda^* - v_{T^K}^K}^2]\nonumber \\
    & + \tfrac{1}{\eta^1}E[\norm{\lambda^* - w_0^1}^2] + \tsum_{k=1}^{K-1}(\tfrac{1}{\eta^{k+1}} - \tfrac{1}{\eta^k})E[\norm{\lambda^* - w_{T^k}^k}^2]  - \tfrac{1}{\eta^K}E[\norm{\lambda^* - w_{T^K}^K}^2] \nonumber\\
    & +  \tfrac{1 - \gamma_\lambda^1}{2\sigma^1}E[\norm{\lambda^* - \lambda_0^1}^2] + \tsum_{k=1}^{K-1}(\tfrac{1 - \gamma_\lambda^{k+1}}{2\sigma^{k+1}} - \tfrac{1 - \gamma_\lambda^k}{2\sigma^k})E[\norm{\lambda^* - \lambda_{T^k}^k}^2] - \tfrac{1 - \gamma_\lambda^K}{4\sigma^K}E[\norm{\lambda^* - \lambda_{T^K}^K}^2] 
     \nonumber\\
    & + \tfrac{\gamma_\lambda^1T^1}{2\sigma^1}E[D_\lambda^1(\lambda^*)] + \tsum_{k=1}^{K-1}(\tfrac{\gamma_\lambda^{k+1}T^{k+1}}{2\sigma^{k+1}} - \tfrac{\gamma_\lambda^kT^k}{2\sigma^k})E[D_\lambda^{k+1}(\lambda^*)] - \tfrac{\gamma_\lambda^K T^K}{2\sigma^K}E[D_\lambda^{K+1}(\lambda^*)]\nonumber\\
    & + \tsum_{k=1}^{K}\tsum_{t=0}^{T^k-1} (\tfrac{1}{\eta^k} + \tfrac{M_g}{2} -  \tfrac{1 - \gamma_\lambda^k}{2\sigma^k})E[\norm{\lambda_t^k - \lambda_{t+1}^k}^2] + \tsum_{k=1}^{K}\tsum_{t=0}^{T^k - 1}(\eta^kM_g - \tfrac{\gamma_\lambda^k}{2\sigma^k})E[D_\lambda^k(\lambda_{t+1}^k)] \nonumber\\
    & + \tsum_{k=1}^{K}\tsum_{t=0}^{T^k - 1}(\tfrac{\norm{\lambda^*}L_g + M_g + L_f + 8L_gD_xD_{\lambda}M + 2\eta^k (L_f^2  + 2M_g^2)}{2} + \tfrac{1 + 2M_g}{2\eta^k}  - \tfrac{1 - \gamma_x^k}{2\tau^k}) E[\norm{x_t^k - x_{t+1}^k}^2]  \nonumber\\
    & + \tfrac{1 - \gamma_x^1}{2\tau^1}E[\norm{x^* - x_0^1}^2] + \tsum_{k=1}^{K-1}(\tfrac{1 - \gamma_x^{k+1}}{2\tau^{k+1}} - \tfrac{1 - \gamma_x^k}{2\tau^k})E[\norm{x^* - x_{T^k}^k}^2] - \tfrac{1 - \gamma_x^K}{2\tau^K}E[\norm{x^* - x_{T^K}^K}^2] \nonumber\\
    & + \tfrac{\gamma_x^1T^1}{2\tau^1}E[D_x^1(x^*)] + \tsum_{k=1}^{K-1}(\tfrac{\gamma_x^{k+1}T^{k+1}}{2\tau^{k+1}} - \tfrac{\gamma_x^kT^k}{2\tau^k})E[D_x^{k+1}(x^*)] - \tfrac{\gamma_x^K T^K}{2\tau^K}E[D_x^{K+1}(x^*)] \nonumber\\
    & + \tsum_{k=1}^{K}\tsum_{t=0}^{T^k - 1}(\eta^kL_f^2- \tfrac{\gamma_x^k}{2\tau^k})E[D_x^k(x_{t+1}^k)] + \tfrac{\sigma^KM_g^2}{1-\gamma_{\lambda}^K}E[\norm{x_{T^K-1}^K - x_{T^K}^K}^2]\bigg). \label{eq:a42}
\end{align}
We will see again $ \tfrac{\sigma^KM_g^2}{1-\gamma_{\lambda}^K}$ could be eliminated by the coefficient of $E[\norm{x_t^k - x_{t+1}^k}^2]$. 
\begin{remark} \label{eq:th4}
    The following step size conditions are necessary for the convergence of Algorithm \ref{alg:alg2} via the gap inequality \eqref{eq:a42}:
    \begin{align}
        \tfrac{1}{\eta^k} + \tfrac{M_g}{2} -  \tfrac{1 - \gamma_\lambda^k}{2\sigma^k} &\leq 0 & \tfrac{M_g}{2} - \tfrac{1 - \gamma_\lambda^K}{2\sigma^K} &\leq 0 & \eta^kM_g- \tfrac{\gamma_\lambda^k}{2\sigma^k} &\leq 0, \label{eq:a47}
    \end{align}
    \begin{align}
        \tfrac{\norm{\lambda^*}L_g + M_g + L_f + 8L_gD_xD_{\lambda}M + 2\eta^k (L_f^2  + 2M_g^2)}{2} + \tfrac{1 + 2M_g}{2\eta^k}  - \tfrac{1 - \gamma_x^k}{2\tau^k} &\leq 0 
        &\eta^kL_f^2- \tfrac{\gamma_x^k}{2\tau^k} \leq 0, \label{eq:a51}
    \end{align}
    In particular, using conditions \eqref{eq:a47} and \eqref{eq:a51}, we simplify the upper bound on the gap in \eqref{eq:a42}.
\end{remark} 

\subsubsection*{Constant step size policy.}
\begin{theorem}
    Assume that the parameters are chosen as follows:
    \begin{align*}
        \eta &= \tfrac{1}{q^{1/2}} & \sigma & = \tfrac{1}{2q^{1/2}(2 + M_g)} &\tau &= \tfrac{1}{2\big(2c^{\prime} + q^{1/2}c\big)} & \gamma_{\lambda} &= \tfrac{1}{q} &\gamma_x &= \tfrac{1}{q} & T &= q, 
    \end{align*}
    where $c =\tfrac{2L_f^2 + 1+2M_g + 4M_g^2}{2}$, and $c' = \tfrac{\norm{\lambda^*}L_g + M_g + L_f + 8L_gD_xD_{\lambda}M + 2 (L_f^2  + 2M_g^2)}{2}$. The expected gap inequality will be:
    \begin{align}
        & E[Q(\wt{Z}^K,z^*)] \leq \tfrac{1}{Kq}\big(q^{1/2}(E[\norm{\lambda^* - v_0^1}^2] + E[\norm{\lambda^* - w_0^1}^2]) + ((\tfrac{2 + M_g)(q-1)}{q^{1/2}})E[\norm{\lambda^* - \lambda_0^1}^2]  \nonumber\\
        & + q^{1/2}(2+M_g)E[D_\lambda^1(\lambda^*)] + \tfrac{(2c^{\prime} + q^{1/2}c)(q - 1)}{q}E[\norm{x^* - x_0^1}^2]  + (2c^{\prime} + q^{1/2}c)E[D_x^1(x^*)] \big) = \epsilon . \label{eq:a55}
    \end{align}
    The gradient computations to reach $\epsilon$-gap for the objective function and the function constraint will be $O(\tfrac{m}{\epsilon q^{1/2}})$ and $O(\tfrac{n}{\epsilon q^{1/2}})$ respectively using the gap inequality \eqref{eq:a42}.
\end{theorem}
\begin{proof}
    Note that step-size parameters above satisfy \eqref{eq:a47} and \eqref{eq:a51}. Using these relations inside \eqref{eq:a42}, then eliminating the negative terms on the RHS and applying the parameter values to remaining terms of the RHS, we obtain \eqref{eq:a55}. The number of gradient computation for the objective and the constraint are $mK+\tsum_{k=1}^{K}T^k=mK+TK = O(\tfrac{m}{\epsilon q^{1/2}})$, and $nK+\tsum_{k=1}^{K}T^k=nK+TK = O(\tfrac{n}{\epsilon q^{1/2}})$ respectively.  
\end{proof} 
\subsubsection*{Non-constant step size policy.}
\begin{theorem}
    Assume we have the following parameters setup:
    \begin{align*}
        \eta^k &= \tfrac{1}{k}  &\gamma_x^k &= \tfrac{1}{(k+1)^2}  &\gamma_\lambda^k &= \tfrac{1}{(k+1)^2} &\sigma^k &= \tfrac{k}{2(1+M_g)(k+1)^2} &\tau^k &= \tfrac{1}{4(c' + Ak)} & T^k &=(k+1)^2,
    \end{align*} 
    where $A = \tfrac{1 + 2M_g + 4M_g^2 + 6L_f^2}{2}$. Then the expected gap will be: 
    \begin{align}
        &E[Q(\wt{Z}^K,z^*)] \leq \tfrac{6}{K(2K^2 + 9K + 13)}\big(E[\norm{\lambda^* - v_0^1}^2] + E[\norm{\lambda^* - w_0^1}^2]  + \big((4 + 2M_g)K - 2 + \tfrac{1+M_g}{K}\big)D_{\lambda}^2   \nonumber\\
        & + 4(1+M_g)E[D_{\lambda}^1(\lambda^*)] + (c'+A)(\tfrac{3}{2}E[\norm{x^* - x_0^1}^2] + 2E[D_x^1(x^*)]) + 2(c' + AK)(2 - \tfrac{1}{(K+1)^2})D_x^2\big) = \epsilon. \label{eq:a307} 
    \end{align}
    The gradient computations for the objective and the constraint will be $O(\tfrac{m}{\epsilon^{1/2}} + \tfrac{1}{\epsilon^{1.5}})$, and $O(\tfrac{n}{\epsilon^{1/2}} + \tfrac{1}{\epsilon^{1.5}})$.
\end{theorem}. 
\begin{proof}
    With the setup in this theorem applied to the inequality \eqref{eq:a42}, the expected gap will be \eqref{eq:a307}, and it satisfies the conditions in \ref{eq:th4}. From this inequality, $K = O(\tfrac{1}{\epsilon^{1/2}})$. The number of gradient computations for the objective and the constraints will be $Km + \tsum_{k=1}^K T^k = O(\tfrac{m}{\epsilon^{1/2}} + \tfrac{1}{\epsilon^{1.5}})$ and $Kn + \tsum_{k=1}^K T^k = O(\tfrac{n}{\epsilon^{1/2}} + \tfrac{1}{\epsilon^{1.5}})$
\end{proof}  
We conclude this section by recalling Remark \ref{rem:manage_lambda_in_step_size} which allows us to ignore $\lambda^*$ related terms in stepsize policies for both constant and nonconstant cases.
\section{Strongly convex objective function.} \label{sec-strong}
In this section, we assume that the objective function is strongly convex and use this assumption to accelerate the algorithm and reduce the gradient computation complexity simultaneously. Here is the algorithm for the general case of a strongly convex objective function. Similar to the convex objective function, we have semi-finite-sum and fully-finite-sum cases. For each case, we have the convergence analysis both for non-constant step size policy. For strong convexity, we give $\alpha_t^k$ weight for each iteration to improve the rate compared with the merely convex case. As it has been shown before, Strongly convex objective function could enable us to improve the rate of the convergence for the algorithms we want to design. For the strongly convex case, we only add the extra assumption of strong convexity of each component of the objective function.
\begin{align*}
    f_i(y) \geq f_i(x) + \inprod{\grad f_i(x)}{y - x} + \tfrac{\mu_f}{2}\norm{x - y}^2, i \in [m], \hspace{0.5em} x,y \in X.
\end{align*}
\begin{algorithm}[!t]
	\begin{algorithmic}[1]
		\State {\bf Input: }$\wt{x}^0 \in X, \wt{\lambda}^0 = 0$. 
		\State Set $(x_0^1, \lambda_0^1) =(\wt{x}^0,\wt{\lambda}^0)$.% and $g(x^{-1}) = g(x^0)$.
		\For{$k=1,2,\dots,K$}
        \State $G_{\lambda} \gets g(\tilde{x}^{k-1})$
		\State $G_x \gets \grad f(\tilde{x}^{k-1}) + \grad g(\tilde{x}^{k-1}) \wt{\lambda}^{k-1}$
        \State In epoch $k$
        \For{ $t = 0,1, \dots, T^k-1$}
        \State {Pick a random subset $j_t^k \in [n]$}
        \State $\wh{\lambda}^k_t \gets (1-\gamma^k_\lambda)\lambda^k_t + \gamma^k_\lambda \wt{\lambda}^{k-1}$
		\State $\xi_t^k \gets g_{j_t^k}(x^k_t) - g_{j_t^k}(\wt{x}^{k-1}) + G_\lambda$
		\State $q^k_t \gets \ell{g_{j_t^k}}(x^k_t) - g_{j_t^k}(x^k_{t-1})$
        \State $\lambda_{t+1}^k \gets \argmin_{\lambda \ge 0} -\inprod{\xi_t^k + \theta_{t+1}^k q^k_t}{\lambda} + \frac{1}{2\sigma_{t+1}^k}\norm{\lambda - \wh{\lambda}^k_t}^2$
        \State {Pick random subsets $i_t^k \in [m]$ and $p_t^k \in [n]$}
        \State $\wh{x}_t^k \gets (1-\gamma_x^k)x^k_t + \gamma_x^k \wt{x}^{k-1}$
        \State $\beta_t^k \gets \grad f_{i_t^k}(x^k_t)+ \lambda_{t+1}^k \grad g_{p_t^k}(x^k_t)  - [\grad f_{i_t^k}(\wt{x}^{k-1}) + \grad g_{p_t^k}(\wt{x}^{k-1})\wt{\lambda}^{k-1}] + G_x$
        \State $x_{t+1}^k \gets \argmin_{x \in X}  \inprod{\beta_t^k}{x} + \frac{1}{2\tau_{t+1}^k} \norm{x-\wh{x}^k_t}^2$
        \EndFor
        \State {\bf Output: }$\wt{x}^k = \tfrac{\tsum_{t=1}^{T^k}\frac{\gamma_x^k\alpha_t^k}{2\tau_t^k}x_t^k}{\tsum_{t=1}^{T^k}\frac{\gamma_x^k\alpha_t^k}{2\tau_t^k}}$ {and} $\wt{\lambda}^k = \tfrac{\tsum_{t=1}^{T^k}\frac{\gamma_{\lambda}^k\alpha_t^k}{2\sigma_t^k}\lambda_t^k}{\tsum_{t=1}^{T^k}\frac{\gamma_{\lambda}^k\alpha_t^k}{2\sigma_t^k}}$
        \State $(x^{k+1}_0, \lambda^{k+1}_0, x^{k+1}_{-1}, \lambda^{k+1}_{-1}) = (x^k_{T^k},\lambda^k_{T^k}, x^k_{T^k-1}, \lambda^k_{T^k-1})$
		\EndFor
	\end{algorithmic}
	\caption{Variance-reduced anchored constraint extrapolation for strongly convex objective} \label{alg:alg3}
\end{algorithm}
The structure of Algorithm \ref{alg:alg3} is similar to Algorithm \ref{alg:alg2}. The only differences are in extrapolation, prox steps, and averaging. We have $\sigma_t^k$, and $\tau_t^k$ that change with every iteration, while in Algorithm \ref{alg:alg2} these change with epochs. Also in extrapolation step of Algorithm \ref{alg:alg3}, $\theta_t^k$ plays a key role, which helps us to eliminate inner products in the analysis, and with this, we improve the convergence rate of the algorithm. We define new averaging for the norms 
\begin{align}
    D_{\lambda}^k(\lambda) & := \tfrac{1}{s^{k-1}}\tsum_{l=1}^{T^{k-1}}\tfrac{\gamma_{\lambda}^{k-1}\alpha_l^{k-1}}{2\sigma_l^{k-1}}\norm{\lambda - \lambda_l^{k-1}}^2 & D_x^k(x) & := \tfrac{1}{r^{k-1}}\tsum_{l=1}^{T^{k-1}}\tfrac{\gamma_x^{k-1}\alpha_l^{k-1}}{2\tau_l^{k-1}}\norm{x - x_l^{k-1}}^2, \label{eq:p12}
\end{align}
where, $s^{k-1} = \tsum_{l=1}^{T^{k-1}}\tfrac{\gamma_{\lambda}^{k-1}\alpha_l^{k-1}}{2\sigma_l^{k-1}}$, and $r^{k-1} = \tsum_{l=1}^{T^{k-1}}\tfrac{\gamma_x^{k-1}\alpha_l^{k-1}}{2\tau_l^{k-1}}$. Also $u^k = \tsum_{t=1}^{T^k}\alpha_t^k$. Using Jensen's inequality for the gap will be
\begin{align}
    E[Q(\wt{z}^k,z^*)] &\leq \tfrac{1}{u^k}\tsum_{t=1}^{T^k} \alpha_t^kE[Q(z_t^k,z^*)] & E[Q(\wt{Z}^K,z^*)] &\leq \tfrac{1}{U^K}\tsum_{k=1}^{K} u^kE[Q(\wt{z}^k,z^*)], \label{eq:p13}
\end{align}
where $\wt{z}^k = \tfrac{1}{u^k}\tsum_{t=1}^{T^k}\alpha_t^k z_t^k$, $U^K = \tsum_{k=1}^{K}u^k$, and $\wt{Z}^K = \tfrac{1}{U^K}\tsum_{k=1}^{K}u^k\wt{z}^k$. We start with the semi-finite case. 
\subsection{Semi-finite-sum case.} \label{sec-semi}
\begin{lemma}
    The gap originating from Algorithm \ref{alg:alg3} is as follows:
    \begin{align}
        & Q(z_{t+1}^k, z) \leq \inprod{\delta_t^{\lambda^k}}{\lambda_{t+1}^k - \lambda} 
        - \inprod{\bar{q}_{t+1}^k}{\lambda_{t+1}^k - \lambda} 
        + \inprod{\theta_{t+1}^k \bar{q}_t^k}{\lambda_t^k - \lambda} + \inprod{\theta_{t+1}^k(q_t^k - \bar{q}_t^k)}{\lambda_t^k - \lambda} \nonumber \\
        & + \tfrac{\theta_{t+1}^k M_g}{2\rho_{t+1}^k}\norm{x_t^k - x_{t-1}^k}^2 + \tfrac{\theta_{t+1}^k\rho_{t+1}^k M_g}{2}\norm{\lambda_{t+1}^k - \lambda_t^k}^2 + \tfrac{1 - \gamma_\lambda^k}{2\sigma_{t+1}^k} \big( \norm{\lambda - \lambda_t^k}^2 - \norm{\lambda - \lambda_{t+1}^k}^2 - \norm{\lambda_t^k - \lambda_{t+1}^k}^2 \big) \nonumber \\
        & + \tfrac{\gamma_\lambda^k}{2\sigma_{t+1}^k} \big( D_\lambda^k(\lambda) 
        - \norm{\lambda - \lambda_{t+1}^k}^2 
        - D_\lambda^k(\lambda_{t+1}^k) \big) + (\tfrac{\norm{\lambda}L_g}{2} + \tfrac{L_f}{2})\norm{x_t^k - x_{t+1}^k}^2 + \inprod{\delta_t^{x^k}}{x - x_{t+1}^k} \nonumber \\
        & + (\tfrac{1-\gamma_x^k}{2\tau_{t+1}^k} - \tfrac{\mu_f}{2})\norm{x - x_t^k}^2 - \tfrac{1-\gamma_x^k}{2\tau_{t+1}^k}\norm{x - x_{t+1}^k}^2 - \tfrac{1-\gamma_x^k}{2\tau_{t+1}^k}\norm{x_{t+1}^k - x_t^k}^2 \nonumber \\
        & + \tfrac{\gamma_x^k}{2\tau_{t+1}^k} \big( D_x^k(x) 
        - \norm{x - x_{t+1}^k}^2 
        - D_x^k(x_{t+1}^k) \big). \label{eq:a58}
    \end{align}
\end{lemma}
\begin{proof}
    According to the update rule of $\lambda_{t+1}^k$ we have:
    \begin{align}
	& \inprod{\xi_t^k + \theta_{t+1}^k q_t^k}{\lambda - \lambda_{t+1}^k} \leq \tfrac{1}{\sigma_{t+1}^k} \inprod{\lambda_{t+1}^k - \wh{\lambda}_t^k}{\lambda - \lambda^k_{t+1}} = \tfrac{1}{\sigma_{t+1}^k} \inprod{\lambda_{t+1}^k - (1-\gamma_\lambda^k) \lambda_t^k - \gamma_\lambda^k \tilde{\lambda}^{k-1}}{\lambda - \lambda_{t+1}^k} \nonumber\\
	& = \tfrac{1 - \gamma^k_\lambda}{\sigma_{t+1}^k} \inprod{\lambda_{t+1}^k - \lambda_t^k}{\lambda - \lambda_{t+1}^k} + \tfrac{\gamma^k_\lambda}{\sigma_{t+1}^ks^{k-1}} \tsum_{l=1}^{T^{k-1}} \tfrac{\gamma_{\lambda}^{k-1}\alpha_l^{k-1}}{2\sigma_l^{k-1}}\inprod{\lambda_{t+1}^k - \lambda_l^{k-1}}{\lambda - \lambda_{t+1}^k} \nonumber\\
	& \leq \tfrac{1 - \gamma^k_\lambda}{2\sigma_{t+1}^k}\big(\norm{\lambda - \lambda_t^k}^2 - \norm{\lambda - \lambda_{t+1}^k}^2 - \norm{\lambda_t^k - \lambda_{t+1}^k}^2 \big) \nonumber\\
	& + \tfrac{\gamma^k_\lambda}{2\sigma_{t+1}^k s^{k-1}} \tsum_{l=1}^{T^{k-1}} \tfrac{\gamma_{\lambda}^{k-1}\alpha_l^{k-1}}{2\sigma_l^{k-1}}\big(\norm{\lambda - \lambda_l^{k-1}}^2 - \norm{\lambda - \lambda_{t+1}^k}^2 - \norm{\lambda_l^{k-1} - \lambda_{t+1}^k}^2 \big). \label{eq:a59}
\end{align}	
The first inequality in \eqref{eq:a59} comes from optimality of $\lambda_{t+1}^k$ in line 12 of Algorithm \ref{alg:alg3}, second equality comes from the definition of $\wh{\lambda}_t^k$ and the last inequality comes from Pythagorean inequality for inner products. Defining $\delta_t^{\lambda^k} := \xi_t^k - g(x^k_t)$ and $\bar{q}^k_t := \ell_g(x^k_t) - g(x^k_{t-1})$ and
\begin{align}
	&\inprod{\xi_t^k + \theta_{t+1}^k q_t^k}{\lambda_{t+1}^k - \lambda} = \inprod{\delta_t^{\lambda^k}}{\lambda_{t+1}^k - \lambda} + \inprod{g(x_t^k)}{\lambda_{t+1}^k - \lambda}  + \inprod{\theta_{t+1}^k q_t^k}{\lambda_t^k - \lambda} + \inprod{\theta_{t+1}^k q_t^k}{\lambda_{t+1}^k - \lambda^k_t} \nonumber \\
	&\Rightarrow \inprod{\ell_g(x_{t+1}^k)}{\lambda - \lambda_{t+1}^k} \leq \inprod{\delta_t^{\lambda^k}}{\lambda_{t+1}^k - \lambda} - \inprod{\bar{q}_{t+1}^k}{\lambda_{t+1}^k - \lambda} + \inprod{\theta_{t+1}^k (q_t^k - \bar{q}_{t}^k)}{\lambda_t^k - \lambda} \nonumber\\
    & + \inprod{\theta_{t+1}^k \bar{q}_t^k}{\lambda_t^k - \lambda} + \inprod{\theta_{t+1}^k q_t^k}{\lambda_{t+1}^k - \lambda_t^k} + \tfrac{1 - \gamma^k_\lambda}{2\sigma_{t+1}^k} \big( \norm{\lambda - \lambda_t^k}^2 - \norm{\lambda - \lambda_{t+1}^k}^2 - \norm{\lambda_t^k - \lambda_{t+1}^k}^2 \big) \nonumber \\
	& + \tfrac{\gamma^k_\lambda}{2\sigma_{t+1}^k} \big( D_\lambda^k(\lambda) - \norm{\lambda - \lambda_{t+1}^k}^2 - D_\lambda^k(\lambda_{t+1}^k) \big). \label{eq:a60}
\end{align}
Also
\begin{align}
	& \norm{q_t^k} = \norm{\ell_{g_{j_t^k}}(x_t^k) - g_{j_t^k}(x_{t-1}^k)} = \norm{g_{j_t^k}(x_{t-1}^k) + \inprod{\grad g_{j_t^k}(x_{t-1}^k)}{x_t^k - x_{t-1}^k} - g_{j_t^k}(x_{t-1}^k)} \nonumber\\
    & \leq M_g\norm{x_t^k - x_{t-1}^k}. \label{eq:a61} 
\end{align}
Changing bound inside \eqref{eq:a60} and using \eqref{eq:p9}, and \eqref{eq:p10} for $\inprod{\theta_{t+1}^k q_t^k}{\lambda_{t+1}^k - \lambda_t^k}$
\begin{align}
	&\inprod{\ell_g(x_{t+1}^k)}{\lambda - \lambda_{t+1}^k} \leq \inprod{\delta_t^{\lambda^k}}{\lambda_{t+1}^k - \lambda} - \inprod{\bar{q}_{t+1}^k}{\lambda_{t+1}^k - \lambda} + \inprod{\theta_{t+1}^k (q_t^k - \bar{q}_{t}^k)}{\lambda_t^k - \lambda} \nonumber \\
	& + \inprod{\theta_{t+1}^k \bar{q}_t^k}{\lambda_t^k - \lambda} + \tfrac{\theta_{t+1}^k M_g}{2\rho_{t+1}^k}\norm{x_t^k - x_{t-1}^k}^2 + \tfrac{\theta_{t+1}^k \rho_{t+1}^k M_g}{2}\norm{\lambda_{t+1}^k - \lambda_t^k}^2 \nonumber \\
	& + \tfrac{1 - \gamma^k_\lambda}{2\sigma_{t+1}^k} \big( \norm{\lambda - \lambda_t^k}^2 - \norm{\lambda - \lambda_{t+1}^k}^2 - \norm{\lambda_t^k - \lambda_{t+1}^k}^2 \big) + \tfrac{\gamma^k_\lambda}{2\sigma_{t+1}^k} \big( D_\lambda^k(\lambda) - \norm{\lambda - \lambda_{t+1}^k}^2 - D_\lambda^k(\lambda_{t+1}^k) \big). \label{eq:a62}
\end{align}
Modification on $\inprod{\ell_g(x_{t+1}^k)}{\lambda - \lambda_{t+1}^k}$
\begin{align*}
\inprod{\ell_g(x_{t+1}^k)}{\lambda - \lambda_{t+1}^k} = \inprod{\ell_g(x_{t+1}^k)}{\lambda} + \inprod{\ell_g(x_{t+1}^k)}{- \lambda_{t+1}^k}. 
\end{align*}
Also: 
\begin{align*}
\inprod{\ell_g(x_{t+1}^k)}{\lambda} = \inprod{\ell_g(x_{t+1}^k) - g(x_{t+1}^k) + g(x_{t+1}^k)}{\lambda}, 
\end{align*}
which by using \eqref{eq:p9} and \eqref{eq:p5}
\begin{align}
\inprod{g(x_{t+1}^k) - \ell_g(x_{t+1}^k) }{\lambda} &\leq \tfrac{\norm{\lambda}L_g}{2}\norm{x_t^k - x_{t+1}^k}^2.  \label{eq:a63}
\end{align}
Taking $\inprod{\ell_g(x_{t+1}^k)}{- \lambda_{t+1}^k}$ to the other side of \eqref{eq:a62} and using \eqref{eq:a63}
\begin{align}
    & \inprod{g(x_{t+1}^k)}{\lambda} \leq \inprod{\delta_t^{\lambda^k}}{\lambda_{t+1}^k - \lambda} 
    - \inprod{\bar{q}_{t+1}^k}{\lambda_{t+1}^k - \lambda} 
    + \inprod{\theta_{t+1}^k\bar{q}_t^k}{\lambda_t^k - \lambda} + \inprod{\theta_{t+1}^k(q_t^k - \bar{q}_t^k)}{\lambda_t^k - \lambda} \nonumber \\
    & + \tfrac{\theta_{t+1}^k M_g}{2\rho_{t+1}^k}\norm{x_t^k - x_{t-1}^k}^2 + \tfrac{\theta_{t+1}^k \rho_{t+1}^k M_g}{2}\norm{\lambda_{t+1}^k - \lambda_t^k}^2 + \tfrac{1 - \gamma_\lambda^k}{2\sigma_{t+1}^k} ( \norm{\lambda - \lambda_t^k}^2 
    - \norm{\lambda - \lambda_{t+1}^k}^2 
    - \norm{\lambda_t^k - \lambda_{t+1}^k}^2 ) \nonumber \\
    & + \tfrac{\gamma_\lambda^k}{2\sigma_{t+1}^k} ( D_\lambda^k(\lambda) 
    - \norm{\lambda - \lambda_{t+1}^k}^2 
    - D_\lambda^k(\lambda_{t+1}^k) ) + \inprod{\ell_g(x_{t+1}^k)}{\lambda_{t+1}^k} 
    + \tfrac{\norm{\lambda}L_g}{2} \norm{x_t^k - x_{t+1}^k}^2.
      \label{eq:a64}
\end{align}
By optimality of $x^k_{t+1}$ from line 16 of Algorithm \ref{alg:alg3}:
\begin{align}
	& \inprod{\beta_t^k}{x_{t+1}^k - x} \leq \frac{1}{\tau_{t+1}^k} \inprod{x_{t+1}^k - \wh{x}_t^k}{x - x_{t+1}^k}  = \frac{1 - \gamma_x^k}{\tau_{t+1}^k} \inprod{x_{t+1}^k - x_t^k}{x - x_{t+1}^k} 
	+ \tfrac{\gamma_x^k}{\tau_{t+1}^k} \inprod{x_{t+1}^k - \tilde{x}^{k-1}}{x - x_{t+1}^k} \nonumber \\
	&\leq \tfrac{1 - \gamma_x^k}{2\tau_{t+1}^k} \big( \norm{x - x_t^k}^2 - \norm{x - x_{t+1}^k}^2 - \norm{x_t^k - x_{t+1}^k}^2 \big) \nonumber \\
	& + \tfrac{\gamma_x^k}{2\tau_{t+1}^k r^{k-1}} \tsum_{l=1}^{T^{k-1}} \tfrac{\gamma_x^{k-1}\alpha_l^{k-1}}{2\tau_l^{k-1}}\big( \norm{x - x_l^{k-1}}^2 - \norm{x - x_{t+1}^k}^2 - \norm{x_l^{k-1} - x_{t+1}^k}^2 \big). \label{eq:a65}
\end{align}
The first inequality comes from optimality of $x_{t+1}^k$, equality comes from definition of $\wh{x}_t^k$, and last inequality comes from Pythagorean inequality for inner products. As $\delta^{x^k}_t := \beta_t^k -  \grad \varphi(x^k_t,  \lambda^k_{t+1})$
\begin{align}
	&\inprod{\delta^{x^k}_t}{x_{t+1}^k - x} \leq \inprod{\nabla_x \varphi(x_t^k,\lambda_{t+1}^k)}{x - x_{t+1}^k} + \tfrac{1-\gamma_x^k}{2\tau_{t+1}^k} \big(\norm{x - x_t^k}^2 - \norm{x - x_{t+1}^k}^2 - \norm{x_{t+1}^k - x_t^k}^2\big) \nonumber\\
    & + \tfrac{\gamma_x^k}{2\tau_{t+1}^k} \big(D_x^k(x) - \norm{x - x_{t+1}^k}^2 - D_x^k(x_{t+1}^k)\big). \label{eq:a66}
\end{align}
For $\inprod{\grad_x \varphi(x_t^k,\lambda_{t+1}^k)}{x - x_{t+1}^k}$, using strong convexity of $\varphi(.,\lambda_{t+1}^k)$
\begin{align}
    & \inprod{\nabla_x \varphi(x_t^k, \lambda_{t+1}^k)}{x - x_{t+1}^k} = \inprod{\nabla_x \varphi(x_t^k, \lambda_{t+1}^k)}{x - x_t^k} 
    + \inprod{\nabla_x \varphi(x_t^k, \lambda_{t+1}^k)}{x_t^k - x_{t+1}^k} \nonumber \\
    & \leq \varphi(x, \lambda_{t+1}^k) - \varphi(x_t^k, \lambda_{t+1}^k) - \tfrac{\mu_f}{2}\norm{x - x_t^k}^2
    + \inprod{\nabla_x \varphi(x_t^k, \lambda_{t+1}^k)}{x_t^k - x_{t+1}^k}. \label{eq:a67}
\end{align}
Using inequality \eqref{eq:a67} in \eqref{eq:a66}
\begin{align}
	&\inprod{\delta^{x^k}_t}{x_{t+1}^k - x} - \varphi(x, \lambda_{t+1}^k) + \varphi(x_t^k, \lambda_{t+1}^k) \leq \inprod{\nabla_x \varphi(x_t^k, \lambda_{t+1}^k)}{x_t^k - x_{t+1}^k} + (\tfrac{1-\gamma_x^k}{2\tau_{t+1}^k}- \tfrac{\mu_f}{2})\norm{x - x_t^k}^2 \nonumber\\
    & - \tfrac{1-\gamma_x^k}{2\tau_{t+1}^k}\norm{x - x_{t+1}^k}^2 - \tfrac{1-\gamma_x^k}{2\tau_{t+1}^k}\norm{x_{t+1}^k - x_t^k}^2 + \tfrac{\gamma_x^k}{2\tau_{t+1}^k} \big(D_x^k(x) - \norm{x - x_{t+1}^k}^2 - D_x^k(x_{t+1}^k) \big).\label{eq:a69} 
\end{align}
Upper bound of $\inprod{\nabla_x \varphi(x_t^k, \lambda_{t+1}^k)}{x_t^k - x_{t+1}^k}$
\begin{align}
    &\inprod{\nabla_x \varphi(x_t^k, \lambda_{t+1}^k)}{x_t^k - x_{t+1}^k} 
    = \inprod{\nabla f(x_t^k)}{x_t^k - x_{t+1}^k} 
    + \inprod{\lambda_{t+1}^k \nabla g(x_t^k)}{x_t^k - x_{t+1}^k} \nonumber \\
    &= \inprod{\nabla f(x_t^k)}{x_t^k - x_{t+1}^k} 
    + \inprod{\lambda_{t+1}^k }{\nabla g(x_t^k)^\top(x_t^k - x_{t+1}^k)} =  \inprod{\nabla f(x_t^k)}{x_t^k - x_{t+1}^k} 
    + \inprod{\lambda_{t+1}^k }{g(x^k_t)-\ell_g(x_{t+1}^k)} \nonumber \\
    &\le f(x^k_t) - f(x^k_{t+1}) + \tfrac{L_f}{2} \gnorm{x^k_{t+1}- x^k_t}{}{2} + \inprod{\lambda_{t+1}^k }{g(x^k_t)-\ell_g(x_{t+1}^k)} \nonumber \\
    &= \varphi(x^k_t, \lambda^k_{t+1}) - f(x^k_{t+1}) + \tfrac{L_f}{2} \gnorm{x^k_{t+1}- x^k_t}{}{2} - \inprod{\lambda_{t+1}^k }{\ell_g(x_{t+1}^k)}, \label{eq:a70}
\end{align}
where in second-last equality we use \eqref{eq:p5}, and the last inequality, \eqref{eq:p2}. Upper bound for \eqref{eq:a69}
\begin{align}
    &\inprod{\delta^{x^k}_t}{x_{t+1}^k - x} - \varphi(x, \lambda_{t+1}^k) + \varphi(x_t^k, \lambda_{t+1}^k) \leq \varphi(x^k_t, \lambda^k_{t+1}) - f(x^k_{t+1}) + \tfrac{L_f}{2} \norm{x^k_{t+1}- x^k_t}^2 - \inprod{\lambda_{t+1}^k }{\ell_g(x_{t+1}^k)} \nonumber\\
    & +  (\tfrac{1-\gamma_x^k}{2\tau_{t+1}^k}- \tfrac{\mu_f}{2})\norm{x - x_t^k}^2 - \tfrac{1-\gamma_x^k}{2\tau_{t+1}^k}\norm{x - x_{t+1}^k}^2 - \tfrac{1-\gamma_x^k}{2\tau_{t+1}^k}\norm{x_{t+1}^k - x_t^k}^2 \nonumber\\
    & + \tfrac{\gamma_x^k}{2\tau_{t+1}^k} \big(D_x^k(x) - \norm{x - x_{t+1}^k}^2 - D_x^k(x_{t+1}^k) \big). \label{eq:a71}
\end{align}
Summing up inequalities \eqref{eq:a71} and \eqref{eq:a64} gives \eqref{eq:a58}
\end{proof}
\begin{lemma}
    The expected gap inequality is as follows:
\begin{align}
     &\tsum_{k=1}^{K}u^kE[Q(\wt{z}^k,z^*)]\leq  \tsum_{k=1}^{K-1}(\alpha_1^{k+1}\theta_1^{k+1} - \alpha_{T^k}^k)E[\inprod{\bar{q}_{T^k}^k}{\lambda_{T^k}^k - \lambda^*}] \nonumber\\ 
     &  + \tsum_{k=1}^{K}\tsum_{t=1}^{T^k-1}(\alpha_{t+1}^k\theta_{t+1}^k - \alpha_t^k)E[\inprod{\bar{q}_t^k}{\lambda_t^k - \lambda^*}]  \nonumber \\
     & + \tsum_{k=1}^{K}\tsum_{t=0}^{T^k-1}(\tfrac{\alpha_{t+1}^k\theta_{t+1}^k\rho_{t+1}^k M_g}{2} - \tfrac{(1 - \gamma_\lambda^k)\alpha_{t+1}^k}{2\sigma_{t+1}^k})E[\norm{\lambda_{t+1}^k - \lambda_t^k}^2] - \tsum_{k=1}^{K}\tsum_{t=0}^{T^k-1}\tfrac{\gamma_\lambda^k\alpha_{t+1}^k}{2\sigma_{t+1}^k}E[D_\lambda^k(\lambda_{t+1}^k)]\nonumber\\ 
     & + \tfrac{(1 - \gamma_\lambda^1)\alpha_1^1}{2\sigma_1^1} E[\norm{\lambda^* - \lambda_0^1}^2] + \tsum_{k=1}^{K-1}(\tfrac{(1 - \gamma_\lambda^{k+1})\alpha_1^{k+1}}{2\sigma_1^{k+1}} - \tfrac{(1 - \gamma_\lambda^k)\alpha_{T^k}^k}{2\sigma_{T^k}^k})E[\norm{\lambda^* - \lambda_{T^k}^k}^2]  \nonumber\\
     & - \tfrac{(1 - \gamma_\lambda^K)\alpha_{T^K}^K}{4\sigma_{T^K}^K}E[\norm{\lambda^* - \lambda_{T^K}^K}^2] + \tsum_{k=1}^{K}\tfrac{1 - \gamma_\lambda^k}{2}\tsum_{t=1}^{T^k-1}(\tfrac{\alpha_{t+1}^k}{\sigma_{t+1}^k} - \tfrac{\alpha_t^k}{\sigma_t^k})E[\norm{\lambda^* - \lambda_t^k}^2] \nonumber \\
     & + s^1E[D_\lambda^1(\lambda^*)] + \tsum_{k=1}^{K-1}(s^{k+1} - s^k)E[D_\lambda^{k+1}(\lambda^*)] - s^KE[D_\lambda^{K+1}(\lambda^*)]  \nonumber\\
     & + \tsum_{k=1}^{K}\tsum_{t=0}^{T^k-1}(\tfrac{\norm{\lambda^*}L_g}{2} + L_f^2\nu_{t+1}^k + \tfrac{1}{2\nu_{t+1}^k} + \tfrac{L_f}{2} + \tfrac{\alpha_{t+2}^k\theta_{t+2}^kM_g}{2\rho_{t+2}^k\alpha_{t+1}^k} - \tfrac{1 - \gamma_x^k}{2\tau_{t+1}^k})\alpha_{t+1}^k E[\norm{x_t^k - x_{t+1}^k}^2] \nonumber\\
     & + (\tfrac{1-\gamma_x^1}{2\tau_1^1} - \tfrac{\mu_f}{2})\alpha_1^1 E[\norm{x^* - x_0^1}^2] + \tsum_{k=1}^{K-1}\big((\tfrac{(1-\gamma_x^{k+1})}{2\tau_1^{k+1}} - \tfrac
     {\mu_f}{2})\alpha_1^{k+1}  - \tfrac{(1-\gamma_x^k)\alpha_{T^k}^k}{2\tau_{T^k}^k}\big)E[\norm{x^* - x_{T^k}^k}^2]    \nonumber \\
     & - \tfrac{(1-\gamma_x^K)\alpha_{T^K}^K}{2\tau_{T^K}^K}E[\norm{x^* - x_{T^K}^K}^2] + \tsum_{k=1}^{K}\tsum_{t=1}^{T^k-1}\big((\tfrac{1-\gamma_x^k}{2\tau_{t+1}^k} - \tfrac{\mu_f}{2})\alpha_{t+1}^k - \tfrac{(1-\gamma_x^k)\alpha_t^k}{2\tau_t^k}\big)E[\norm{x^* - x_t^k}^2] \nonumber\\
     & + r^1 E[D_x^1(x^*)] + \tsum_{k=1}^{K-1}(r^{k+1} - r^k)E[D_x^{k+1}(x^*)] - r^KE[D_x^{K+1}(x^*)]  \nonumber\\
     & + \tfrac{\alpha_1^1\theta_1^1 M_g}{2\rho_1^1}E[\norm{x_0^1 - x_{-1}^1}^2] +  \tsum_{k=1}^{K-1}\big(\tfrac{\alpha_1^{k+1}\theta_1^{k+1} M_g}{2\rho_1^{k+1}} - \tfrac{\alpha_{T^k}^k\theta_{T^k}^k M_g}{2\rho_{T^k}^k}\big)E[\norm{x_{T^k}^k - x_{T^k-1}^k}^2] \nonumber\\
     & + (\tfrac{M_g^2\sigma_{T^K}^K\alpha_{T^K}^K\theta_{T^K}^K}{1-\gamma_{\lambda}^K}- \tfrac{\alpha_{T^K}^K\theta_{T^K}^K M_g}{2\rho_{T^K}^K})E[\norm{x_{T^K}^K - x_{T^K-1}^K}^2] + \tsum_{k=1}^{K}\tsum_{t=0}^{T^k-1}(\tfrac{4L_f^2\nu_{t+1}^k}{2} - \tfrac{\gamma_x^k}{2\tau_{t+1}^k})\alpha_{t+1}^k E[D_x^k(x_{t+1}^k)].\label{eq:d2}
\end{align}
\end{lemma}
\begin{proof}
Starting with \eqref{eq:a58}, Setting $(x,\lambda,z)=(x^*,\lambda^*,z^*)$, taking expectation on both sides of the inequality and using Lemma \ref{eq:C1}, multiplying both sides by $\alpha_{t+1}^k$, summing over index $t$ and accumulating the similar terms: 
\begin{align}
    &\tsum_{t=0}^{T^k-1}\alpha_{t+1}^k E[Q(z_{t+1}^k, z^*)] \leq  \alpha_1^k\theta_1^kE[\inprod{\bar{q}_0^k}{\lambda_0^k - \lambda^*}] + \tsum_{t=1}^{T^k-1}(\alpha_{t+1}^k\theta_{t+1}^k - \alpha_t^k)E[\inprod{\bar{q}_t^k}{\lambda_t^k - \lambda^*}] \nonumber \\
    & - \alpha_{T^k}^kE[\inprod{\bar{q}_{T^k}^k}{\lambda_{T^k}^k - \lambda^*}] + \tsum_{t=0}^{T^k-1}(\tfrac{\alpha_{t+1}^k\theta_{t+1}^k\rho_{t+1}^kM_g}{2} - \tfrac{1 - \gamma_\lambda^k}{2}\tfrac{\alpha_{t+1}^k}{\sigma_{t+1}^k})E[\norm{\lambda_{t+1}^k - \lambda_t^k}^2] \nonumber \\
    & + \tfrac{1 - \gamma_\lambda^k}{2}\tfrac{\alpha_1^k}{\sigma_1^k}E[\norm{\lambda^* - \lambda_0^k}^2] + \tfrac{1 - \gamma_\lambda^k}{2}\tsum_{t=1}^{T^k-1}(\tfrac{\alpha_{t+1}^k}{\sigma_{t+1}^k} - \tfrac{\alpha_t^k}{\sigma_t^k})E[\norm{\lambda^* - \lambda_t^k}^2] - \tfrac{1 - \gamma_\lambda^k}{2}\tfrac{\alpha_{T^k}^k}{\sigma_{T^k}^k}E[\norm{\lambda^* - \lambda_{T^k}^k}^2] \nonumber\\
    &  +  s^k E[D_\lambda^k(\lambda^*)] 
    -  s^k E[D_\lambda^{k+1}(\lambda^*)] 
    - \tsum_{t=0}^{T^k-1}\tfrac{\gamma_\lambda^k\alpha_{t+1}^k}{2\sigma_{t+1}^k}E[D_\lambda^k(\lambda_{t+1}^k)]  \nonumber \\
    & + \tsum_{t=0}^{T^k-1}(\tfrac{\norm{\lambda^*}L_g}{2} + L_f^2\nu_{t+1}^k + \tfrac{1}{2\nu_{t+1}^k} + \tfrac{L_f}{2} + \tfrac{\alpha_{t+2}^k\theta_{t+2}^kM_g}{2\rho_{t+2}^k\alpha_{t+1}^k} - \tfrac{1 - \gamma_x^k}{2\tau_{t+1}^k})\alpha_{t+1}^k E[\norm{x_t^k - x_{t+1}^k}^2] \nonumber \\
    & + \tfrac{\alpha_1^k\theta_1^kM_g}{2\rho_1^k}E[\norm{x_0^k - x_{-1}^k}^2] - \tfrac{\alpha_{T^k}^k\theta_{T^k}^kM_g}{2\rho_{T^k}^k}E[\norm{x_{T^k}^k - x_{T^k-1}^k}^2] + \tsum_{t=0}^{T^k-1}(2L_f^2\nu_{t+1}^k - \tfrac{\gamma_x^k}{2\tau_{t+1}^k})\alpha_{t+1}^kE[D_x^k(x_{t+1}^k)]\nonumber \\
    & + (\tfrac{1-\gamma_x^k}{2\tau_1^k} - \tfrac{\mu_f}{2})\alpha_1^kE[\norm{x^* - x_0^k}^2] + \tsum_{t=1}^{T^k-1}\big(\tfrac{1-\gamma_x^k}{2\tau_{t+1}^k}\alpha_{t+1}^k  - \tfrac{1-\gamma_x^k}{2\tau_t^k}\alpha_t^k - \tfrac{\mu_f}{2}\alpha_{t+1}^k\big)E[\norm{x^* - x_t^k}^2]  \nonumber\\
    & - \tfrac{\alpha_{T^k}^k}{\tau_{T^k}^k}\tfrac{1-\gamma_x^k}{2}E[\norm{x^* - x_{T^k}^k}^2] + \tsum_{t=0}^{T^k-1}\tfrac{\gamma_x^k\alpha_{t+1}^k}{2\tau_{t+1}^k} \big( E[D_x^k(x^*)]  - E[\norm{x^* - x_{t+1}^k}^2] \big). \label{eq:a72}
\end{align}
summing up over index $k$ after using \eqref{eq:p13}, the result will be \eqref{eq:d2}. 
\end{proof}
Let $m_{t+1}^k  = \tfrac{\norm{\lambda^*}L_g + L_f + 2L_f^2\nu_{t+1}^k}{2} + \tfrac{1}{2\nu_{t+1}^k} + \tfrac{\alpha_{t+2}^k\theta_{t+2}^kM_g}{2\rho_{t+2}^k\alpha_{t+1}^k} - \tfrac{1 - \gamma_x^k}{2\tau_{t+1}^k}$, and $n_{t+1}^k = \tfrac{1-\gamma_x^k}{2\tau_{t+1}^k}\alpha_{t+1}^k  - \tfrac{1-\gamma_x^k}{2\tau_t^k}\alpha_t^k - \tfrac{\mu_f}{2}\alpha_{t+1}^k$.
\begin{remark} \label{eq:l3}
    The following conditions are necessary for the convergence.
    \begin{align*}
        \alpha_{t+1}^k\theta_{t+1}^k - \alpha_t^k &=0 &\alpha_1^{k+1}\theta_1^{k+1} - \alpha_{T^k}^k & = 0 &m_{t+1}^k &\leq 0 & n_{t+1}^k &\leq 0  & \tfrac{\alpha_{t+1}^k}{\sigma_{t+1}^k} - \tfrac{\alpha_t^k}{\sigma_t^k} &\leq 0,
    \end{align*}
    \begin{align*}
        2\nu_{t+1}^kL_f^2 - \tfrac{\gamma_x^k}{2\tau_{t+1}^k} &\leq 0 &\tfrac{\alpha_{t+1}^k\theta_{t+1}^k\rho_{t+1}^k M_g}{2} - \tfrac{(1 - \gamma_\lambda^k)\alpha_{t+1}^k}{2\sigma_{t+1}^k} & \leq 0  & \tfrac{M_g^2\sigma_{T^K}^K\alpha_{T^K}^K}{1-\gamma_{\lambda}^K}- \tfrac{\alpha_{T^K}^K\theta_{T^K}^K M_g}{2\rho_{T^K}^K} &\leq 0.
    \end{align*}
\end{remark}
The condition $2\nu_{t+1}^kL_f^2 - \tfrac{\gamma_x^k}{2\tau_{t+1}^k} \leq 0$ may not be satisfied for the first epochs or iterations but by passing through more epochs, it will be satisfied and since $E[D_x^k(x_{t+1}^k)]$ is bounded we can upper bound the first terms that this condition is not satisfied with some number but for the sake of simplicity, we assume that this condition will be satisfied from the beginning.
\begin{theorem} \label{eq:l2}
    Assume that the following setup for parameters holds:
    \begin{align*}
        A_t^k &=t_0 + \tfrac{k(k-1)}{2} + t &\alpha_t^k &= \tfrac{w}{\mu_f}A_t^k & \theta_t^k &= \tfrac{\alpha_{t-1}^k}{\alpha_t^k} & \rho_t^k &= \tfrac{4M_g w}{\mu_f(A_t^k - t_0)}   &\sigma_t^k &= \tfrac{\mu_f(A_t^k - t_0)}{6M_g^2wk}, 
    \end{align*}
    \begin{align*}
         \gamma_{\lambda}^k &= \tfrac{1}{k} &\nu_t^k &= \tfrac{4w}{\mu_f(A_t^k - t_0)} & \tau_t^k &= \tfrac{w}{2\mu_fA_t^k} & \gamma_x^k &= \tfrac{1}{k} & T^k = k,
    \end{align*}
    where $t_0 = \tfrac{4w}{\mu_f}(\tfrac{\norm{\lambda^*}L_g + L_f + 3L_f^2 + 4}{2})$, and $w \geq 9$. The gradient computation for objective function to reach $\epsilon$-gap will be $\frac{m}{\epsilon^{1/3}} + \tfrac{1}{\epsilon^{2/3}}$.
\end{theorem}
\begin{proof}
It can be verified that the parameters in Theorem \ref{eq:l2} satisfy the conditions in Remark \ref{eq:l3}. $n_t^k$ is a negative coefficient using parameters of Theorem \ref{eq:l2} but this negative number could be used to eliminate the positive term $r^{k+1} - r^k$. The cumulative coefficient of $E[\norm{x^* - x_t^k}^2]$
\begin{align}
    \tfrac{r^{k+1} - r^k}{r^{k+1}}\tfrac{\gamma_x^k\alpha_t^k}{2\tau_t^k} + n_t^k \leq \tfrac{r^{k+1} - r^k}{r^{k+1}}\tfrac{A_t^k}{k} + \tfrac{k-1}{k}(2A_t^k - 1) - \tfrac{w}{2}A_t^k \leq 0. \label{eq:a90}
\end{align}
 \eqref{eq:a90} holds if $w \geq 9$. New gap inequality
\begin{align}
     &\tsum_{k=1}^{K}u^kE[Q(\wt{z}^k,z^*)] \leq  \tfrac{(1 - \gamma_\lambda^1)\alpha_1^1}{2\sigma_1^1} E[\norm{\lambda^* - \lambda_0^1}^2] + \tsum_{k=1}^{K-1}(\tfrac{(1 - \gamma_\lambda^{k+1})\alpha_1^{k+1}}{2\sigma_1^{k+1}} - \tfrac{(1 - \gamma_\lambda^k)\alpha_{T^k}^k}{2\sigma_{T^k}^k})E[\norm{\lambda^* - \lambda_{T^k}^k}^2]  \nonumber\\
     & - \tfrac{(1 - \gamma_\lambda^K)\alpha_{T^K}^K}{4\sigma_{T^K}^K}E[\norm{\lambda^* - \lambda_{T^K}^K}^2]  + s^1E[D_\lambda^1(\lambda^*)] + \tsum_{k=1}^{K-1}(s^{k+1} - s^k)E[D_\lambda^{k+1}(\lambda^*)] - s^KE[D_\lambda^{K+1}(\lambda^*)] \nonumber\\
     & + (\tfrac{1-\gamma_x^1}{2\tau_1^1} - \tfrac{\mu_f}{2})\alpha_1^1 E[\norm{x^* - x_0^1}^2] + r^1 E[D_x^1(x^*)] + \tfrac{\alpha_1^1\theta_1^1 M_g}{2\rho_1^1}E[\norm{x_0^1 - x_{-1}^1}^2]. \label{eq:c92}
\end{align}
Considering $p = \min\{\tfrac{(1 - \gamma_\lambda^K)\alpha_{T^K}^K}{4\sigma_{T^K}^K},s^K\}$ and $q = \max_{k \in \{1,\dots,K-1\}}\{\tfrac{(1 - \gamma_\lambda^{k+1})\alpha_1^{k+1}}{2\sigma_1^{k+1}} - \tfrac{(1 - \gamma_\lambda^k)\alpha_{T^k}^k}{2\sigma_{T^k}^k},s^{k+1} - s^k\}$,
and since $\tsum_{k=1}^{K}u^kE[Q(\wt{z}^k,z^*)] \geq 0$ 
\begin{align}
     & E[\norm{\lambda^* - \lambda_{T^K}^K}^2] + E[D_\lambda^{K+1}(\lambda^*)] \leq \tfrac{1}{p}\bigg(\tfrac{(1 - \gamma_\lambda^1)\alpha_1^1}{2\sigma_1^1} E[\norm{\lambda^* - \lambda_0^1}^2] + q \tsum_{k=1}^{K-1}(E[\norm{\lambda^* - \lambda_{T^k}^k}^2] + E[D_\lambda^{k+1}(\lambda^*)])  \nonumber\\
     & + \gamma_\lambda^1s^1E[D_\lambda^1(\lambda^*)] + (\tfrac{1-\gamma_x^1}{2\tau_1^1} - \tfrac{\mu_f}{2})\alpha_1^1 E[\norm{x^* - x_0^1}^2] + r^1 E[D_x^1(x^*)] + \tfrac{\mu_f t_0 }{8M_gw}E[\norm{x_0^1 - x_{-1}^1}^2]\bigg). \label{eq:a93}
\end{align}
Naming $a = \tfrac{(1 - \gamma_\lambda^1)\alpha_1^1}{2\sigma_1^1} E[\norm{\lambda^* - \lambda_0^1}^2] + \gamma_\lambda^1s^1E[D_\lambda^1(\lambda^*)] + (\tfrac{1-\gamma_x^1}{2\tau_1^1} - \tfrac{\mu_f}{2})\alpha_1^1 E[\norm{x^* - x_0^1}^2] + r^1 E[D_x^1(x^*)] + \tfrac{\alpha_1^1\theta_1^1 M_g}{4\rho_1^1}E[\norm{x_0^1 - x_{-1}^1}^2]$ and using Lemma \ref{eq:a35} 
\begin{align}
      E[\norm{\lambda^* - \lambda_{T^K}^K}^2] + E[D_\lambda^{K+1}(\lambda^*)] \leq \tfrac{1}{p}\big(a + q \tsum_{k=1}^{K-1}(E[\norm{\lambda^* - \lambda_{T^k}^k}^2] + E[D_\lambda^{k+1}(\lambda^*)])\big) \leq \tfrac{a}{p}(1 + \tfrac{q}{p})^K. \label{eq:a94}
\end{align}
Finally $\lim_{K \rightarrow \infty} \tfrac{a}{p}(1 + \tfrac{q}{p})^K = 0$ and $\tfrac{qK}{p} \leq 2 + t_0$ so $E[\norm{\lambda^* - \lambda_{T^K}^K}^2] + E[D_\lambda^{K+1}(\lambda^*)] \leq a e^{t_0+2}$. This proves a constant bound for the Lagrangian terms. Back to \eqref{eq:c92}, using \eqref{eq:p13}, dividing both sides of the inequality by $U^K$, the final expected gap  
\begin{align}
     &E[Q(\wt{Z}^K,z^*)] \leq  \tfrac{1}{U^K}\bigg(2qK e^{t_0 + 2} + t_0^2 D_x^1(x^*) + \tfrac{\alpha_1^1\theta_1^1 M_g}{4\rho_1^1}E[\norm{x_0^1 - x_{-1}^1}^2] \bigg) = \epsilon. \label{eq:a95} 
\end{align}
\eqref{eq:a95} results in $K = O(\tfrac{1}{\epsilon^{1/3}})$. The gradient computations for the objective function are $Km + \tsum_{k=1}^{K}T^k = \tfrac{m}{\epsilon^{1/3}} + \tfrac{1}{\epsilon^{2/3}}$ to reach $\epsilon$-gap. This concludes the proof of theorem \ref{eq:l2}. 
\end{proof}

\subsection{Fully-finite-sum case.} \label{sec-full}
The algorithm for this case is Algorithm \ref{alg:alg3}, and the analysis to reach a gap inequality is identical to the semi-finite-sum case. Let start from \eqref{eq:a58}. We use upper bound \eqref{eq:b24} for $\inprod{\delta_t^{x^k}}{x - x_{t+1}^k}$ to find a constant upper bound for some Lagrangian expressions. 
\begin{align}
     &Q(z_{t+1}^k, z) \leq \inprod{\delta_t^{\lambda^k}}{\lambda_{t+1}^k - \lambda} - \inprod{\bar{q}_{t+1}^k}{\lambda_{t+1}^k - \lambda} 
    + \inprod{\theta_{t+1}^k \bar{q}_t^k}{\lambda_t^k - \lambda} + \inprod{\theta_{t+1}^k(q_t^k - \bar{q}_t^k)}{\lambda_t^k - \lambda} \nonumber \\
    & + \tfrac{\theta_{t+1}^k M_g}{2\rho_{t+1}^k}\norm{x_t^k - x_{t-1}^k}^2 + \tfrac{\theta_{t+1}^k\rho_{t+1}^k M_g}{2}\norm{\lambda_{t+1}^k - \lambda_t^k}^2  + \tfrac{1 - \gamma_\lambda^k}{2\sigma_{t+1}^k} \big( \norm{\lambda - \lambda_t^k}^2 - \norm{\lambda - \lambda_{t+1}^k}^2 - \norm{\lambda_t^k - \lambda_{t+1}^k}^2 \big) \nonumber \\
    & + \tfrac{\gamma_\lambda^k}{2\sigma_{t+1}^k} \big( D_\lambda^k(\lambda)         - \norm{\lambda - \lambda_{t+1}^k}^2 - D_\lambda^k(\lambda_{t+1}^k) \big) + M_g\nu_{t+1}^kD_{\lambda}^k(\lambda_{t+1}^k) + 2L_gD_x\nu_{t+1}^kD_{\lambda}^k(\lambda) + 2L_gD_x\nu_{t+1}^k\norm{\lambda}^2\nonumber\\ 
    &  + (\tfrac{\norm{\lambda}L_g + 5L_f}{2} + \tfrac{M_g + 2L_gD_x + 2L_f}{\nu_{t+1}^k})\norm{x_t^k - x_{t+1}^k}^2 + L_f\nu_{t+1}^kD_x^k(x_{t+1}^k) + (\tfrac{1-\gamma_x^k}{2\tau_{t+1}^k} - \tfrac{\mu_f}{2})\norm{x - x_t^k}^2 \nonumber\\
    & - \tfrac{1-\gamma_x^k}{2\tau_{t+1}^k}\norm{x - x_{t+1}^k}^2 - \tfrac{1-\gamma_x^k}{2\tau_{t+1}^k}\norm{x_{t+1}^k - x_t^k}^2  + \tfrac{\gamma_x^k}{2\tau_{t+1}^k} \big( D_x^k(x) - \norm{x - x_{t+1}^k}^2 
    - D_x^k(x_{t+1}^k) \big). \label{eq:a80}
\end{align}
Notice that, in \eqref{eq:a80}, We replaced $\eta^k$ with $\nu_{t+1}^k$ for the upper bound \eqref{eq:b24}. 
\begin{lemma}
    The final expected gap will be:
\begin{align}
     &\tsum_{k=1}^{K}u^kE[Q(\wt{z}^k,z^*)]\leq  \alpha_1^1\theta_1^1E[\inprod{\bar{q}_0^1}{\lambda_0^1 - \lambda^*}] + \tsum_{k=1}^{K-1}(\alpha_1^{k+1}\theta_1^{k+1} - \alpha_{T^k}^k)E[\inprod{\bar{q}_{T^k}^k}{\lambda_{T^k}^k - \lambda^*}] \nonumber\\ 
     & - \alpha_{T^K}^KE[\inprod{\bar{q}_{T^K}^K}{\lambda_{T^K}^K - \lambda^*}] +\tsum_{k=1}^{K}\tsum_{t=1}^{T^k-1}(\alpha_{t+1}^k\theta_{t+1}^k - \alpha_t^k)E[\inprod{\bar{q}_t^k}{\lambda_t^k - \lambda^*}] + \tsum_{k=1}^{K}v^k\norm{\lambda^*}^2 \nonumber \\
     & + \tfrac{\alpha_1^1\theta_1^1}{2\eta_1^1}E[\norm{\lambda_0^1 -\lambda_{-1}^1}^2] + \tsum_{k=1}^{K-1}(\tfrac{\alpha_1^{k+1}\theta_1^{k+1}}{2\eta_1^{k+1}} - \tfrac{\alpha_{T^k}^k\theta_{T^k}^k}{2\eta_{T^k}^k})E[\norm{\lambda_{T^k}^k - \lambda_{T^k-1}^k}^2] - \tfrac{\alpha_{T^K}^K\theta_{T^K}^K}{2\eta_{T^K}^K}E[\norm{\lambda_{T^K}^K - \lambda_{T^K-1}^K}^2]\nonumber\\ 
     & + \tsum_{k=1}^{K}\tsum_{t=0}^{T^k-1}(\tfrac{\alpha_{t+1}^k\theta_{t+1}^k\rho_{t+1}^k M_g}{2} + \tfrac{\alpha_{t+1}^k}{2\omega_{t+1}^k} + \tfrac{\alpha_{t+2}^k\theta_{t+2}^k}{2\eta_{t+2}^k} - \tfrac{(1 - \gamma_\lambda^k)\alpha_{t+1}^k}{2\sigma_{t+1}^k})E[\norm{\lambda_{t+1}^k - \lambda_t^k}^2] \nonumber\\ 
     & + \tfrac{(1 - \gamma_\lambda^1)\alpha_1^1}{2\sigma_1^1} E[\norm{\lambda^* - \lambda_0^1}^2] + \tsum_{k=1}^{K-1}(\tfrac{(1 - \gamma_\lambda^{k+1})\alpha_1^{k+1}}{2\sigma_1^{k+1}} - \tfrac{(1 - \gamma_\lambda^k)\alpha_{T^k}^k}{2\sigma_{T^k}^k})E[\norm{\lambda^* - \lambda_{T^k}^k}^2]  \nonumber\\
     & - \tfrac{(1 - \gamma_\lambda^K)\alpha_{T^K}^K}{4\sigma_{T^K}^K}E[\norm{\lambda^* - \lambda_{T^K}^K}^2] + \tsum_{k=1}^{K}\tfrac{1 - \gamma_\lambda^k}{2}\tsum_{t=1}^{T^k-1}(\tfrac{\alpha_{t+1}^k}{\sigma_{t+1}^k} - \tfrac{\alpha_t^k}{\sigma_t^k})E[\norm{\lambda^* - \lambda_t^k}^2] \nonumber \\
     & + (s^1 + v^1)E[D_\lambda^1(\lambda^*)] + \tsum_{k=1}^{K-1}(s^{k+1} - s^k + v^{k+1})E[D_\lambda^{k+1}(\lambda^*)] - s^KE[D_\lambda^{K+1}(\lambda^*)] \nonumber\\
     & + \tsum_{k=1}^{K}\tsum_{t=0}^{T^k-1}(M_g\nu_{t+1}^k - \tfrac{\gamma_\lambda^k}{2\sigma_{t+1}^k})\alpha_{t+1}^k E[D_\lambda^k(\lambda_{t+1}^k)] + \tsum_{k=1}^{K}\tsum_{t=0}^{T^k-2}(\tfrac{\norm{\lambda^*}L_g + 5L_f}{2} + \tfrac{M_g^2\alpha_{t+2}^k\theta_{t+2}^k\eta_{t+2}^k}{2\alpha_{t+1}^k} \nonumber\\
     & + M_g^2\omega_{t+1}^k + \tfrac{\alpha_{t+2}^k\theta_{t+2}^k M_g}{2\rho_{t+2}^k\alpha_{t+1}^k} + \tfrac{M_g + 2L_gD_x + L_f}{\nu_{t+1}^k} - \tfrac{(1-\gamma_x^k)}{2\tau_{t+1}^k})\alpha_{t+1}^kE[\norm{x_t^k - x_{t+1}^k}^2] \nonumber\\
     & + (\tfrac{1-\gamma_x^1}{2\tau_1^1} - \tfrac{\mu_f}{2})\alpha_1^1E[\norm{x^* - x_0^1}^2] + \tsum_{k=1}^{K-1}\big((\tfrac{(1-\gamma_x^{k+1})}{2\tau_1^{k+1}} - \tfrac
     {\mu_f}{2})\alpha_1^{k+1}  - \tfrac{(1-\gamma_x^k)\alpha_{T^k}^k}{2\tau_{T^k}^k}\big)E[\norm{x^* - x_{T^k}^k}^2]    \nonumber \\
     & - \tfrac{(1-\gamma_x^K)\alpha_{T^K}^K}{2\tau_{T^K}^K}E[\norm{x^* - x_{T^K}^K}^2] + \tsum_{k=1}^{K}\tsum_{t=1}^{T^k-1}\big((\tfrac{1-\gamma_x^k}{2\tau_{t+1}^k} - \tfrac{\mu_f}{2})\alpha_{t+1}^k - \tfrac{(1-\gamma_x^k)\alpha_t^k}{2\tau_t^k}\big)E[\norm{x^* - x_t^k}^2] \nonumber\\
     & + r^1E[D_x^1(x^*)] + \tsum_{k=1}^{K-1}(r^{k+1} - r^k)E[D_x^{k+1}(x^*)] - r^KE[D_x^{K+1}(x^*)]  \nonumber\\
     & + \tsum_{k=1}^{K}\tsum_{t=0}^{T^k-1}(M_g^2\omega_{t+1}^k + L_f\nu_{t+1}^k - \tfrac{\gamma_x^k}{2\tau_{t+1}^k})\alpha_{t+1}^kE[D_x^k(x_{t+1}^k)] + (\tfrac{M_g^2\alpha_1^1\theta_1^1\eta_1^1}{2} + \tfrac{\alpha_1^1\theta_1^1 M_g}{2\rho_1^1})E[\norm{x_0^1 - x_{-1}^1}^2] \nonumber\\
     & + (\tfrac{M_g^2\sigma_{T^K}^K\alpha_{T^K}^K\theta_{T^K}^K}{1-\gamma_{\lambda}^K} - \tfrac{M_g^2\alpha_{T^K}^K\theta_{T^K}^K\eta_{T^K}^K}{2} - \tfrac{\alpha_{T^K}^K\theta_{T^K}^K M_g}{2\rho_{T^K}^K})E[\norm{x_{T^K}^K - x_{T^K-1}^K}^2] \nonumber\\
     & + \tsum_{k=1}^{K-1}\big(\tfrac{M_g^2\alpha_1^{k+1}\theta_1^{k+1}\eta_1^{k+1}}{2} + \tfrac{\alpha_1^{k+1}\theta_1^{k+1} M_g}{2\rho_1^{k+1}} - \tfrac{M_g^2\alpha_{T^k}^k\theta_{T^k}^k\eta_{T^k}^k}{2} - \tfrac{\alpha_{T^k}^k\theta_{T^k}^k M_g}{2\rho_{T^k}^k}\big)E[\norm{x_{T^k}^k - x_{T^k-1}^k}^2], \label{eq:a82}
\end{align}
where $(\alpha_{T^k+1}^k,\rho_{T^k+1}^k,\eta_{T^k+1}^k,\theta_{T^k+1}^k) = (\alpha_1^{k+1},\rho_1^{k+1},\eta_1^{k+1},\theta_1^{k+1})$.
\end{lemma}
\begin{proof}
Multiplying both sides of \eqref{eq:a80} by $\alpha_{t+1}^k$, putting $(x,\lambda,z)=(x^*,\lambda^*,z^*)$, applying upper bound \eqref{eq:b24}, and taking expectation on both sides of the inequality
\begin{align}
     &\alpha_{t+1}^kE[Q(z_{t+1}^k, z^*)] \leq \alpha_{t+1}^kE[\inprod{\delta_t^{\lambda^k}}{\lambda_{t+1}^k - \lambda_t^k + \lambda_t^k - \lambda^*}] + \alpha_{t+1}^kE[\inprod{\theta_{t+1}^k \bar{q}_t^k}{\lambda_t^k - \lambda^*}] \nonumber \\
     &- \alpha_{t+1}^kE[\inprod{\bar{q}_{t+1}^k}{\lambda_{t+1}^k - \lambda^*}] + \alpha_{t+1}^kE[\inprod{\theta_{t+1}^k(q_t^k - \bar{q}_t^k)}{\lambda_t^k - \lambda_{t-1}^k + \lambda_{t-1}^k - \lambda^*}] + 2L_gD_x\nu_{t+1}^k\alpha_{t+1}^k\norm{\lambda^*}^2 \nonumber\\ 
     & + (\tfrac{\theta_{t+1}^k\rho_{t+1}^k M_g}{2} - \tfrac{(1 - \gamma_\lambda^k)}{2\sigma_{t+1}^k})\alpha_{t+1}^kE[\norm{\lambda_{t+1}^k - \lambda_t^k}^2] + \tfrac{(1 - \gamma_\lambda^k)\alpha_{t+1}^k}{2\sigma_{t+1}^k} \big( E[\norm{\lambda^* - \lambda_t^k}^2] - E[\norm{\lambda^* - \lambda_{t+1}^k}^2]\big) \nonumber \\
     & + (\tfrac{\gamma_\lambda^k}{2\sigma_{t+1}^k} + 2L_gD_x\nu_{t+1}^k)\alpha_{t+1}^kE[D_\lambda^k(\lambda^*)] - \tfrac{\gamma_\lambda^k\alpha_{t+1}^k}{2\sigma_{t+1}^k}E[\norm{\lambda^* - \lambda_{t+1}^k}^2] + (M_g \nu_{t+1}^k - \tfrac{\gamma_\lambda^k}{2\sigma_{t+1}^k})\alpha_{t+1}^kE[D_\lambda^k(\lambda_{t+1}^k)] \nonumber \\
     & + \big(\tfrac{\norm{\lambda^*}L_g + 5L_f}{2}  + \tfrac{M_g + 2L_gD_x + L_f}{\nu_{t+1}^k} - \tfrac{(1-\gamma_x^k)}{2\tau_{t+1}^k}\big)\alpha_{t+1}^kE[\norm{x_t^k - x_{t+1}^k}^2] \nonumber\\
     & + \tfrac{\alpha_{t+1}^k\theta_{t+1}^k M_g}{2\rho_{t+1}^k}E[\norm{x_t^k - x_{t-1}^k}^2] + (\tfrac{1-\gamma_x^k}{2\tau_{t+1}^k} - \tfrac{\mu_f}{2})\alpha_{t+1}^kE[\norm{x^* - x_t^k}^2] - \tfrac{(1-\gamma_x^k)\alpha_{t+1}^k}{2\tau_{t+1}^k}E[\norm{x^* - x_{t+1}^k}^2]  \nonumber \\
     & + \tfrac{\gamma_x^k\alpha_{t+1}^k}{2\tau_{t+1}^k} \big( E[D_x^k(x^*)] 
     - E[\norm{x^* - x_{t+1}^k}^2] 
     \big) + (L_f \nu_{t+1}^k  - \tfrac{\gamma_x^k}{2\tau_{t+1}^k})\alpha_{t+1}^kE[D_x^k(x_{t+1}^k)]. \label{eq:a81}
\end{align}
 Applying \eqref{eq:p9}, and \eqref{eq:p10} for $\inprod{\delta_t^{\lambda^k}}{\lambda_{t+1}^k - \lambda_t^k}$ and $\inprod{\theta_{t+1}^k(q_t^k - \bar{q}_t^k)}{\lambda_t^k - \lambda_{t-1}^k}$, summing up with respect to the index $t$ and applying the upper bounds in \eqref{eq:a38}, rearranging the terms, and using \eqref{eq:p13}
\begin{align}
     &u^kE[Q(\wt{z}^k,z^*)] \leq  \alpha_1^k\theta_1^kE[\inprod{\bar{q}_0^k}{\lambda_0^k - \lambda^*}] + \tsum_{t=1}^{T^k-1}(\alpha_{t+1}^k\theta_{t+1}^k - \alpha_t^k)E[\inprod{\bar{q}_t^k}{\lambda_t^k - \lambda^*}] \nonumber \\
     & - \alpha_{T^k}^kE[\inprod{\bar{q}_{T^k}^k}{\lambda_{T^k}^k - \lambda^*}] + \tfrac{\alpha_1^k\theta_1^k}{2\eta_1^k}E[\norm{\lambda_0^k - \lambda_{-1}^k}^2] - \tfrac{\alpha_{T^k}^k\theta_{T^k}^k}{2\eta_{T^k}^k}E[\norm{\lambda_{T^k}^k - \lambda_{T^k-1}^k}^2] + v^k\norm{\lambda^*}^2\nonumber\\ 
     & + \tsum_{t=0}^{T^k-1}(\tfrac{\alpha_{t+1}^k\theta_{t+1}^k\rho_{t+1}^k M_g}{2} + \tfrac{\alpha_{t+1}^k}{2\omega_{t+1}^k} + \tfrac{\alpha_{t+2}^k\theta_{t+2}^k}{2\eta_{t+2}^k} - \tfrac{(1 - \gamma_\lambda^k)\alpha_{t+1}^k}{2\sigma_{t+1}^k})E[\norm{\lambda_{t+1}^k - \lambda_t^k}^2] \nonumber\\ 
     & + \tfrac{(1 - \gamma_{\lambda}^k)\alpha_1^k}{2\sigma_1^k}E[\norm{\lambda^* - \lambda_0^k}^2] + \tfrac{1 - \gamma_{\lambda}^k}{2}\tsum_{t=1}^{T^k-1}(\tfrac{\alpha_{t+1}^k}{\sigma_{t+1}^k} - \tfrac{\alpha_t^k}{\sigma_t^k})E[\norm{\lambda^* - \lambda_t^k}^2] - \tfrac{(1 - \gamma_{\lambda}^k)\alpha_{T^k}^k}{2\sigma_{T^k}^k}E[\norm{\lambda^* - \lambda_{T^k}^k}^2] \nonumber \\
     & + (s^k + v^k)E[D_\lambda^k(\lambda^*)] - s^kE[D_\lambda^{k+1}(\lambda^*)] + \tsum_{t=0}^{T^k-1}(M_g\nu_{t+1}^k - \tfrac{\gamma_\lambda^k}{2\sigma_{t+1}^k})\alpha_{t+1}^kE[D_\lambda^k(\lambda_{t+1}^k)] \nonumber \\
     & + \tsum_{t=0}^{T^k-1}\big(\tfrac{\norm{\lambda^*}L_g + 5L_f}{2} + \tfrac{M_g^2\alpha_{t+2}^k\theta_{t+2}^k\eta_{t+2}^k}{2\alpha_{t+1}^k} + \tfrac{\alpha_{t+2}^k\theta_{t+2}^k M_g}{2\rho_{t+2}^k\alpha_{t+1}^k}  \nonumber\\
     & + M_g^2\omega_{t+1}^k + \tfrac{M_g + 2L_gD_x + L_f}{\nu_{t+1}^k} - \tfrac{(1-\gamma_x^k)}{2\tau_{t+1}^k}\big)\alpha_{t+1}^kE[\norm{x_t^k - x_{t+1}^k}^2] \nonumber\\
     & + (\tfrac{1-\gamma_x^k}{2\tau_1^k} - \tfrac{\mu_f}{2})\alpha_1^kE[\norm{x^* - x_0^k}^2] + \tsum_{t=1}^{T^k-1}\big((\tfrac{1-\gamma_x^k}{2\tau_{t+1}^k} - \tfrac{\mu_f}{2})\alpha_{t+1}^k - \tfrac{(1-\gamma_x^k)\alpha_t^k}{2\tau_t^k}\big)E[\norm{x^* - x_t^k}^2] \nonumber\\
     & - \tfrac{(1-\gamma_x^k)\alpha_{T^k}^k}{2\tau_{T^k}^k}E[\norm{x^* - x_{T^k}^k}^2]  + r^kE[D_x^k(x^*)] - r^kE[D_x^{k+1}(x^*)] \nonumber\\
     & + \tsum_{t=0}^{T^k-1}(M_g^2\omega_{t+1}^k + L_f\nu_{t+1}^k - \tfrac{\gamma_x^k}{2\tau_{t+1}^k})\alpha_{t+1}^kE[D_x^k(x_{t+1}^k)]  \nonumber\\
     & + (\tfrac{M_g^2\alpha_1^k\theta_1^k\eta_1^k}{2} + \tfrac{\alpha_1^k\theta_1^k M_g}{2\rho_1^k})E[\norm{x_0^k - x_{-1}^k}^2] - (\tfrac{M_g^2\alpha_{T^k}^k\theta_{T^k}^k\eta_{T^k}^k}{2} + \tfrac{\alpha_{T^k}^k\theta_{T^k}^k M_g}{2\rho_{T^k}^k})E[\norm{x_{T^k}^k - x_{T^k-1}^k}^2], 
\end{align}
where $v^k = \sum_{t=0}^{T^k - 1}2L_gD_x\nu_{t+1}^k\alpha_{t+1}^k$. Summing up the one-epoch inequality over all epochs and using \eqref{eq:p15} gives \eqref{eq:a82}.
\end{proof}
\begin{remark} \label{eq:l4}
     The following conditions for step sizes are necessary to find a constant upper bound for the expected Lagrangian terms of Algorithm \ref{alg:alg3} using the gap inequality \eqref{eq:a82}:
    \begin{align*}
        \alpha_1^{k+1}\theta_1^{k+1} - \alpha_{T^k}^k &= 0 & \alpha_{t+1}^k\theta_{t+1}^k - \alpha_t^k &= 0 & \tfrac{\alpha_1^{k+1}\theta_1^{k+1}}{2\eta_1^{k+1}} - \tfrac{\alpha_{T^k}^k\theta_{T^k}^k}{2\eta_{T^k}^k} & \leq 0,
    \end{align*}
    
    \begin{align*}
        \tfrac{\theta_{t+1}^k\rho_{t+1}^k M_g}{2} + \tfrac{1}{2\omega_{t+1}^k} + \tfrac{\alpha_{t+2}^k\theta_{t+2}^k}{2\eta_{t+2}^k\alpha_{t+1}^k} - \tfrac{1 - \gamma_\lambda^k}{2\sigma_{t+1}^k} & \leq 0 & \tfrac{\alpha_{t+1}^k}{\sigma_{t+1}^k} - \tfrac{\alpha_t^k}{\sigma_t^k} & \leq 0 & M_g\nu_{t+1}^k - \tfrac{\gamma_\lambda^k}{2\sigma_{t+1}^k} \leq 0, 
    \end{align*}
    \begin{align}
        \tfrac{\norm{\lambda^*}L_g + 5L_f}{2} + \tfrac{M_g^2\alpha_{t+2}^k\theta_{t+2}^k\eta_{t+2}^k}{2\alpha_{t+1}^k} + \tfrac{\alpha_{t+2}^k\theta_{t+2}^k M_g}{2\rho_{t+2}^k\alpha_{t+1}^k} + M_g^2\omega_{t+1}^k + \tfrac{M_g + 2L_gD_x + L_f}{\nu_{t+1}^k} - \tfrac{(1-\gamma_x^k)}{2\tau_{t+1}^k} = m_{t+1}^k \leq 0, \label{eq:c13}
    \end{align}
    
    \begin{align*}
        (\tfrac{(1-\gamma_x^{k+1})}{2\tau_1^{k+1}} - \tfrac
     {\mu_f}{2})\alpha_1^{k+1}  - \tfrac{(1-\gamma_x^k)\alpha_{T^k}^k}{2\tau_{T^k}^k} & \leq 0 & (\tfrac{1-\gamma_x^k}{2\tau_{t+1}^k} - \tfrac{\mu_f}{2})\alpha_{t+1}^k - \tfrac{(1-\gamma_x^k)\alpha_t^k}{2\tau_t^k} & \leq 0,  
    \end{align*}
    \begin{align*}
         M_g^2\omega_{t+1}^k + L_f\nu_{t+1}^k - \tfrac{\gamma_x^k}{2\tau_{t+1}^k} & \leq 0 & \tfrac{M_g^2\sigma_{T^K}^K\alpha_{T^K}^K\theta_{T^K}^K}{1-\gamma_{\lambda}^K} - \tfrac{M_g^2\alpha_{T^K}^K\theta_{T^K}^K\eta_{T^K}^K}{2} - \tfrac{\alpha_{T^K}^K\theta_{T^K}^K M_g}{2\rho_{T^K}^K} & \leq 0.
    \end{align*}
\end{remark}
\begin{theorem} \label{eq:th6}
    If the following set up for the parameters holds:
    \begin{align*}
        A_t^k &=t_0 + \tfrac{k(k-1)}{2} + t &\alpha_t^k &= \tfrac{w}{\mu_f}A_t^k & \theta_t^k &= \tfrac{\alpha_{t-1}^k}{\alpha_t^k} & \rho_t^k &= \tfrac{4M_g w}{\mu_f(A_t^k - t_0)}   &\sigma_t^k &= \tfrac{\mu_f(A_t^k - t_0)}{32C_{\nu}M_gwk^2}, 
    \end{align*}
    \begin{align*}
         \gamma_{\lambda}^k &= \tfrac{1}{k+1} &\nu_t^k &= \tfrac{16C_{\nu}w}{\mu_f(A_t^k - t_0)} & \tau_t^k &= \tfrac{w}{2\mu_fA_t^k} & \gamma_x^k &= \tfrac{1}{k+1} & T^k = k,
    \end{align*}
    \begin{align*}
        \eta_t^k & = \tfrac{\mu_f(A_t^k - t_0)}{8M_g^2w} & \omega_t^k &= \tfrac{\mu_f(A_t^k - t_0)}{2M_g^2wk} ,
    \end{align*}
    where $t_0 = \tfrac{w(\norm{\lambda^*}L_g + L_f + 8L_fD_xD_{\lambda}M + 2L_f^2)}{2\mu_f}$, $C_{\nu} = M_g + 2L_gD_x + L_f$ and $w \geq 9$. With this setup, $E[\norm{\lambda^* - \lambda_{T^K}^K}^2]$ and $E[D_\lambda^{K+1}(\lambda^*)]$ will have an upper bound:
    \begin{align*}
        D_{\lambda}^2 := \sup_{K \geq 1} \tfrac{2R}{q^K}.
    \end{align*}  
    where
    \begin{align*}
        &R = \tfrac{\alpha_1^1\theta_1^1}{2\eta_1^1}E[\norm{\lambda_0^1 -\lambda_{-1}^1}^2] + \tsum_{k=1}^{K}v^k\norm{\lambda^*}^2 + \tfrac{(1 - \gamma_\lambda^1)\alpha_1^1}{2\sigma_1^1} E[\norm{\lambda^* - \lambda_0^1}^2] + (s^1 + v^1)E[D_\lambda^1(\lambda^*)]  \nonumber\\
        & + (\tfrac{1-\gamma_x^1}{2\tau_1^1} - \tfrac{\mu_f}{2})\alpha_1^1E[\norm{x^* - x_0^1}^2] + r^1E[D_x^1(x^*)] + (\tfrac{M_g^2\alpha_1^1\theta_1^1\eta_1^1}{2} + \tfrac{\alpha_1^1\theta_1^1 M_g}{2\rho_1^1})E[\norm{x_0^1 - x_{-1}^1}^2],
     \end{align*}
    and $q^K=\min\{\frac{16C_\nu M_gw^2}{\mu_f^2}\cdot\frac{K(K+1)}{2}(1+\frac{2t_0}{K(K+1)}),\frac{16C_\nu M_gw^2}{\mu_f^2}(K+1)(K+1+t_0\sum_{t=1}^{K+1}\frac{1}{\frac{K(K+1)}{2}+t})\}$.
\end{theorem}
\begin{proof}
The paramerters in Theorem \ref{eq:th6} satisfy the conditions in Remark \ref{eq:l4}. Using the parameters in Theorem \ref{eq:th6} 
the cumulative coefficient of $E[\norm{x^* - x_t^k}^2]$
\begin{align}
    \tfrac{r^{k+1} - r^k}{r^{k+1}}\tfrac{\gamma_x^k\alpha_t^k}{2\tau_t^k} + n_t^k \leq \tfrac{r^{k+1} - r^k}{r^{k+1}}\tfrac{A_t^k}{k} + \tfrac{k-1}{k}(2A_t^k - 1) - \tfrac{w}{2}A_t^k \leq 0. \label{eq:c8}
\end{align}
Using $w \geq 9$ in \eqref{eq:c8} could satisfy the inequality sign.
 $\tfrac{M_g^2\alpha_1^{k+1}\theta_1^{k+1}\eta_1^{k+1}}{2} + \tfrac{\alpha_1^{k+1}\theta_1^{k+1} M_g}{2\rho_1^{k+1}} - \tfrac{M_g^2\alpha_{T^k}^k\theta_{T^k}^k\eta_{T^k}^k}{2} - \tfrac{\alpha_{T^k}^k\theta_{T^k}^k M_g}{2\rho_{T^k}^k}$ is positive but $m_{T^k}^k\alpha_{T^k}^k$ could eliminate it. New expected gap inequality
\begin{align}
     &\tsum_{k=1}^{K}u^kE[Q(\wt{z}^k,z^*)] \leq \tfrac{\alpha_1^1\theta_1^1}{2\eta_1^1}E[\norm{\lambda_0^1 -\lambda_{-1}^1}^2] + \tfrac{(1 - \gamma_\lambda^1)\alpha_1^1}{2\sigma_1^1} E[\norm{\lambda^* - \lambda_0^1}^2] + \tsum_{k=1}^{K}v^k\norm{\lambda^*}^2 \nonumber\\ 
     & + \tsum_{k=1}^{K-1}(\tfrac{(1 - \gamma_\lambda^{k+1})\alpha_1^{k+1}}{2\sigma_1^{k+1}} - \tfrac{(1 - \gamma_\lambda^k)\alpha_{T^k}^k}{2\sigma_{T^k}^k})E[\norm{\lambda^* - \lambda_{T^k}^k}^2] - \tfrac{(1 - \gamma_\lambda^K)\alpha_{T^K}^K}{4\sigma_{T^K}^K}E[\norm{\lambda^* - \lambda_{T^K}^K}^2] \nonumber\\
     & + (s^1 + v^1)E[D_\lambda^1(\lambda^*)] + \tsum_{k=1}^{K-1}(s^{k+1} - s^k + v^{k+1})E[D_\lambda^{k+1}(\lambda^*)] - s^KE[D_\lambda^{K+1}(\lambda^*)] \nonumber\\
     & + (\tfrac{1-\gamma_x^1}{2\tau_1^1} - \tfrac{\mu_f}{2})\alpha_1^1E[\norm{x^* - x_0^1}^2] + r^1E[D_x^1(x^*)] + (\tfrac{M_g^2\alpha_1^1\theta_1^1\eta_1^1}{2} + \tfrac{\alpha_1^1\theta_1^1 M_g}{2\rho_1^1})E[\norm{x_0^1 - x_{-1}^1}^2]. \label{eq:a85}
\end{align}
Taking $- \tfrac{(1 - \gamma_\lambda^K)\alpha_{T^K}^K}{4\sigma_{T^K}^K}E[\norm{\lambda^* - \lambda_{T^K}^K}^2]$ and $- s^KE[D_\lambda^{K+1}(\lambda^*)]$ to the left hand side of \eqref{eq:a85}, since $\tsum_{k=1}^{K}u^k E[Q(\wt{z}^k,z^*)] \geq 0 $, by dividing both sides of inequality \eqref{eq:a85} by $q^K = \min\{\tfrac{(1 - \gamma_\lambda^K)\alpha_{T^K}^K}{4\sigma_{T^K}^K},s^K\}$, and taking $p^k = \max\{\tfrac{(1 - \gamma_\lambda^{k+1})\alpha_1^{k+1}}{2\sigma_1^{k+1}} - \tfrac{(1 - \gamma_\lambda^k)\alpha_{T^k}^k}{2\sigma_{T^k}^k},s^{k+1} - s^k + v^{k+1}\}$ new inequality:  
\begin{align}
     &E[\norm{\lambda^* - \lambda_{T^K}^K}^2] + E[D_\lambda^{K+1}(\lambda^*)] \leq \tfrac{1}{q^K}\big(\tfrac{\alpha_1^1\theta_1^1}{2\eta_1^1}E[\norm{\lambda_0^1 -\lambda_{-1}^1}^2] + \tsum_{k=1}^{K}v^k\norm{\lambda^*}^2 + \tfrac{(1 - \gamma_\lambda^1)\alpha_1^1}{2\sigma_1^1} E[\norm{\lambda^* - \lambda_0^1}^2]\nonumber\\ 
     & + (s^1 + v^1)E[D_\lambda^1(\lambda^*)] + \tsum_{k=1}^{K-1}p^k(E[\norm{\lambda^* - \lambda_{T^k}^k}^2] + E[D_\lambda^{k+1}(\lambda^*)]) \nonumber\\
     & + (\tfrac{1-\gamma_x^1}{2\tau_1^1} - \tfrac{\mu_f}{2})\alpha_1^1E[\norm{x^* - x_0^1}^2] + r^1E[D_x^1(x^*)] + (\tfrac{M_g^2\alpha_1^1\theta_1^1\eta_1^1}{2} + \tfrac{\alpha_1^1\theta_1^1 M_g}{2\rho_1^1})E[\norm{x_0^1 - x_{-1}^1}^2]\big). \label{eq:a86}
\end{align}
Using Lemma \ref{eq:a35} for \eqref{eq:a86}, $E[\norm{\lambda^* - \lambda_{T^k}^k}^2] + E[D_\lambda^{k+1}(\lambda^*)] \leq \tfrac{R}{q^K}(1 + \tfrac{\max_{k\in [K]}p^k}{q^K})^K$. $(1 + \tfrac{\max_{k\in [K]}p^k}{q^K})^K$ will converge to $1$. Taking the worst case for $(1 + \tfrac{\max_{k\in [K]}p^k}{q^K})^K$, we have
\begin{align}
    &E[\norm{\lambda^* - \lambda_{T^K}^K}^2] + E[D_\lambda^{K+1}(\lambda^*)] \leq \sup_{K \geq 1} \tfrac{R}{q^K}(1 + \tfrac{\max_{k\in [K]}p^k}{q^K})^K \leq \sup_{K \geq 1} \tfrac{2R}{q^K} := D_{\lambda}^2,\label{eq:a88}
\end{align}
\end{proof}
\begin{lemma} 
    The new expected gap inequality, via using the upper bound \eqref{eq:b35}, is as follows:
\begin{align}
     &\tsum_{k=1}^{K}u^kE[Q(\wt{z}^k,z^*)] \leq \alpha_1^1\theta_1^1E[\inprod{\bar{q}_0^1}{\lambda_0^1 - \lambda^*}] + \tsum_{k=1}^{K-1}(\alpha_1^{k+1}\theta_1^{k+1} - \alpha_{T^k}^k)E[\inprod{\bar{q}_{T^k}^k}{\lambda_{T^k}^k - \lambda^*}] \nonumber\\
     & + \tsum_{k=1}^{K}\tsum_{t=1}^{T^k-1}(\alpha_{t+1}^k\theta_{t+1}^k - \alpha_t^k)E[\inprod{\bar{q}_t^k}{\lambda_t^k - \lambda^*}] + \tfrac{\alpha_1^1\theta_1^1}{2\eta_1^1}E[\norm{\lambda_0^1 -\lambda_{-1}^1}^2] \nonumber \\
     & + \tsum_{k=1}^{K-1}\big(\tfrac{\alpha_1^{k+1}\theta_1^{k+1}}{2\eta_1^{k+1}} - \tfrac{\alpha_{T^k}^k\theta_{T^k}^k}{2\eta_{T^k}^k}\big)E[\norm{\lambda_{T^k}^k - \lambda_{T^k-1}^k}^2] - \tfrac{\alpha_{T^K}^K\theta_{T^K}^K}{2\eta_{T^K}^K}E[\norm{\lambda_{T^K}^K - \lambda_{T^K-1}^K}^2]\nonumber\\ 
     & + \tsum_{k=1}^{K}\tsum_{t=0}^{T^k-1}(\tfrac{\theta_{t+2}^k\alpha_{t+2}^k}{2\eta_{t+2}^k} + \tfrac{\theta_{t+1}^k\alpha_{t+1}^k\rho_{t+1}^k M_g}{4} + \tfrac{\alpha_{t+1}^k}{2\omega_{t+1}^k} - \tfrac{(1 - \gamma_\lambda^k)\alpha_{t+1}^k}{2\sigma_{t+1}^k})E[\norm{\lambda_{t+1}^k - \lambda_t^k}^2] \nonumber\\ 
     & + \tfrac{(1 - \gamma_\lambda^1)\alpha_1^1}{2\sigma_1^1} E[\norm{\lambda^* - \lambda_0^1}^2] + \tsum_{k=1}^{K-1}(\tfrac{(1 - \gamma_\lambda^{k+1})\alpha_1^{k+1}}{2\sigma_1^{k+1}} - \tfrac{(1 - \gamma_\lambda^k)\alpha_{T^k}^k}{2\sigma_{T^k}^k})E[\norm{\lambda^* - \lambda_{T^k}^k}^2]  \nonumber\\
     & - \tfrac{(1 - \gamma_\lambda^K)\alpha_{T^K}^K}{4\sigma_{T^K}^K}E[\norm{\lambda^* - \lambda_{T^K}^K}^2] + \tsum_{k=1}^{K}\tfrac{1 - \gamma_\lambda^k}{2}\tsum_{t=1}^{T^k-1}(\tfrac{\alpha_{t+1}^k}{\sigma_{t+1}^k} - \tfrac{\alpha_t^k}{\sigma_t^k})E[\norm{\lambda^* - \lambda_t^k}^2] \nonumber \\
     & + s^1E[D_\lambda^1(\lambda^*)] + \tsum_{k=1}^{K-1}(s^{k+1} - s^k)E[D_\lambda^{k+1}(\lambda^*)] - s^KE[D_\lambda^{K+1}(\lambda^*)] \nonumber\\
     & + \tsum_{k=1}^{K}\tsum_{t=0}^{T^k-1}(M_g^2\nu_{t+1}^k - \tfrac{\gamma_\lambda^k}{2\sigma_{t+1}^k})\alpha_{t+1}^k E[D_\lambda^k(\lambda_{t+1}^k)] + \tsum_{k=1}^{K}\tsum_{t=0}^{T^k-1}(\tfrac{\norm{\lambda^*}L_g}{2} + \tfrac{L_f}{2} + \tfrac{M_g^2\alpha_{t+2}^k\theta_{t+2}^k\eta_{t+2}^k}{2\alpha_{t+1}^k} \nonumber\\
     & + M_g^2\omega_{t+1}^k+ \tfrac{\alpha_{t+2}^k\theta_{t+2}^k M_g}{4\rho_{t+2}^k\alpha_{t+1}^k} + \tfrac{1}{2\nu_{t+1}^k} + \nu_{t+1}^k L_f^2 + \tfrac{M_g^2}{\nu_{t+1}^k} + 4L_fD_xD_{\lambda}M - \tfrac{1-\gamma_x^k}{2\tau_{t+1}^k})\alpha_{t+1}^kE[\norm{x_t^k - x_{t+1}^k}^2] \nonumber\\
     & + (\tfrac{1-\gamma_x^1}{2\tau_1^1} - \tfrac{\mu_f}{2})\alpha_1^1E[\norm{x^* - x_0^1}^2] + \tsum_{k=1}^{K-1}\big((\tfrac{(1-\gamma_x^{k+1})}{2\tau_1^{k+1}} - \tfrac
     {\mu_f}{2})\alpha_1^{k+1}  - \tfrac{(1-\gamma_x^k)\alpha_{T^k}^k}{2\tau_{T^k}^k}\big)E[\norm{x^* - x_{T^k}^k}^2]    \nonumber \\
     & - \tfrac{(1-\gamma_x^K)\alpha_{T^K}^K}{2\tau_{T^K}^K}E[\norm{x^* - x_{T^K}^K}^2] + \tsum_{k=1}^{K}\tsum_{t=1}^{T^k-1}\big((\tfrac{1-\gamma_x^k}{2\tau_{t+1}^k} - \tfrac{\mu_f}{2})\alpha_{t+1}^k - \tfrac{(1-\gamma_x^k)\alpha_t^k}{2\tau_t^k}\big)E[\norm{x^* - x_t^k}^2] \nonumber\\
     & + r^1E[D_x^1(x^*)] + \tsum_{k=1}^{K-1}(r^{k+1} - r^k)E[D_x^{k+1}(x^*)] - r^KE[D_x^{K+1}(x^*)]  \nonumber\\
     & + \tsum_{k=1}^{K}\tsum_{t=0}^{T^k-1}(M_g^2\omega_{t+1}^k + L_f^2\nu_{t+1}^k - \tfrac{\gamma_x^k}{2\tau_{t+1}^k})\alpha_{t+1}^kE[D_x^k(x_{t+1}^k)]  + (\tfrac{M_g^2\alpha_1^1\theta_1^1\eta_1^1}{2} + \tfrac{\alpha_1^1\theta_1^1 M_g}{4\rho_1^1})E[\norm{x_0^1 - x_{-1}^1}^2]\nonumber\\
     & + (\tfrac{M_g^2\sigma_{T^K}^K\alpha_{T^K}^K\theta_{T^K}^K}{1-\gamma_{\lambda}^K} - \tfrac{M_g^2\alpha_{T^K}^K\theta_{T^K}^K\eta_{T^K}^K}{2} - \tfrac{\alpha_{T^K}^K\theta_{T^K}^K M_g}{2\rho_{T^K}^K})E[\norm{x_{T^K}^K - x_{T^K-1}^K}^2] \nonumber\\
     & + \tsum_{k=1}^{K-1}\big(\tfrac{M_g^2\alpha_1^{k+1}\theta_1^{k+1}\eta_1^{k+1}}{2} + \tfrac{\alpha_1^{k+1}\theta_1^{k+1} M_g}{2\rho_1^{k+1}} - \tfrac{M_g^2\alpha_{T^k}^k\theta_{T^k}^k\eta_{T^k}^k}{2} - \tfrac{\alpha_{T^k}^k\theta_{T^k}^k M_g}{2\rho_{T^k}^k}\big)E[\norm{x_{T^k}^k - x_{T^k-1}^k}^2]. \label{eq:a200}
\end{align}
\end{lemma}
\begin{proof}
    Starting with gap inequality \eqref{eq:a58}
    \begin{align}
        &Q(z_{t+1}^k, z) \leq \inprod{\delta_t^{\lambda^k}}{\lambda_{t+1}^k - \lambda} 
        - \inprod{\bar{q}_{t+1}^k}{\lambda_{t+1}^k - \lambda} 
        + \inprod{\theta_{t+1}^k \bar{q}_t^k}{\lambda_t^k - \lambda} \nonumber \\
        & + \inprod{\theta_{t+1}^k(q_t^k - \bar{q}_t^k)}{\lambda_t^k - \lambda} 
        + \tfrac{\theta_{t+1}^k M_g}{4\rho_{t+1}^k}\norm{x_t^k - x_{t-1}^k}^2 + \tfrac{\theta_{t+1}^k\rho_{t+1}^k M_g}{4}\norm{\lambda_{t+1}^k - \lambda_t^k}^2 \nonumber\\ 
        & + \tfrac{1 - \gamma_\lambda^k}{2\sigma_{t+1}^k} \big( \norm{\lambda - \lambda_t^k}^2 - \norm{\lambda - \lambda_{t+1}^k}^2 - \norm{\lambda_t^k - \lambda_{t+1}^k}^2 \big) + \tfrac{\gamma_\lambda^k}{2\sigma_{t+1}^k} \big( D_\lambda^k(\lambda) - \norm{\lambda - \lambda_{t+1}^k}^2 
        - D_\lambda^k(\lambda_{t+1}^k) \big)\nonumber \\
        & + \tfrac{\norm{\lambda}L_g + L_f}{2}\norm{x_t^k - x_{t+1}^k}^2 + \inprod{\delta_t^{x^k}}{x - x_{t+1}^k}  \nonumber \\
        & + (\tfrac{1-\gamma_x^k}{2\tau_{t+1}^k} - \tfrac{\mu_f}{2})\norm{x - x_t^k}^2 - \tfrac{1-\gamma_x^k}{2\tau_{t+1}^k}\norm{x - x_{t+1}^k}^2 - \tfrac{1-\gamma_x^k}{2\tau_{t+1}^k}\norm{x_{t+1}^k - x_t^k}^2 \nonumber \\
        & + \tfrac{\gamma_x^k}{2\tau_{t+1}^k} \big( D_x^k(x) 
        - \norm{x - x_{t+1}^k}^2 
        - D_x^k(x_{t+1}^k) \big). \label{eq:a201}
    \end{align}
    For $E[\inprod{\delta_t^{x^k}}{x - x_{t+1}^k}]$, we use upper bound \eqref{eq:b35} and for $E[\inprod{\delta_t^{\lambda^k}}{\lambda_{t+1}^k - \lambda}]$:
    \begin{align}
        & E[\inprod{\delta_t^{\lambda^k}}{\lambda_{t+1}^k - \lambda}] = E[\inprod{\delta_t^{\lambda^k}}{\lambda_{t+1}^k - \lambda_t^k}] \leq E[\norm{\delta_t^{\lambda^k}}\norm{\lambda_{t+1}^k - \lambda_t^k}] \nonumber\\
        & \leq M_g^2\omega_{t+1}^k(E[D_x^k(x_{t+1}^k)]  + E[\norm{x_{t+1}^k - x_{t}^k}^2]) + \tfrac{1}{2\omega_{t+1}^k}E[\norm{\lambda_{t+1}^k - \lambda_t^k}^2], \label{eq:a202}
    \end{align}
    where in \eqref{eq:a202}, the equality comes from $E[\inprod{\delta_t^{\lambda^k}}{\lambda_t^k - \lambda}] = 0$, the first inequality comes from \eqref{eq:p9}, the second inequality comes from \eqref{eq:p10}, and \eqref{eq:a39} for $E[\norm{\delta_t^{\lambda^k}}^2]$. For $\inprod{\theta_{t+1}^k(q_t^k - \bar{q}_t^k)}{\lambda_t^k - \lambda}$:
    \begin{align}
        & E[\inprod{\theta_{t+1}^k(q_t^k - \bar{q}_t^k)}{\lambda_t^k - \lambda}] = E[\inprod{\theta_{t+1}^k(q_t^k - \bar{q}_t^k)}{\lambda_t^k - \lambda_{t-1}^k}] \leq E[\norm{\theta_{t+1}^k(q_t^k - \bar{q}_t^k)}\norm{\lambda_t^k - \lambda_{t-1}^k}] \nonumber\\
        & \leq  \tfrac{\eta_{t+1}^k\theta_{t+1}^kM_g^2}{2}E[\norm{x_t^k - x_{t-1}^k}^2] + \tfrac{\theta_{t+1}^k}{2\eta_{t+1}^k}E[\norm{\lambda_t^k - \lambda_{t-1}^k}^2], \label{eq:a203}
    \end{align}
    where in \eqref{eq:a203}, the equality comes from $E[\inprod{\theta_{t+1}^k(q_t^k - \bar{q}_t^k)}{\lambda_{t-1}^k - \lambda}] = 0$, the first inequality comes from \eqref{eq:p9} ,the second inequality comes from \eqref{eq:p10}, and upper bound \eqref{eq:a41} for $E[\norm{q_t^k - \bar{q}_t^k}^2]$. Putting $(x,\lambda,z)=(x^*,\lambda^*,z^*)$, taking expectation on both sides of \eqref{eq:a201}, using upper bounds \eqref{eq:a203}, \eqref{eq:a202} and \eqref{eq:b35}, and multiplying both sides of inequality by $\alpha_{t+1}^k$:
    \begin{align}
        & \alpha_{t+1}^kE[Q(z_{t+1}^k, z^*)] \leq M_g^2\omega_{t+1}^k\alpha_{t+1}^k(E[D_x^k(x_{t+1}^k)]  + E[\norm{x_{t+1}^k - x_{t}^k}^2]) \nonumber\\
        & + (\tfrac{\alpha_{t+1}^k}{2\omega_{t+1}^k} + \tfrac{\theta_{t+1}^k\rho_{t+1}^k M_g}{4})E[\norm{\lambda_{t+1}^k - \lambda_t^k}^2] - \alpha_{t+1}^kE[\inprod{\bar{q}_{t+1}^k}{\lambda_{t+1}^k - \lambda^*}] 
        + \alpha_{t+1}^kE[\inprod{\theta_{t+1}^k \bar{q}_t^k}{\lambda_t^k - \lambda^*}] \nonumber \\
        & + (\tfrac{\eta_{t+1}^k\theta_{t+1}^k\alpha_{t+1}^kM_g^2}{2} + \tfrac{\theta_{t+1}^k \alpha_{t+1}^kM_g}{4\rho_{t+1}^k})E[\norm{x_t^k - x_{t-1}^k}^2] + \tfrac{\theta_{t+1}^k\alpha_{t+1}^k}{2\eta_{t+1}^k}E[\norm{\lambda_t^k - \lambda_{t-1}^k}^2] 
          \nonumber\\ 
        & + \tfrac{(1 - \gamma_\lambda^k)\alpha_{t+1}^k}{2\sigma_{t+1}^k} \big( E[\norm{\lambda^* - \lambda_t^k}^2] - E[\norm{\lambda^* - \lambda_{t+1}^k}^2] - E[\norm{\lambda_t^k - \lambda_{t+1}^k}^2] \big) \nonumber \\
        & + \tfrac{\gamma_\lambda^k\alpha_{t+1}^k}{2\sigma_{t+1}^k} \big( E[D_\lambda^k(\lambda^*)] 
        - E[\norm{\lambda^* - \lambda_{t+1}^k}^2] 
        \big) + (M_g^2 \nu_{t+1}^k - \tfrac{\gamma_\lambda^k}{2\sigma_{t+1}^k})\alpha_{t+1}^kE[D_\lambda^k(\lambda_{t+1}^k)]\nonumber \\
        & + (\tfrac{\norm{\lambda^*}L_g + L_f + 8L_fD_xD_{\lambda}M + 2\nu_{t+1}^k L_f^2}{2} + \tfrac{1}{2\nu_{t+1}^k} + \tfrac{M_g^2}{\nu_{t+1}^k} + M_g^2\omega_{t+1}^k - \tfrac{1-\gamma_x^k}{2\tau_{t+1}^k})\alpha_{t+1}^kE[\norm{x_t^k - x_{t+1}^k}^2] \nonumber\\
        & + (\tfrac{1-\gamma_x^k}{2\tau_{t+1}^k} - \tfrac{\mu_f}{2})\alpha_{t+1}^kE[\norm{x^* - x_t^k}^2] - \tfrac{(1-\gamma_x^k)\alpha_{t+1}^k}{2\tau_{t+1}^k}E[\norm{x^* - x_{t+1}^k}^2]  \nonumber \\
        & + \tfrac{\gamma_x^k\alpha_{t+1}^k}{2\tau_{t+1}^k} \big( E[D_x^k(x^*)] 
        - E[\norm{x^* - x_{t+1}^k}^2]\big) + (M_g^2\omega_{t+1}^k + \nu_{t+1}^k L_f^2 - \tfrac{\gamma_x^k}{2\tau_{t+1}^k})\alpha_{t+1}^kE[D_x^k(x_{t+1}^k)]. \label{eq:a204}
    \end{align}
    Notice that we have replaced $\eta^k$ with $\nu_{t+1}^k$ for \eqref{eq:b35}. Summing up \eqref{eq:a204} with respect to $t$ and using \eqref{eq:p13}
    \begin{align}
     &u^kE[Q(\wt{z}^k,z^*)] \leq  \alpha_1^k\theta_1^kE[\inprod{\bar{q}_0^k}{\lambda_0^k - \lambda^*}] + \tsum_{t=1}^{T^k-1}(\alpha_{t+1}^k\theta_{t+1}^k - \alpha_t^k)E[\inprod{\bar{q}_t^k}{\lambda_t^k - \lambda^*}] \nonumber \\
     & - \alpha_{T^k}^kE[\inprod{\bar{q}_{T^k}^k}{\lambda_{T^k}^k - \lambda^*}] + \tfrac{\alpha_{1}^k\theta_1^k}{2\eta_1^k}E[\norm{\lambda_0^k - \lambda_{-1}^k}^2] - \tfrac{\alpha_{T^k}^k\theta_{T^k}^k}{2\eta_{T^k}^k}E[\norm{\lambda_{T^k}^k - \lambda_{T^k-1}^k}^2] \nonumber\\ 
     & + \tsum_{t=0}^{T^k-2}(\tfrac{\theta_{t+2}^k\alpha_{t+2}^k}{2\eta_{t+2}^k} + \tfrac{\theta_{t+1}^k\alpha_{t+1}^k\rho_{t+1}^k M_g}{4} + \tfrac{\alpha_{t+1}^k}{2\omega_{t+1}^k} - \tfrac{(1 - \gamma_\lambda^k)\alpha_{t+1}^k}{2\sigma_{t+1}^k})E[\norm{\lambda_{t+1}^k - \lambda_t^k}^2] \nonumber\\ 
     & + \tfrac{(1 - \gamma_{\lambda}^k)\alpha_1^k}{2\sigma_1^k}E[\norm{\lambda^* - \lambda_0^k}^2] + \tfrac{1 - \gamma_{\lambda}^k}{2}\tsum_{t=1}^{T^k-1}(\tfrac{\alpha_{t+1}^k}{\sigma_{t+1}^k} - \tfrac{\alpha_t^k}{\sigma_t^k})E[\norm{\lambda^* - \lambda_t^k}^2] - \tfrac{(1 - \gamma_{\lambda}^k)\alpha_{T^k}^k}{2\sigma_{T^k}^k}E[\norm{\lambda^* - \lambda_{T^k}^k}^2] \nonumber \\
     & + s^kE[D_\lambda^k(\lambda^*)] - s^kE[D_\lambda^{k+1}(\lambda^*)] + \tsum_{t=0}^{T^k-1}(M_g^2\nu_{t+1}^k - \tfrac{\gamma_\lambda^k}{2\sigma_{t+1}^k})\alpha_{t+1}^kE[D_\lambda^k(\lambda_{t+1}^k)] \nonumber \\
     & + \tsum_{t=0}^{T^k-1}\big(\tfrac{\norm{\lambda^*}L_g + L_f + 8L_fD_xD_{\lambda}M + 2\nu_{t+1}^k L_f^2}{2} + \tfrac{M_g^2\alpha_{t+2}^k\theta_{t+2}^k\eta_{t+2}^k}{2\alpha_{t+1}^k} + \tfrac{\alpha_{t+2}^k\theta_{t+2}^k M_g}{4\rho_{t+2}^k\alpha_{t+1}^k} + \tfrac{1}{2\nu_{t+1}^k} \nonumber\\
     & + M_g^2\omega_{t+1}^k + \tfrac{M_g^2}{\nu_{t+1}^k} - \tfrac{1-\gamma_x^k}{2\tau_{t+1}^k}\big)\alpha_{t+1}^kE[\norm{x_t^k - x_{t+1}^k}^2] + (\tfrac{1-\gamma_x^k}{2\tau_1^k} - \tfrac{\mu_f}{2})\alpha_1^kE[\norm{x^* - x_0^k}^2] \nonumber\\
     &  + \tsum_{t=1}^{T^k-1}\big((\tfrac{1-\gamma_x^k}{2\tau_{t+1}^k} - \tfrac{\mu_f}{2})\alpha_{t+1}^k - \tfrac{(1-\gamma_x^k)\alpha_t^k}{2\tau_t^k}\big)E[\norm{x^* - x_t^k}^2] - \tfrac{(1-\gamma_x^k)\alpha_{T^k}^k}{2\tau_{T^k}^k}E[\norm{x^* - x_{T^k}^k}^2]\nonumber\\
     & + r^kE[D_x^k(x^*)] - r^kE[D_x^{k+1}(x^*)] + \tsum_{t=0}^{T^k-1}(M_g^2\omega_{t+1}^k + L_f^2\nu_{t+1}^k - \tfrac{\gamma_x^k}{2\tau_{t+1}^k})\alpha_{t+1}^kE[D_x^k(x_{t+1}^k)]  \nonumber\\
     & + (\tfrac{M_g^2\alpha_1^k\theta_1^k\eta_1^k}{2} + \tfrac{\alpha_1^k\theta_1^k M_g}{4\rho_1^k})E[\norm{x_0^k - x_{-1}^k}^2] - (\tfrac{M_g^2\alpha_{T^k}^k\theta_{T^k}^k\eta_{T^k}^k}{2} + \tfrac{\alpha_{T^k}^k\theta_{T^k}^k M_g}{4\rho_{T^k}^k})E[\norm{x_{T^k}^k - x_{T^k-1}^k}^2]. 
\end{align}
Summing up with respect to index $k$, and applying \eqref{eq:p15} will give \eqref{eq:a200}.
\end{proof}
\begin{remark} \label{eq:l5}
    The following conditions are necessary to get the convergence for Algorithm \ref{alg:alg3} using the gap inequality \eqref{eq:a200}.
    \begin{align*}
        \alpha_1^{k+1}\theta_1^{k+1} - \alpha_{T^k}^k &= 0 & \alpha_{t+1}^k\theta_{t+1}^k - \alpha_t^k &= 0 & \tfrac{\alpha_1^{k+1}\theta_1^{k+1}}{2\eta_1^{k+1}} - \tfrac{\alpha_{T^k}^k\theta_{T^k}^k}{2\eta_{T^k}^k} & \leq 0 & M_g\nu_{t+1}^k - \tfrac{\gamma_\lambda^k}{2\sigma_{t+1}^k} \leq 0,
    \end{align*}
    \begin{align*}
        \tfrac{\theta_{t+1}^k\rho_{t+1}^k M_g}{2} + \tfrac{1}{2\omega_{t+1}^k} + \tfrac{\alpha_{t+2}^k\theta_{t+2}^k}{2\eta_{t+2}^k\alpha_{t+1}^k} - \tfrac{1 - \gamma_\lambda^k}{2\sigma_{t+1}^k} & \leq 0 & \tfrac{\alpha_{t+1}^k}{\sigma_{t+1}^k} - \tfrac{\alpha_t^k}{\sigma_t^k} & \leq 0 , 
    \end{align*}
    \begin{align*}
        \tfrac{\norm{\lambda^*}L_g + L_f + 8L_fD_xD_{\lambda}M + 2\nu_{t+1}^k L_f^2}{2} + \tfrac{M_g^2\alpha_{t+2}^k\theta_{t+2}^k\eta_{t+2}^k}{2\alpha_{t+1}^k} + M_g^2\omega_{t+1}^k+ \tfrac{\alpha_{t+2}^k\theta_{t+2}^k M_g}{2\rho_{t+2}^k\alpha_{t+1}^k} + \tfrac{1 + 2M_g^2}{2\nu_{t+1}^k} - \tfrac{1-\gamma_x^k}{2\tau_{t+1}^k} = m_{t+1}^k \leq 0, 
    \end{align*}
    \begin{align*}
        (\tfrac{(1-\gamma_x^{k+1})}{2\tau_1^{k+1}} - \tfrac
     {\mu_f}{2})\alpha_1^{k+1}  - \tfrac{(1-\gamma_x^k)\alpha_{T^k}^k}{2\tau_{T^k}^k} & \leq 0 & (\tfrac{1-\gamma_x^k}{2\tau_{t+1}^k} - \tfrac{\mu_f}{2})\alpha_{t+1}^k - \tfrac{(1-\gamma_x^k)\alpha_t^k}{2\tau_t^k} & \leq 0 & M_g^2\omega_{t+1}^k + L_f^2\nu_{t+1}^k - \tfrac{\gamma_x^k}{2\tau_{t+1}^k} & \leq 0.
    \end{align*}
    \begin{align*}
         \tfrac{M_g^2\sigma_{T^K}^K\alpha_{T^K}^K\theta_{T^K}^K}{1-\gamma_{\lambda}^K} - \tfrac{M_g^2\alpha_{T^K}^K\theta_{T^K}^K\eta_{T^K}^K}{2} - \tfrac{\alpha_{T^K}^K\theta_{T^K}^K M_g}{2\rho_{T^K}^K}  \leq 0.
    \end{align*}
\end{remark}
\begin{theorem} \label{eq:a206}
    Assume the following setup for the parameters holds:
    \begin{align*}
        A_t^k &=t_0 + \tfrac{k(k-1)}{2} + t &\alpha_t^k &= \tfrac{w}{\mu_f}A_t^k & \theta_t^k &= \tfrac{\alpha_{t-1}^k}{\alpha_t^k} & \rho_t^k &= \tfrac{4M_g w}{\mu_f(A_t^k - t_0)}   &\sigma_t^k &= \tfrac{\mu_f(A_t^k - t_0)}{32C_{\nu}M_gwk}, 
    \end{align*}
    \begin{align*}
         \gamma_{\lambda}^k &= \tfrac{1}{k+1} &\nu_t^k &= \tfrac{16C_{\nu}w}{\mu_f(A_t^k - t_0)} & \tau_t^k &= \tfrac{w}{2\mu_fA_t^k} & \gamma_x^k &= \tfrac{1}{k+1} & T^k = k,
    \end{align*}
    \begin{align*}
        \eta_t^k & = \tfrac{\mu_f(A_t^k - t_0)}{8M_g^2w} & \omega_t^k &= \tfrac{\mu_f(A_t^k - t_0)}{2M_g^2wk} ,
    \end{align*}
    where $t_0 = \tfrac{w(\norm{\lambda^*}L_g + L_f + 8L_fD_xD_{\lambda}M + 2L_f^2)}{2\mu_f}$, $C_{\nu} = M_g + 2L_gD_x + L_f$ and $w \geq 9$. The gradients computation for Algorithm \ref{alg:alg3} using the gap inequality \eqref{eq:a200} to reach $\epsilon$-gap for objective function and function constraint will be $O(\tfrac{m}{\epsilon^{1/3}} + \tfrac{1}{\epsilon^{2/3}})$ and $O(\tfrac{n}{\epsilon^{1/3}} + \tfrac{1}{\epsilon^{2/3}})$ respectively. 
\end{theorem}
\begin{proof}
    The setup in Theorem \ref{eq:a206} satisfies the conditions in Remark \ref{eq:l5}. Using parameters in Theorem \ref{eq:a206}, new gap inequality
\begin{align}
     &\tsum_{k=1}^{K}u^kE[Q(\wt{z}^k,z^*)] \leq  \alpha_1^1\theta_1^1E[\inprod{\bar{q}_0^1}{\lambda_0^1 - \lambda^*}] + \tfrac{\alpha_1^1\theta_1^1}{2\eta_1^1}E[\norm{\lambda_0^1 -\lambda_{-1}^1}^2] - \tfrac{\alpha_{T^K}^K\theta_{T^K}^K}{2\eta_{T^K}^K}E[\norm{\lambda_{T^K}^K - \lambda_{T^K-1}^K}^2]\nonumber\\ 
     & + \tfrac{(1 - \gamma_\lambda^1)\alpha_1^1}{2\sigma_1^1} E[\norm{\lambda^* - \lambda_0^1}^2] + \tsum_{k=1}^{K-1}(\tfrac{(1 - \gamma_\lambda^{k+1})\alpha_1^{k+1}}{2\sigma_1^{k+1}} - \tfrac{(1 - \gamma_\lambda^k)\alpha_{T^k}^k}{2\sigma_{T^k}^k})E[\norm{\lambda^* - \lambda_{T^k}^k}^2]   \nonumber\\
     & - \tfrac{(1 - \gamma_\lambda^K)\alpha_{T^K}^K}{4\sigma_{T^K}^K}E[\norm{\lambda^* - \lambda_{T^K}^K}^2] + s^1E[D_\lambda^1(\lambda^*)] + \tsum_{k=1}^{K-1}(s^{k+1} - s^k)E[D_\lambda^{k+1}(\lambda^*)] - s^KE[D_\lambda^{K+1}(\lambda^*)] \nonumber\\
     & + (\tfrac{1-\gamma_x^1}{2\tau_1^1} - \tfrac{\mu_f}{2})\alpha_1^1E[\norm{x^* - x_0^1}^2] + \tsum_{k=1}^{K-1}\big((\tfrac{(1-\gamma_x^{k+1})}{2\tau_1^{k+1}} - \tfrac
     {\mu_f}{2})\alpha_1^{k+1}  - \tfrac{(1-\gamma_x^k)\alpha_{T^k}^k}{2\tau_{T^k}^k}\big)E[\norm{x^* - x_{T^k}^k}^2]    \nonumber \\
     & - \tfrac{(1-\gamma_x^K)\alpha_{T^K}^K}{2\tau_{T^K}^K}E[\norm{x^* - x_{T^K}^K}^2] + \tsum_{k=1}^{K}\tsum_{t=1}^{T^k-1}\big((\tfrac{1-\gamma_x^k}{2\tau_{t+1}^k} - \tfrac{\mu_f}{2})\alpha_{t+1}^k - \tfrac{(1-\gamma_x^k)\alpha_t^k}{2\tau_t^k}\big)E[\norm{x^* - x_t^k}^2] + r^1E[D_x^1(x^*)] \nonumber\\
     & + \tsum_{k=1}^{K-1}(r^{k+1} - r^k)E[D_x^{k+1}(x^*)] - r^KE[D_x^{K+1}(x^*)] + \tsum_{k=1}^{K}\tsum_{t=0}^{T^k-1}m_{t+1}^k\alpha_{t+1}^kE[\norm{x_t^k - x_{t+1}^k}^2]   \nonumber\\
     & + (\tfrac{M_g^2\alpha_1^1\theta_1^1\eta_1^1}{2} + \tfrac{\alpha_1^1\theta_1^1 M_g}{2\rho_1^1})E[\norm{x_0^1 - x_{-1}^1}^2] + (\tfrac{M_g^2\sigma_{T^K}^K\alpha_{T^K}^K\theta_{T^K}^K}{1-\gamma_{\lambda}^K} - \tfrac{M_g^2\alpha_{T^K}^K\theta_{T^K}^K\eta_{T^K}^K}{2} - \tfrac{\alpha_{T^K}^K\theta_{T^K}^K M_g}{2\rho_{T^K}^K})E[\norm{x_{T^K}^K - x_{T^K-1}^K}^2] \nonumber\\
     & + \tsum_{k=1}^{K-1}\big(\tfrac{M_g^2\alpha_1^{k+1}\theta_1^{k+1}\eta_1^{k+1}}{2} + \tfrac{\alpha_1^{k+1}\theta_1^{k+1} M_g}{2\rho_1^{k+1}} - \tfrac{M_g^2\alpha_{T^k}^k\theta_{T^k}^k\eta_{T^k}^k}{2} - \tfrac{\alpha_{T^k}^k\theta_{T^k}^k M_g}{2\rho_{T^k}^k}\big)E[\norm{x_{T^k}^k - x_{T^k-1}^k}^2]. \label{eq:a207} 
\end{align}
The cumulative coefficient of $E[\norm{x^* - x_t^k}^2]$:
\begin{align*}
    &\tfrac{r^{k+1} - r^k}{r^{k+1}}\tfrac{\gamma_x^k\alpha_t^k}{2\tau_t^k} - \big[(\tfrac{1-\gamma_x^k}{2\tau_{t+1}^k} - \tfrac{\mu_f}{2})\alpha_{t+1}^k - \tfrac{(1-\gamma_x^k)\alpha_t^k}{2\tau_t^k}\big] \nonumber\\
    & = \tfrac{2}{k^2}(t_0 + \tfrac{k(k-1)}{2} + t)^2 + \tfrac{k^2-1}{2k^2}(2t_0 + k(k-1) + 2t +1) - \tfrac{w}{2}(t_0 + \tfrac{k(k-1)}{2} + t +1) \leq 0, 
\end{align*}
when $w \geq 9$. $m_{t+1}^k\alpha_{t+1}^k$ can eliminate $\tfrac{\alpha_{T^K}^K\sigma_{T^K}^K}{1-\gamma_{\lambda}^K}M_g^2$. Also 
\begin{align*}
    \alpha_1^1\theta_1^1E[\inprod{\bar{q}_0^1}{\lambda_0^1 - \lambda^*}]  \leq \alpha_1^1\theta_1^1M_gE[\norm{x_0^1 - x_{-1}^1}]\norm{\lambda^* -\lambda_0^1}] = 0,
\end{align*}
since $x_0^1 = x_{-1}^1$. $\tfrac{M_g^2\alpha_1^{k+1}\theta_1^{k+1}\eta_1^{k+1}}{2} + \tfrac{\alpha_1^{k+1}\theta_1^{k+1} M_g}{4\rho_1^{k+1}} - \tfrac{M_g^2\alpha_{T^k}^k\theta_{T^k}^k\eta_{T^k}^k}{2} - \tfrac{\alpha_{T^k}^k\theta_{T^k}^k M_g}{4\rho_{T^k}^k}$ could be eliminated with $m_{T^k}^k\alpha_{T^k}^k$. New gap inequality
\begin{align}
     &\tsum_{k=1}^{K}u^kE[Q(\wt{z}^k,z^*)] \leq \tfrac{\alpha_1^1\theta_1^1}{2\eta_1^1}E[\norm{\lambda_0^1 -\lambda_{-1}^1}^2] + \tfrac{(1 - \gamma_\lambda^1)\alpha_1^1}{2\sigma_1^1} E[\norm{\lambda^* - \lambda_0^1}^2]  \nonumber\\
     & + \tsum_{k=1}^{K-1}(\tfrac{(1 - \gamma_\lambda^{k+1})\alpha_1^{k+1}}{2\sigma_1^{k+1}} - \tfrac{(1 - \gamma_\lambda^k)\alpha_{T^k}^k}{2\sigma_{T^k}^k})E[\norm{\lambda^* - \lambda_{T^k}^k}^2] + s^1E[D_\lambda^1(\lambda^*)] + \tsum_{k=1}^{K-1}(s^{k+1} - s^k)E[D_\lambda^{k+1}(\lambda^*)] \nonumber\\
     & + (\tfrac{1-\gamma_x^1}{2\tau_1^1} - \tfrac{\mu_f}{2})\alpha_1^1E[\norm{x^* - x_0^1}^2] + r^1E[D_x^1(x^*)] + (\tfrac{M_g^2\alpha_1^1\theta_1^1\eta_1^1}{2} + \tfrac{\alpha_1^1\theta_1^1 M_g}{4\rho_1^1})E[\norm{x_0^1 - x_{-1}^1}^2]. \label{eq:a209} 
\end{align}
Applying \eqref{eq:p13}, and dividing both sides of \eqref{eq:a209} by $U^K$, and using the constant upper bound \eqref{eq:a88}, the final gap inequality will be:
\begin{align}
     E[Q(\wt{Z}^K,z^*)] \leq \tfrac{1}{U^K}\big((4M_g^2 + 2M_g + 2 + t_0)KD_{\lambda}^2 + t_0^2D_x^2\big) = \epsilon \Rightarrow K = O(\tfrac{1}{\epsilon^{1/3}}). \label{eq:a210} 
\end{align}
The gradient computations to reach $\epsilon$-gap for the objective function and function constraint will be $Km + \tsum_{k=1}^{K}T^k = Km + \tfrac{K(K+1)}{2} = O(\tfrac{m}{\epsilon^{1/3}} + \tfrac{1}{\epsilon^{2/3}})$ and $Kn + \tsum_{k=1}^{K}T^k = Kn + \tfrac{K(K+1)}{2} = O(\tfrac{n}{\epsilon^{1/3}} + \tfrac{1}{\epsilon^{2/3}})$ respectively. This concludes the proof of Theorem \ref{eq:a206}. 
\end{proof}
\section{Numerical results} \label{sec-num}

The test problem is a Neyman--Pearson classification problem with the smooth hinge-loss function $\phi$. The loss function is defined as follows: for $z\leq 0$, we use $\phi(z)=\frac{1}{2}-z$; for $0<z\leq 1$, we use $\phi(z)=\frac{1}{2}(1-z)^2$; and for $z>1$, we use $\phi(z)=0$. The dataset used in the experiments is the rcv1 training dataset from the LIBSVM library\footnote{\url{http://www.csie.ntu.edu.tw/~cjlin/libsvmtools/datasets/binary.html}}. It contains $n=20,242$ data points with dimension $d=47,236$. Among them, $n_0=10,491$ positive data points are used in the objective function and $n_1=9,751$ negative data points are used in the constraint. We compare Algorithm \ref{alg:alg2}, denoted by $\vrconex$, with ConEx \cite{boob2023stochastic} and the stochastic primal-dual queue method of \cite{yu2017online}, denoted by Yu. The Neyman--Pearson problem has the following form:
\begin{align*}
    & \min_{\|x\|\leq r} f(x):= \tfrac{1}{n_0}\tsum_{j=1}^{n_0}\phi(x^TD_{0j}) \\
    & \text{s.t } g(x):=\tfrac{1}{n_1}\tsum_{j=1}^{n_1}\phi(-x^TD_{1j})-s\leq 0.
\end{align*}
Here, $x$ is the decision variable, $D_{0j}$ and $D_{1j}$ are data vectors, and we set $r=5$ and $s=0.1$.

For the repeated comparison experiment, each algorithm was executed
independently for 10 stochastic runs. A fresh master random seed
was generated from operating-system entropy at the beginning of each
execution. Independent random streams were then assigned to
\(\vrconex\), ConEx, Yu, and the generation of the initial point in
each run.

The primal initial point was randomized across runs. More precisely, for run \(r\), a vector \(u_r\sim\mathcal{N}(0,I_d)\) was generated and the initial point was defined as $x_0^{(r)} = \frac{u_r}{\|u_r\|}$, hence $\|x_0^{(r)}\|=1<r=5$ and every initial point was feasible. To ensure a fair comparison, all three algorithms in run \(r\) were initialized from the same primal point \(x_0^{(r)}\). However, a different initial point and different mini-batch sequences were used in the next run. The initial dual variables were set to $\lambda_0=0, y_0=0, Q_0=0$ for \(\vrconex\), ConEx, and Yu, respectively. No unconstrained warm-start solution was used. For all three algorithms, the mini-batch sizes were $b_{\mathrm{obj}}=256,b_{\mathrm{con}}=512$. The quantity used to avoid numerical difficulties in the logarithmic constraint-violation plots was $\epsilon_{\mathrm{plot}}=10^{-12}$. For \(\vrconex\), we used $T_0=300,T^k=T_0$
together with $\eta_{\tau}=0.60,\eta_{\sigma}=3.50,\gamma_0=0.18$. The primal and dual averaging parameters were $\gamma_x^k=\max\{0.005,\frac{\gamma_0}{k+1}\},\gamma_{\lambda}^k=\max\{0.005,\frac{\gamma_0}{k+1}\}$.
The primal step size was $\tau_k =\frac{\tau_0}{\sqrt{k+1}},\tau_0 = \frac{\eta_{\tau}}{L_f+10^{-8}}$. The dual step size was $\sigma_k = \frac{0.90\sigma_0}{\sqrt{k+1}}, k<50$ and $\sigma_k =\frac{1.50\sigma_0}{\sqrt{k+1}}, k\geq 50$, where $\sigma_0 = \frac{\eta_{\sigma}}{L_g+10^{-8}}$.
The computational budget was \(10{,}000\) data passes.
Since one epoch of \(\vrconex\) has gradient-operation cost $(n_0+n_1)+T_0( b_{\mathrm{obj}}+2b_{\mathrm{con}})$, the number of data passes per epoch was $\frac{(n_0+n_1) + T_0( b_{\mathrm{obj}}+2b_{\mathrm{con}})}{n_0+n_1}=19.970457$. Consequently, the number of epochs was selected as $K=\lceil\frac{10{,}000}{19.970457}\rceil=501$. The resulting final data-pass count of \(\vrconex\) was approximately 10005.199190. For ConEx, the parameters were $B=5.0,\eta_{\mathrm{scale}}=1.50,\tau_{\mathrm{scale}}=0.20,\theta=2.25,y_{\max}=100$. The corresponding primal and dual parameters were $\tau= \tau_{\mathrm{scale}}\frac{2L_g\sqrt{T}}{\sqrt{B}},\eta=\frac{\eta_{\mathrm{scale}}}{BL_f+10^{-12}}$. For Yu's method, we used $V=0.80\sqrt{T},\alpha=0.25T$. The queue update was $Q_{t+1}=\max\{Q_t+1.5(g_t(x_t)+\nabla g_t(x_t)^T(x_{t+1}-x_t)),0\}$. The effective multiplier associated with Yu's method was taken as $\lambda_t=\frac{Q_t}{V}$. For ConEx and Yu, the number of iterations was selected so that their
final data-pass budget matched the actual budget of \(\vrconex\).
Since one baseline iteration uses $b_{\mathrm{obj}}+b_{\mathrm{con}}$ sample-gradient operations, the baseline data-pass count was computed
as $\mathrm{data\ passes}=\frac{T(b_{\mathrm{obj}}+b_{\mathrm{con}})}{n_0+n_1}$. This resulted in $T=263705$ baseline iterations, corresponding to approximately 10005.208971 data passes. Therefore, the difference between the final computational budgets of \(\vrconex\) and the two baseline methods was negligible.
Let $f_{k,r}^{(a)}$, denote the objective value of algorithm \(a\) at the \(k\)-th
data-pass checkpoint in run \(r\), and let $v_{k,r}^{(a)} = \max\{g_{k,r}^{(a)},0\}$, denote the corresponding constraint violation. After verifying that all ten runs had identical data-pass checkpoints, the plotted curves were computed point-wise as $\overline{f}_k^{(a)}=\frac{1}{R}\sum_{r=1}^{R}f_{k,r}^{(a)},R=10$, and $\overline{v}_k^{(a)}=\frac{1}{R}\sum_{r=1}^{R}v_{k,r}^{(a)}$. Thus, the constraint curve in the averaged figure is the mean of the ten nonnegative constraint-violation curves, rather than $\max\{\frac{1}{R}\sum_{r=1}^{R}g_{k,r}^{(a)},0\}$.
\begin{figure}[ht]
    \centering
    \includegraphics[
        width=1\linewidth
    ]{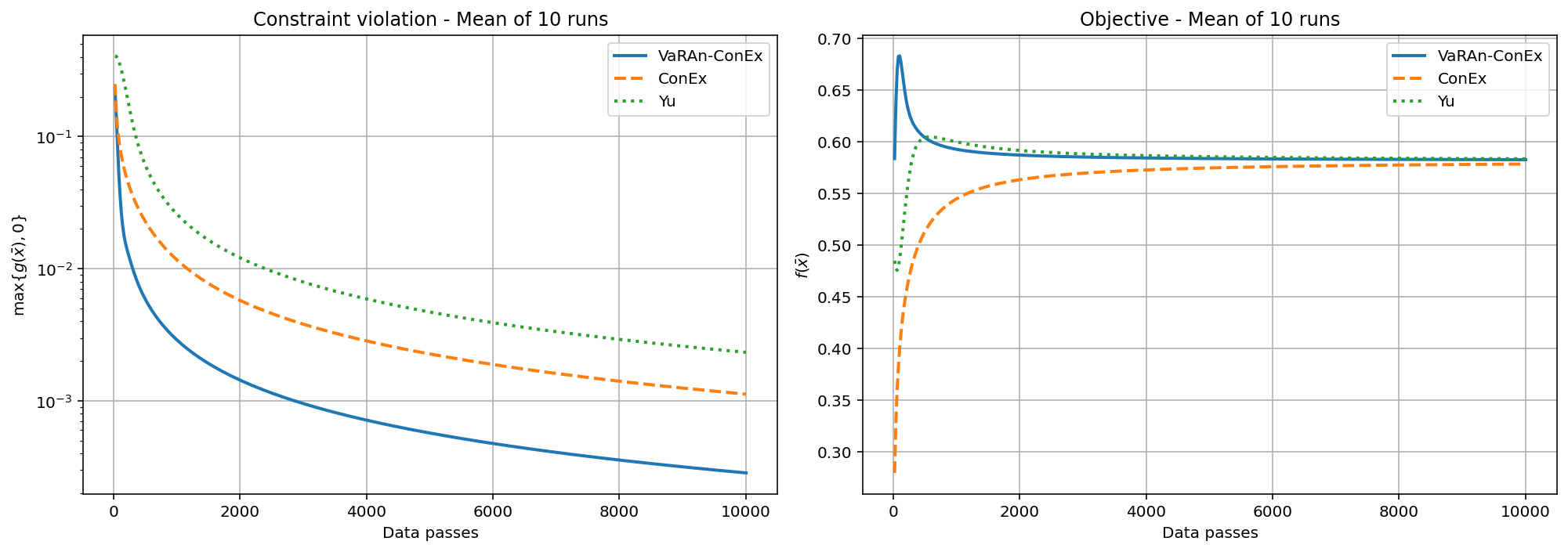}
    \caption{
        Point-wise mean objective value and constraint violation for
        \(\vrconex\), ConEx, and Yu over 10 independent stochastic
        runs. In every run, the methods use independent mini-batch
        sequences but share the same randomly generated feasible
        initial point.
    }
    \label{fig:comparison_mean_10_runs}
\end{figure}

The terminal values reported below are the arithmetic means of the
final values obtained in the ten independent runs. For
\(\vrconex\), the mean terminal values were $\wb{f}=\MeanFVR,\wb{g}=
\MeanGVR,\wb{\max\{g,0\}}=\MeanViolVR$, with mean multiplier $\wb{\lambda}=\MeanLambdaVR$. For ConEx, the corresponding mean values were $\overline{f}=\MeanFConEx,\overline{g}=\MeanGConEx,\overline{\max\{g,0\}}=\MeanViolConEx$, with $\overline{y}=\MeanYConEx$. For Yu, the mean terminal values were $\overline{f}=\MeanFYu,\overline{g}=\MeanGYu,\overline{\max\{g,0\}}=\MeanViolYu$ and $\overline{Q}=\MeanQYu,\overline{\frac{Q}{V}}=\MeanQOverVYu$. Because the gap expression is linear in \(f\), \(g\), and the
effective multiplier, this value is also equal to the arithmetic mean
of the ten run-wise raw gaps.

\begin{table}[ht]
    \centering
    \footnotesize
    \caption{
        Mean terminal numerical results over \(10\) independent
        stochastic runs under approximately the same data-pass budget.
    }
    \label{tab:mean_final_comparison}
    \begin{tabular}{|c|c|c|c|c|c|}
        \hline
        Algorithm
        &
        Mean \(f(x)\)
        &
        Mean \(\max\{g(x),0\}\)
        &
        Mean multiplier
        \\
        \hline

        \(\vrconex\)
        &
        \(\MeanFVR\)
        &
        \(\MeanViolVR\)
        &
        \(\MeanLambdaVR\)
        \\
        \hline

        ConEx
        &
        \(\MeanFConEx\)
        &
        \(\MeanViolConEx\)
        &
        \(\MeanYConEx\)
        \\
        \hline

        Yu
        &
        \(\MeanFYu\)
        &
        \(\MeanViolYu\)
        &
        \(\MeanQOverVYu\)
        \\
        \hline
    \end{tabular}
\end{table}

%% file: sec-conclude.tex
\section{Conclusion.} \label{sec-con}
This paper focused on variance-reduced first-order methods for function constrained convex optimization. We showed primal-dual frameworks based stochastic first-order methods where we anchored the proximal point update for both primal and dual variables. Moreover, we also carefully designed linearized constraint extrapolation to mitigate the impact of unbounded Lagrange multipliers on the Lipschitz smoothness constants.  We showed accelerated complexity for both constant and nonconstant policies. Moreover, for the fully-finite-sum case, we showed the dual remains bounded which allowed us to show the variance reduction and faster convergence. We also showed that similar algorithmic design yields even faster convergence when the finite-sum objective $f$ is strongly convex. 
%Finally, we demonstrated that all our methods can be applied to convex-concave saddle point problems with coupled function constraints, while maintaining the same complexity results for all methods. 
This shows the first such comprehensive complexity result for the most finite-sum function-constraints in the literature.